\documentclass[final,3p,times]{elsarticle}
\usepackage{amsmath}
\usepackage{graphicx}
\usepackage{algorithmic}
\usepackage{algorithm}
\usepackage{amssymb}
\usepackage[tight,footnotesize]{subfigure}
\usepackage{array}
\usepackage{color}
\usepackage{multirow} \usepackage{rotating}

\DeclareGraphicsExtensions{.eps}

\newcolumntype{C}[1]{>{\raggedright\let\newline\\\arraybackslash}p{#1}}

\def\x{{\mathbf x}}
\def\s{{\mathbf s}}
\def\L{{\cal L}}
\def\A{{\cal A}}
\def\CF{{\mathcal{S}}}
\def\M {{\mathbf{M}}}         
\def\OM {{\mathbf{\Omega}}}    
\def\P {{\mathbf{P}}}         
\def\Q {{\mathbf{Q}}}         
\def\x {{\mathbf{x}}}        
\def\e {{\mathbf{e}}}        
\def\y {{\mathbf{y}}}        
\def\z {{\mathbf{z}}}        
\def\mdim {m}   
\def\sdim {d}   
\def\pdim {p}   
\def\ddim {n}   
\def\spar {k}   
\def\cosp {\ell} 
\def \Dict {\mathbf{D}}         
\def \aspace {\mathcal{W}}      
\def\subjectto{\quad\text{subject to}\quad}

\def\deltaB{\delta_{2\cosp - \pdim}}
\def\deltaC{\delta_{3\cosp - 2\pdim}}
\def\deltaD{\delta_{4\cosp - 3\pdim}}

\newtheorem{thm}{Theorem}[section]
\newtheorem{cor}[thm]{Corollary}
\newtheorem{lem}[thm]{Lemma}

\newtheorem{defn}[thm]{Definition}
\newtheorem{rem}[thm]{Remark}

\newcommand{\ceil}[1]{\left\lceil#1\right\rceil}
\newcommand{\norm}[1]{\left\Vert#1\right\Vert}
\newcommand{\abs}[1]{\left\vert#1\right\vert}
\newcommand{\inn}[2]{\langle #1, #2 \rangle}
\newcommand\vect[1]{{\bf#1}}
\newcommand\matr[1]{{\bf#1}}

\newcommand\alphabf{{\boldsymbol{\alpha}}}
\newcommand{\argmin}{\operatornamewithlimits{argmin}}

\newcommand{\Real}{\mathbb R}
\newcommand\RR[1]{\mathbb{R}^{#1}}

\DeclareMathOperator{\supp}{supp}
\DeclareMathOperator{\cosupp}{cosupp}

\DeclareMathOperator{\range}{range}
\DeclareMathOperator{\rank}{rank}
\DeclareMathOperator{\spn}{span}

\journal{Special Issue in LAA on Sparse Approximate Solution of Linear Systems}

\begin{document}

\begin{frontmatter}

\title{Greedy-Like Algorithms for the Cosparse Analysis Model}


\author[csa]{R.~Giryes\corref{cor1}}
\author[inf]{S.~Nam}
\author[csa]{M.~Elad}
\author[inf]{R.Gribonval}
\author[eee]{M.~E.~Davies}

\cortext[cor1]{Corresponding author}

\address[csa]{The Department of Computer Science,
        Technion -- Israel Institute of Technology,
        Haifa 32000, Israel}
\address[inf]{INRIA Rennes - Bretagne Atlantique,
Campus de Beaulieu,
F-35042 Rennes Cedex, France}
\address[eee]{School of Engineering and Electronics,
The University of Edinburgh, \\
The King's Buildings,
Mayfield Road,
Edinburgh EH9 3JL, UK}

%




\begin{abstract}
The cosparse analysis model has been introduced recently as
an interesting alternative to the standard sparse synthesis approach.
A prominent question brought up by this new construction is the analysis pursuit problem --
the need to find a signal belonging to this model,
given a set of corrupted measurements of it.
Several pursuit methods have already been proposed based on  $\ell_1$ relaxation and a greedy approach.
In this work we pursue this question further,
and propose a new family of pursuit algorithms for the cosparse analysis model,
mimicking the greedy-like methods --
compressive sampling matching pursuit (CoSaMP), subspace pursuit (SP), iterative hard
thresholding (IHT) and hard thresholding pursuit (HTP).
Assuming the availability of a near optimal projection scheme
that finds the nearest cosparse subspace to any vector,
we provide performance guarantees for these algorithms.
Our theoretical study relies on a restricted isometry property
adapted to the context of the cosparse analysis model.
We explore empirically the performance of these algorithms
by adopting a plain thresholding projection, demonstrating their good performance.
\end{abstract}


\end{frontmatter}



\section{Introduction}
\label{sec:intro}

Many natural signals and images have been observed to be inherently
low dimensional despite their possibly very high ambient signal dimension.
It is by now well understood that this phenomenon lies at the heart of the success
of numerous methods of signal and image processing.
Sparsity-based models for signals offer an elegant and clear way to enforce such inherent low-dimensionality,
explaining their high popularity in recent years.
These models consider the signal $\x\in \RR{\sdim} $ as belonging to
a finite union of subspaces of dimension $k\ll \sdim$ \cite{Lu08Theory}.
In this paper we shall focus on one such approach --
the cosparse analysis model -- and develop pursuit methods for it.

Before we dive into the details of the model assumed and the pursuit problem,
let us first define the following generic inverse problem that will accompany us throughout the paper:
For some unknown signal $\x \in \RR{\sdim}$, an incomplete set of linear observations
$\y \in \RR{\mdim}$ (incomplete implies $\mdim<\sdim$) is available via
\begin{equation}
\label{eq:inverseProblem}
\y = \M \x + \e,
\end{equation}
where $\e \in \RR{\mdim}$ is an additive bounded noise that satisfies $\norm{\e}_2^2 \leq \epsilon^2$.
The task is to recover or approximate $\x$. In the noiseless setting where $\e = 0$, this amounts to solving $\y = \M \x$.
Of course, a simple fact in linear algebra tells us that
this problem admits infinitely many solutions  (since $\mdim<\sdim$).
Therefore, when all we
have is the observation $\y$ and the measurement/observation matrix $\M\in \RR{\mdim\times \sdim}$,
we are in a hopeless situation to recover $\x$.

\subsection{The Synthesis Approach}
This is where `sparse signal models' come into play.
In the sparse synthesis model, the signal $\x$ is assumed to have a very
sparse representation in a given fixed dictionary $\Dict \in \RR{\sdim\times\ddim}$.
In other words, there exists $\alphabf$ with few nonzero entries,
as counted by the ``$\ell_{0}$-norm'' $\norm{\alphabf}_{0}$, such that
\begin{equation}
\label{eq:synthesisModel}
\x = \Dict \alphabf, \quad \text{and} \quad \spar := \norm{\alphabf}_0 \ll \sdim.
\end{equation}
Having this knowledge we solve (\ref{eq:inverseProblem}) using 
\begin{equation}
\label{eq:synthesisL0}
\hat\x_{\ell_0} = \Dict \hat{\alphabf}_{\ell_0}, \quad \text{and} \quad \hat{\alphabf}_{\ell_0} = \argmin_{\alphabf} \norm{\alphabf}_0 \subjectto \norm{\y - \M \Dict \alphabf}_2 \le \epsilon.
\end{equation}
More details about the properties of this problem can be found in \cite{Donoho02Optimally,Gribonval03spars}.

Since solving \eqref{eq:synthesisL0} is an NP-complete problem \cite{Davis97Adaptive}, approximation techniques are required for recovering $\vect{x}$. One strategy is by using relaxation, replacing the $\ell_0$ with $\ell_1$ norm, resulting with the $\ell_1$-synthesis problem
\begin{eqnarray}
\label{eq:l1_synthesis}
\hat\x_{\ell_1} = \Dict \hat{\alphabf}_{\ell_1}, \quad \text{and} \quad \hat{\alphabf}_{\ell_1} = \argmin \norm{\alphabf}_1 & s.t. & \norm{\vect{y} - \matr{M}\matr{D}\alphabf}_2 \le \epsilon.
\end{eqnarray}
For a unitary matrix $\matr{D}$ and a vector $\vect{x}$ with $k$-sparse representation $\alphabf$, if $\delta_{2k} < \delta_{\ell_1}$ then
\begin{eqnarray}
\label{eq:l1_rec_error}
\norm{\hat{\vect{x}}_{\ell_1} - \vect{x}}_2 \le C_{\ell_1}\norm{\e}_2,
\end{eqnarray}
where $\hat{\vect{x}}_{\ell_1} = \matr{D}\hat{\alphabf}_{\ell_1}$, $\delta_{2k}$ is the constant of the restricted isometry property (RIP) of $\matr{MD}$ for 2k sparse signals,
$C_{\ell_1}$ is a constant greater than $\sqrt{2}$ and $\delta_{\ell_1}$ ($\simeq 0.4931$) is a reference constant \cite{Candes06Near, foucart10Sparse,Mo11New}.
Note that this result implies a perfect recovery in the absence of noise.
The above statement was extended also for incoherent redundant dictionaries \cite{Rauhut08Compressed}.

Another option for approximating \eqref{eq:synthesisL0} is using a greedy strategy,
like in the thresholding technique or orthogonal matching pursuit (OMP) \cite{Pati93OMP, MallatZhang93}.
A different related approach is the greedy-like family of algorithms.
Among those we have compressive sampling matching pursuit (CoSaMP) \cite{Needell09CoSaMP},
subspace pursuit (SP) \cite{Dai09Subspace}, iterative hard thresholding (IHT) \cite{Blumensath09Iterative}
and hard thresholding pursuit (HTP) \cite{Foucart11Hard}.
CoSaMP and SP were the first greedy methods shown to have guarantees in the form of
\eqref{eq:l1_rec_error} assuming $\delta_{4k}< \delta_{\text{\tiny CoSaMP}}$
and $\delta_{3k} \le \delta_{\text{\tiny SP}}$ \cite{Needell09CoSaMP,Dai09Subspace,foucart10Sparse,Giryes12RIP}.
Following their work, iterative hard thresholding (IHT)
and hard thresholding pursuit (HTP) were shown to have similar guarantees
under similar conditions \cite{Blumensath09Iterative, Foucart11Hard, Garg09Gradient, foucart10Sparse}.
Recently, a RIP based guarantee was developed also for OMP \cite{Zhang11Sparse}.

\subsection{The Cosparse Analysis Model}
Recently, a new signal model called \emph{cosparse analysis model}
was proposed in \cite{Nam12Cosparse,Nam11CosparseConf}.
The model can be summarized as follows:
For a fixed analysis operator $\OM\in\RR{\pdim\times\sdim}$ referred to as the analysis dictionary,
a signal $\x\in\RR{\sdim}$ belongs to the cosparse analysis model with cosparsity $\cosp$ if
\begin{equation}
\label{eq:analysisModel}
\cosp := \pdim - \norm{\OM\x}_0.
\end{equation}
The quantity $\cosp$ is the number of rows in $\OM$ that are orthogonal to the signal.
The signal $\x$ is said to be $\cosp$-cosparse, or simply cosparse.
We denote the indices of the zeros of the analysis representation as the \emph{cosupport} $\Lambda$ and the sub-matrix that contains the rows
from $\OM$ that belong to $\Lambda$ by $\OM_\Lambda$.
As the definition of cosparsity suggests, the emphasis of the cosparse analysis
model is on the zeros of the analysis representation vector $\OM\x$.
This contrasts the emphasis on `few non-zeros'
in the synthesis model \eqref{eq:synthesisModel}.
It is clear that in the case where
every $\cosp$ rows in $\OM$ are independent, $\x$ resides in a subspace
of dimension $\sdim - \cosp$ that consists of vectors orthogonal to the rows of $\OM_\Lambda$.
In the general case where dependencies occur between the rows of $\OM$, the dimension
is $d$ minus the rank of $\OM_\Lambda$.
This is similar to the behavior in the synthesis case where a $k$-sparse signal
lives in a $k$-dimensional space.  Thus, for this model to be effective, we assume a large value of $\ell$.

In the analysis model, recovering $\x$ from the corrupted measurements
is done by  solving the following minimization problem \cite{elad07Analysis}:
\begin{equation}
\label{eq:analysisL0}
\x_{A-\ell_0} = \argmin_{\x} \norm{\OM\x}_0 \subjectto \norm{\y - \M \x}_2 \le \epsilon.
\end{equation}
Solving this problem is NP-complete \cite{Nam12Cosparse}, just as in the synthesis case, and thus approximation methods are required.
As before, we can use an $\ell_1$ relaxation to \eqref{eq:analysisL0},
replacing the $\ell_0$ with $\ell_1$ in \eqref{eq:analysisL0}, resulting with the $\ell_1$-analysis problem \cite{Nam12Cosparse, elad07Analysis,Candes11Compressed, Vaiter12Robust}.
Another option is the greedy approach. A greedy algorithm called
Greedy Analysis Pursuit (GAP) has been developed in \cite{Nam12Cosparse, Nam11CosparseConf, Nam11GAPN}
that somehow mimics Orthogonal Matching Pursuit~\cite{Pati93OMP, MallatZhang93}
with a form of iterative reweighted least Squares (IRLS) \cite{IRLS10Deubechies}.
Other alternatives for OMP, backward greedy (BG) and orthogonal BG (OBG), were presented in \cite{Rubinstein12Cosparse} for the case that
$\matr{M}$ is the identity. For the same case, the parallel to the thresholding technique was analyzed in \cite{Peleg12Performance}.

\subsection{This Work}

Another avenue exists for the development of analysis pursuit algorithms
-- constructing methods that will imitate the family of greedy-like algorithms.
Indeed, we have recently presented preliminary and simplified versions of
analysis IHT (AIHT), analysis HTP (AHTP), analysis CoSaMP (ACoSaMP) and Analysis SP (ASP)
in \cite{Giryes11Iterative, Giryes12CoSaMP} as analysis versions of
the synthesis counterpart methods. This paper re-introduces these algorithms in a more general form,
ties them to their synthesis origins, and analyze their expected performance.
The main contribution of the paper is our result
on the stability of these analysis pursuit algorithms.
We show that after a finite number of iterations and for a given constant $c_0$,
the reconstruction result $\hat\x$ of AIHT, AHTP, ACoSaMP and ASP all satisfy
\begin{eqnarray}
\norm{\x - \hat\x}_2 \le c_0\norm{\vect{e}}_2,
\end{eqnarray}
under a RIP-like condition on $\M$ and the assumption that we are given a good near optimal projection scheme.
A bound is also given for the case where $\x$ is only nearly $\cosp$-cosparse.
Similar results for the $\ell_1$ analysis appear in \cite{Candes11Compressed,Vaiter12Robust}.
More details about the relation between these papers and our results will be given in Section~\ref{sec:guarantees}.
In addition to our theoretical results we demonstrate the performance of the four pursuit methods
under a thresholding based simple projection scheme.
Both our theoretical and empirical results show that linear dependencies in $\OM$
that result with a larger cosparsity in the signal $\x$,
lead to a better reconstruction performance.
{\em This suggests that, as opposed to the synthesis case, strong linear dependencies within $\OM$ are desired.}

This paper is organized as follows:
\begin{itemize}
\item In Section~\ref{sec:notation} we present the notation used along the paper.
\item In Section~\ref{sec:omega_RIP} we define a RIP-like property, the $\OM$-RIP, for the analysis model,
proving that it has similar characteristics like the regular RIP.
In Section~\ref{sec:near_opt_proj} the notion of near optimal projection
is proposed and some nontrivial operators for which a tractable optimal projection exists are exhibited.
Both the $\OM$-RIP and the near optimal projection are used throughout this paper as a main force for deriving our theoretical results.
\item In Section~\ref{sec:analysis_alg} the four pursuit algorithms for the cosparse analysis framework are defined,
adapted to the general format of the pursuit problem we have defined above.
\item In Section~\ref{sec:guarantees} we derive the success guarantees for all the above algorithms in a unified way.
Note that the provided results can be easily adapted to other
union-of-subspaces models given near optimal projection schemes for them,
in the same fashion done for IHT with an optimal projection scheme in \cite{Blumensath09Sampling}.
The relation between the obtained results and existing work appears in this section as well.
\item Empirical performance of these algorithms is demonstrated in Section~\ref{sec:exp}
in the context of the cosparse signal recovery problem.
We use a simple thresholding as the near optimal projection scheme in the greedy-like techniques.
\item Section~\ref{sec:conc} discuss the presented results and concludes our work.
\end{itemize}

\section{Notations and Preliminaries}
\label{sec:notation}

We use the following notation in our work:
\begin{itemize}
  \item $\sigma_{\matr{M}}$ is the largest singular value of $\matr{M}$, i.e., $\sigma_{\matr{M}}^2 = \norm{\matr{M}^*\matr{M}}_2$.
  \item $\norm{\cdot}_2$ is the euclidian norm for vectors and the spectral norm for matrices. $\norm{\cdot}_1$ is the $\ell_1$
    norm that sums the absolute values of a vector and $\norm{\cdot}_0$, though not really a norm,
    is the $\ell_0$-norm which counts the number of non-zero elements in a vector.
  \item Given a cosupport set $\Lambda$, $\matr{\Omega}_{\Lambda}$ is a sub-matrix of $\matr{\Omega}$ with the {\em rows} that belong to $\Lambda$.
  \item For given vectors $\vect{v},\vect{z}\in \RR{d}$ and an analysis dictionary $\OM$,
        $\cosupp(\OM\vect{v})$ returns the cosupport of $\OM\vect{v}$ and $\cosupp(\OM\vect{z}, \cosp)$ returns
        the index set of $\cosp$ smallest (in absolute value) elements in $\OM\vect{z}$. If more than $\ell$ elements are zero all of them are returned.
        In the case where the $\cosp$-th smallest entry is equal to the $\cosp+1$ smallest entry, one of them is chosen arbitrarily.
  \item In a similar way, in the synthesis case $\matr{D}_T$ is a sub-matrix of $\matr{D}$
  with {\em columns}\footnote{By the abuse of notation we use the same notation for the selection sub-matrices of rows and columns. The selection will be clear from the context since in analysis the focus is always on the rows and in synthesis on the columns.} corresponding to the set of indices $T$,
        $\supp(\cdot)$ returns the support of a vector,
        $\supp(\cdot, k)$ returns the set of $k$-largest elements
        and $\lceil \cdot \rceil_k$ preserves the $k$-largest elements in a vector. In the case where the $k$-th largest entry is equal to the $k+1$ largest entry, one of them is chosen arbitrarily.
  \item $\Q_\Lambda = \matr{I} - \OM_\Lambda^\dag\OM_\Lambda$ is the orthogonal projection onto the orthogonal
        complement of $\range(\matr{\Omega}_\Lambda^*)$.
  \item $\P_\Lambda = \matr{I} - \Q_\Lambda = \OM_\Lambda^\dag\OM_\Lambda$ is the orthogonal projection onto $\range(\matr{\Omega}_\Lambda^*)$.
  \item $\hat\x_{\text{\tiny AIHT}}$/$\hat\x_{\text{\tiny AHTP}}$/$\hat\x_{\text{\tiny ACoSaMP}}$/$\hat\x_{\text{\tiny ASP}}$
        are the reconstruction results of AIHT/ AHTP/ ACoSaMP/ ASP respectively. Sometimes when it is clear from the context to which algorithms we refer,
        we abuse notations and use $\hat\x$ to denote the reconstruction result.
  \item A cosupport $\Lambda$ has a corank $r$ if $\rank(\OM_\Lambda)=r$. A vector $\vect{v}$ has a corank $r$ if its cosupport has a corank $r$.
  \item $[p]$ denotes the set of integers $[1 \dots p]$.
  \item $\L_{\matr{\Omega},\cosp} = \{\Lambda \subseteq [p], \abs{\Lambda}\ge\cosp\}$ is the set of $\cosp$-cosparse cosupports and
        $\L_{\matr{\Omega},r}^{\text{corank}} = \{\Lambda \subseteq [p], \rank(\matr{\Omega_{\Lambda}})\ge r\}$ is the set of all cosupports with corresponding corank $r$.
  \item $\aspace_\Lambda =\spn^\perp(\OM_\Lambda) = \{ \Q_\Lambda\vect{z}, \vect{z}\in \RR{\sdim}\}$ is the subspace spanned by a cosparsity set $\Lambda$.
  \item $\A_{\matr{\Omega},\cosp} = \bigcup_{\Lambda \in \L_{\matr{\Omega},\cosp}}\aspace_\Lambda$ is the union of subspaces of $\cosp$-cosparse vectors and
        $\A_{\matr{\Omega},r}^{\text{corank}} = \bigcup_{\Lambda \in \L_{\matr{\Omega},r}^{\text{corank}}}\aspace_\Lambda$ is the union of subspaces of all vectors with corank $r$.
        In the case that every $\cosp$ rows of $\matr{\Omega}$ are independent it is clear that $\A_{\matr{\Omega},\cosp} = \A_{\matr{\Omega},r}^{\text{corank}}$.
        When it will be clear from the context, we will remove $\matr{\Omega}$ from the subscript.
  \item $\x \in \RR{\sdim}$  denotes the original unknown $\cosp$-cosparse vector and $\Lambda_{\x}$ its cosupport.
  \item $\vect{v},\vect{u} \in \A_\cosp$ are used to denote general $\cosp$-cosparse vectors and
        $\vect{z} \in \RR{\sdim}$ is used to denote a general vector.

\end{itemize}

\section{$\OM$-RIP Definition and its Properties}
\label{sec:omega_RIP}

We now turn to define the
$\OM$-RIP, which parallels the regular RIP as used in \cite{Candes06Near}.
This property is a very important property for the analysis of the algorithms which holds
for a large family of matrices $\matr{M}$ as we will see hereafter.

\begin{defn}
\label{def:omega_RIP}
A matrix $\matr{M}$ has the $\matr{\Omega}$-RIP property with a constant $\delta_\cosp$,
if $\delta_\cosp$ is the smallest constant that satisfies
\begin{eqnarray}
\label{eq:omega_RIP}
&& (1-\delta_{\cosp})\norm{\vect{v}}_2^2 \le \norm{\matr{M}\vect{v}}_2^2 \le (1+\delta_\cosp)\norm{\vect{v}}_2^2,
\end{eqnarray}
whenever $\OM\vect{v}$ has at least $\cosp$ zeroes.
\end{defn}

Note that though $\delta_\cosp$ is also a function of $\OM$
we abuse notation and use the same symbol for the $\OM$-RIP as the regular RIP.
It will be clear from the context to which of them we refer and what $\OM$ is in use with the $\OM$-RIP.
A similar property that looks at the corank of the vectors can be defined
\begin{defn}
\label{def:omega_RIP_corank}
A matrix $\matr{M}$ has the corank-$\matr{\Omega}$-RIP property with a constant $\delta_{r}^{\text{corank}}$,
if $\delta_{r}^{\text{corank}}$ is the smallest constant that satisfies
\begin{eqnarray}
\label{eq:omega_RIP_rank}
&& (1-\delta_{r}^{\text{corank}})\norm{{\vect{u}}}_2^2 \le \norm{\matr{M}{\vect{u}}}_2^2 \le (1+\delta_{r}^{\text{corank}})\norm{{\vect{u}}}_2^2
\end{eqnarray}
whenever the corank of $\vect{u}$ with respect to $\OM$ is greater or equal to $r$.
\end{defn}

The $\OM$-RIP, like the regular RIP, inherits several key properties.
We present only those related to $\delta_\cosp$, while very similar characteristics can be derived also for the corank-$\OM$-RIP.
The first property we pose is an immediate corollary of the $\delta_\cosp$ definition.

\begin{cor}
\label{cor:MQ_RIP_norm}
If $\matr{M}$ satisfies the $\matr{\Omega}$-RIP with a constant $\delta_{\cosp}$ then 
\begin{eqnarray}
\label{eq:MQ_RIP_norm}
&& \norm{\matr{M}\Q_{\Lambda}}_2^2 \le 1+\delta_{\cosp}
\end{eqnarray}
for any $\Lambda \in \L_{\cosp}$.
\end{cor}
{\em Proof:} Any $\vect{v} \in \A_{\cosp}$ can be represented as $\vect{v} = \Q_{\Lambda}\vect{z}$ with $\Lambda \in \L_{\cosp}$ and $\vect{z} \in \RR{\sdim}$.
Thus, the $\OM$-RIP in \eqref{eq:omega_RIP} can be reformulated as
\begin{eqnarray}
\label{eq:omega_RIP_Qz}
&& (1-\delta_{\cosp})\norm{\Q_\Lambda\vect{z}}_2^2 \le \norm{\matr{M}\Q_\Lambda\vect{z}}_2^2 \le (1+\delta_\cosp)\norm{\Q_\Lambda\vect{z}}_2^2
\end{eqnarray}
for any $\vect{z} \in \RR{\sdim}$ and $\Lambda \in \L_{\cosp}$. Since $\Q_\Lambda$ is a projection $\norm{\Q_\Lambda\vect{z}}_2^2 \le \norm{\vect{z}}_2^2$.
Combining this with the right inequality in \eqref{eq:omega_RIP_Qz} gives
\begin{eqnarray}
\label{eq:omega_RIP_Qz2}
\norm{\matr{M}\Q_\Lambda\vect{z}}_2^2 \le (1+\delta_\cosp)\norm{\vect{z}}_2^2
\end{eqnarray}
for any $\vect{z} \in \RR{\sdim}$ and $\Lambda \in \L_{\cosp}$.
The first inequality in \eqref{eq:MQ_RIP_norm} follows from \eqref{eq:omega_RIP_Qz2} by the definition of the spectral norm.
\hfill $\Box$ \bigskip

\begin{lem}
\label{lem:l_inequality}
For $\tilde{\cosp} \le \cosp$ 
it holds that $\delta_{\cosp} \le \delta_{\tilde{\cosp}}$.
\end{lem}
{\em Proof:}
Since $\A_{\cosp} \subseteq \A_{\tilde{\cosp}}$ 
the claim is immediate.
\hfill $\Box$ \bigskip

\begin{lem}
\label{lem:omega_RIP_norm}
$\matr{M}$ satisfies the $\matr{\Omega}$-RIP if and only if
\begin{eqnarray}
\label{eq:omega_RIP_norm}
&& \norm{\Q_{\Lambda}(\matr{I} - \matr{M}^*\matr{M})\Q_{\Lambda}}_2 \le \delta_\cosp
\end{eqnarray}
for any $\Lambda \in \L_{\cosp}$.
\end{lem}
{\em Proof:}
The proof is similar to the one of the regular RIP as appears in \cite{foucart10Sparse}.
 As a first step we observe that Definition~\ref{def:omega_RIP} is equivalent to requiring
\begin{eqnarray}
\abs{\norm{\matr{M}\vect{v}}_2^2  - \norm{\vect{v}}_2^2} \le \delta_{\cosp}\norm{\vect{v}}_2^2
\end{eqnarray}
for any $\vect{v} \in \A_{\cosp}$. The latter is equivalent to
\begin{eqnarray}
\abs{\norm{\matr{M}\Q_{\Lambda}{\vect{z}}}_2^2  - \norm{\Q_{\Lambda}{\vect{z}}}_2^2} \le \delta_{\cosp}\norm{\Q_{\Lambda}{\vect{z}}}_2^2
\end{eqnarray}
for any set $\Lambda \in \L_{\cosp}$ and any ${\vect{z}} \in \Real^d$, since $\Q_{\Lambda}{\vect{z}} \in \A_{\cosp}$.
Next we notice that
\begin{eqnarray*}
&& \norm{\matr{M}\Q_{\Lambda}{\vect{z}}}_2^2  - \norm{\Q_{\Lambda}{\vect{z}}}_2^2
= {\vect{z}}^*\Q_{\Lambda}\matr{M}^*\matr{M}\Q_{\Lambda}{\vect{z}} - {\vect{z}}^*\Q_{\Lambda}{\vect{z}}
= \langle\Q_{\Lambda}(\matr{M}^*\matr{M} - \matr{I})\Q_{\Lambda}{\vect{z}}, {\vect{z}} \rangle.
\end{eqnarray*}
Since $\Q_{\Lambda}(\matr{M}^*\matr{M} - \matr{I})\Q_{\Lambda}$ is Hermitian we have that
\begin{eqnarray}
&&\max_{{\vect{z}}}\frac{\abs{\langle\Q_{\Lambda}(\matr{M}^*\matr{M} - \matr{I})\Q_{\Lambda}{\vect{z}}, {\vect{z}} \rangle}}{\norm{{\vect{z}}}_2}
= \norm{\Q_{\Lambda}(\matr{M}^*\matr{M} - \matr{I})\Q_{\Lambda}}_2.
\end{eqnarray}
Thus we have that Definition~\ref{def:omega_RIP} is equivalent to
\eqref{eq:omega_RIP_norm} for any set $\Lambda \in \L_{\cosp}$.
\hfill $\Box$ \bigskip

\begin{cor}
\label{cor:omega_RIP_norm_diff}
If $\matr{M}$ satisfies the $\matr{\Omega}$-RIP then
\begin{eqnarray}
\label{eq:omega_RIP_norm_diff}
&& \norm{\Q_{\Lambda_1}(\matr{I} - \matr{M}^*\matr{M})\Q_{\Lambda_2}}_2 \le \delta_\cosp,
\end{eqnarray}
for any $\Lambda_1$ and $\Lambda_2$ such that  $\Lambda_1 \cap \Lambda_2 \in \L_{\cosp}$.
\end{cor}
{\em Proof:}
Since $\Lambda_1 \cap \Lambda_2 \subseteq \Lambda_1$ and $\Lambda_1 \cap \Lambda_2 \subseteq \Lambda_2$
\begin{eqnarray}
\nonumber
\norm{\Q_{\Lambda_1}(\matr{I} - \matr{M}^*\matr{M})\Q_{\Lambda_2}}_2 \le \norm{\Q_{\Lambda_2 \cap \Lambda_1}(\matr{I} - \matr{M}^*\matr{M})\Q_{\Lambda_2 \cap \Lambda_1}}_2.
\end{eqnarray}
Using Lemma~\ref{lem:omega_RIP_norm} completes the proof.
\hfill $\Box$ \bigskip

As we will see later, we require the $\OM$-RIP to be small. Thus, we are interested to know for what matrices this hold true.
In the synthesis case, where $\matr{\Omega}$ is unitary and the $\OM$-RIP is identical to the RIP, it was shown for certain family of random matrices, such as matrices with Bernoulli or Subgaussian ensembles, that for any value of $\epsilon_k$ if $m \ge C_{\epsilon_k} k \log(\frac{m}{k\epsilon_k})$ then $\delta_k \le \epsilon_k$ \cite{Candes06Near, Rauhut08Compressed,Mendelson08Uniform}, where $\delta_k$ is the RIP constant and $C_{\epsilon_k}$ is a constant depending on $\epsilon_k$ and $\matr{M}$.
A similar result for the same family of random matrices holds for the analysis case.
The result is a special case of the result presented in \cite{Blumensath09Sampling}.

\begin{thm}[Theorem 3.3 in \cite{Blumensath09Sampling}]
\label{thm:analysis_RIP_cond}
Let $\matr{M}\in \RR{m \times d}$ be a random matrix such that for any $\vect{z}\in \RR{d}$ and $0<\tilde{\epsilon} \le \frac{1}{3}$
it satisfies
\begin{eqnarray}
P\left(\abs{\norm{\matr{M}\vect{z}}_2^2 - \norm{\vect{z}}_2^2} \ge \tilde{\epsilon}\norm{\vect{z}}_2^2\right) \le e^{-\frac{C_{\matr{M}}m\tilde{\epsilon}}{2}},
\end{eqnarray}
where $C_{\matr{M}} >0$ is a constant.
For any value of $\epsilon_\cosp >0$, if
\begin{eqnarray}
\label{eq:omega_RIP_m_size_cond}
m \ge \frac{32}{C_M\epsilon_r^2}\left( \log(\abs{\L_{r}^{\text{corank}}}) + (d-r)\log({9}/{\epsilon_r})+t\right),
\end{eqnarray}
then $\delta_{r}^{\text{corank}} \le \epsilon_r$ with probability exceeding $1-e^{-t}$.
\end{thm}

The above theorem is important since it shows that the $\OM$-RIP holds with a small constant for a large family of matrices --
the same family that satisfy the RIP property. In a recent work it was even shown that by randomizing the signs of the columns in the matrices that satisfy
the RIP we get new matrices that also satisfy the RIP \cite{Krahmer11New}. Thus, requiring the $\OM$-RIP constant to be small, as will be done hereafter, is legitimate.

For completeness we present a proof for theorem~\ref{thm:analysis_RIP_cond} in \ref{sec:analysis_RIP_proof} based on \cite{Rauhut08Compressed, Mendelson08Uniform, Baraniuk08Simple}.
We include in it also the proof of Theorem~\ref{thm:analysis_RIP_cond_dependencies} to follow.
In the case that $\matr{\Omega}$ is in general position $\abs{\L_{r}^{\text{corank}}} = {p \choose r} \le (\frac{ep}{p-r})^{p-r}$ (inequality is by Stirling's formula)
and thus $m \ge (p-r)\log(\frac{ep}{p-r})$. Since we want $m$ to be smaller than $d$ we need $p -\cosp$ to be smaller than $d$.
This limits the size of $p$ for $\OM$ since $r$ cannot be greater than $d$.
Thus, we present a variation of the theorem which states the results in terms of $\delta_{\cosp}$ and $\cosp$ instead of $\delta_{r}^{\text{corank}}$ and $r$.
The following theorem is also important because of the fact that our theoretical results are in terms of $\delta_{\cosp}$ and not $\delta_{r}^{\text{corank}}$.
It shows that $\delta_{\cosp}$ is small in the same family of matrices that guarantees $\delta_{r}^{\text{corank}}$ to be small.

\begin{thm}
\label{thm:analysis_RIP_cond_dependencies}
Under the same setup of Theorem~\ref{thm:analysis_RIP_cond}, for any $\epsilon_\cosp > 0$ if
\begin{eqnarray}
\label{eq:omega_RIP_m_size_cond_dependencies}
m \ge \frac{32}{C_M\epsilon_\cosp^2}\left( (p-\cosp)\log\left(\frac{9p}{(p-\cosp)\epsilon_\cosp}\right)+t \right),
\end{eqnarray}
then $\delta_\cosp \le \epsilon_\cosp$ with probability exceeding $1-e^{-t}$.
\end{thm}

Remark that when $\OM$ is in general position
$\cosp$ cannot be greater than $d$ and thus $p$ cannot be greater than $2d$ \cite{Nam12Cosparse}.
For this reason, if we want to have large values for $p$ we should allow linear dependencies between the rows of $\OM$.
In this case the cosparsity of the signal can be greater than $d$.
This explains why linear dependencies are a favorable thing in analysis dictionaries \cite{Rubinstein12Cosparse}.
In Section~\ref{sec:exp} we shall see that also empirically we get a better recovery when $\OM$ contains linear dependencies.

\section{Near Optimal Projection}
\label{sec:near_opt_proj}
As we will see hereafter, in the proposed algorithms we will face the following problem:
Given a general vector $\vect{z} \in \RR{\sdim}$, we would like to find an $\cosp$-cosparse vector that is closest to it in the $\ell_2$-norm sense.
In other words, we would like to project the vector to the closest $\cosp$-cosparse subspace.
Given the cosupport $\Lambda$ of this space the solution is simply $\Q_\Lambda\vect{z}$.
Thus, the problem of finding the closest $\cosp$-cosparse vector turns to be
the problem of finding the cosupport of the closest $\cosp$-cosparse subspace.
We denote the procedure of finding this cosupport by
\begin{eqnarray}
\label{eq:optimal_cosparse_projection}
\CF^*_{\cosp}(\vect{z}) = \argmin_{\Lambda \in \L_{\cosp}} \norm{\vect{z} - \Q_{\Lambda}\vect{z}}_2^2.
\end{eqnarray}
In the representation domain in the synthesis case, the support of the closest $k$-sparse subspace is found simply by hard thresholding, i.e., taking the support of the $k$-largest elements.
However, in the analysis case calculating \eqref{eq:optimal_cosparse_projection} is NP-complete with no efficient method for doing it
for a general $\OM$ \cite{Gribonval12Projection}. Thus an approximation procedure $\hat{\CF}_\cosp$ is needed. For this purpose we introduce the definition of a near-optimal projection \cite{Giryes11Iterative}.
\begin{defn}
\label{def:C_optimal_proj}
A procedure $\hat{\CF}_\cosp$ implies a near-optimal projection $\Q_{\hat{\CF}_\cosp(\cdot)}$ with a constant $C_\cosp$ if for any $\vect{z} \in \Real^d$
\begin{eqnarray}
\label{eq:C_optimal_proj}
&& \norm{\vect{z}-\Q_{\hat\CF_\cosp(\vect{z})}\vect{z}}_2^2 \le C_\cosp\norm{\vect{z} - \Q_{\CF^*_\cosp(\vect{z})}\vect{z}}_2^2.
\end{eqnarray}
\end{defn}
A clear implication of this definition is that if $\hat\CF_\cosp$ implies a near-optimal projection with a constant $C_\cosp$ then for any vector $\vect{z}\in \Real^d$ and an $\cosp$-cosparse vector $\vect{v} \in \Real^d$
\begin{eqnarray}
\label{eq:C_optimal_ineq}
&& \norm{\vect{z}-\Q_{\hat\CF_\cosp(\vect{z})}{\vect{z}}}_2^2 \le C_\cosp \norm{\vect{z}-\vect{v}}_2^2.
\end{eqnarray}

Similarly to the $\OM$-RIP, the above discussion can be directed also for finding the closest vector with corank $r$
defining $\CF_{r}^{\text{corank}*}$ and near optimal projection for this case in a very similar way to
\eqref{eq:optimal_cosparse_projection} and Definition~\ref{def:C_optimal_proj} respectively.

Having a near-optimal cosupport selection scheme  for a general operator is still an open problem
and we leave it for a future work. It is possible that this is also NP-complete.
We start by describing a simple thresholding rule that can be used with any operator.
Even though it does not have any known (near) optimality guarantee besides the case of unitary operators,
the numerical section will show it performs well in practice.
Then we present two tractable algorithms for finding the optimal cosupport
for two non-trivial analysis operators,
the one dimensional finite difference operator $\OM_{\text{1D-DIF}}$ \cite{Han04Optimal}
and the fused Lasso operator $\OM_{\text{FUS}}$ \cite{Tibshirani05Sparsity}.

Later in the paper, we propose theoretical guarantees for algorithms that use operators that has an
optimal or a near-optimal cosupport selection scheme. We leave the theoretical study of the thresholding technique for
a future work but demonstrate its performance empirically in Section~\ref{sec:exp} where this rule is used
showing that also when near-optimality is not at hand reconstruction is feasible.

\subsection{Cosupport Selection by Thresholding}
\label{sec:thresholding}

One intuitive option for cosupport selection is the simple thresholding
\begin{eqnarray}
\label{eq:thresh_cosupp_selection}
\hat\CF_\cosp(\vect{z}) = \cosupp(\OM\vect{z},\cosp),
\end{eqnarray}
which selects as a cosupport the indices of the $\cosp$-smallest elements after applying $\OM$ on $\z$.
As mentioned above, this selection method is optimal for unitary analysis operators where it coincides with
the hard thresholding used in synthesis.
However, in the general case this selection method is not guaranteed to give the optimal cosupport.
Its near optimality constant $C_\cosp$ is not close to one and is equal to the fraction
of the largest and smallest eigenvalues (which are not zero) of the submatrices composed of $\cosp$ rows from $\OM$ \cite{Giryes11Iterative}.

\begin{figure*}[!t]
{\subfigure[The signal $\vect{z}$.]{\includegraphics[width=1.5in]{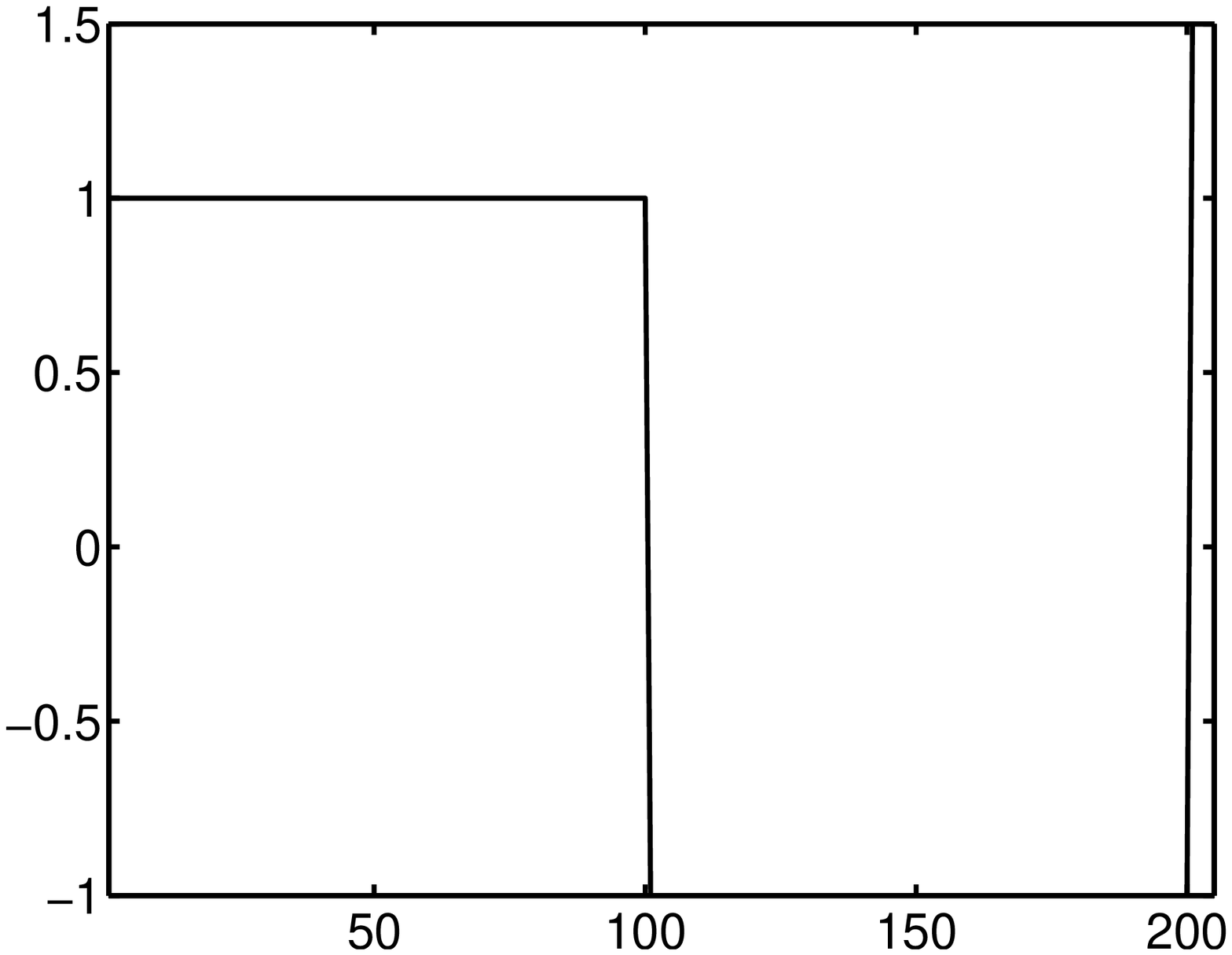}\label{fig:z_sig}}
\hfil
\subfigure[$\OM_{\text{1D-DIF}}\vect{z}$.]{\includegraphics[width=1.5in]{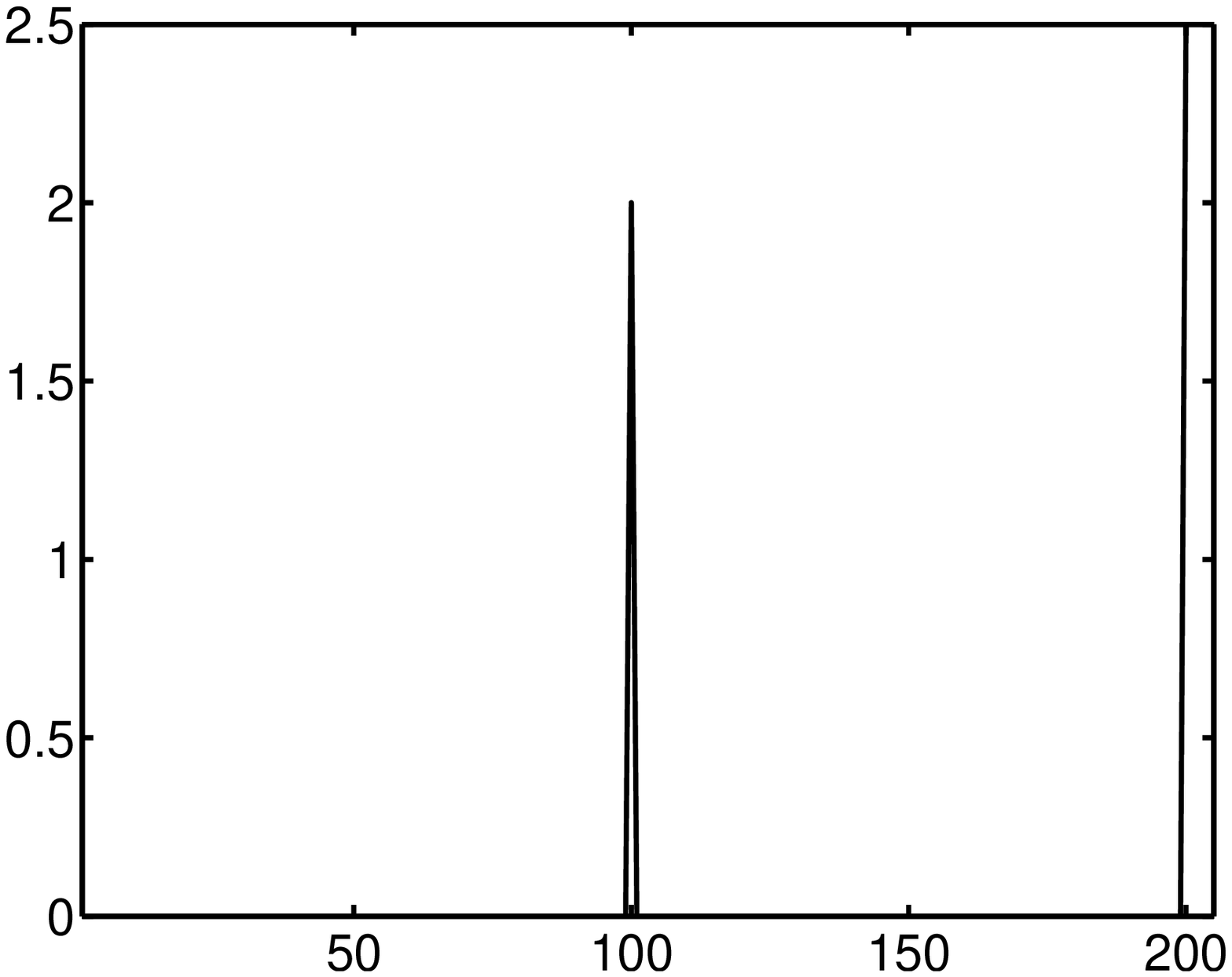}\label{fig:omega_z_sig}}%
\hfil
\subfigure[Projection using thresholding cosupport selection. The projection $\ell_2$-norm error is $\sqrt{200}$.]{\includegraphics[width=1.5in]{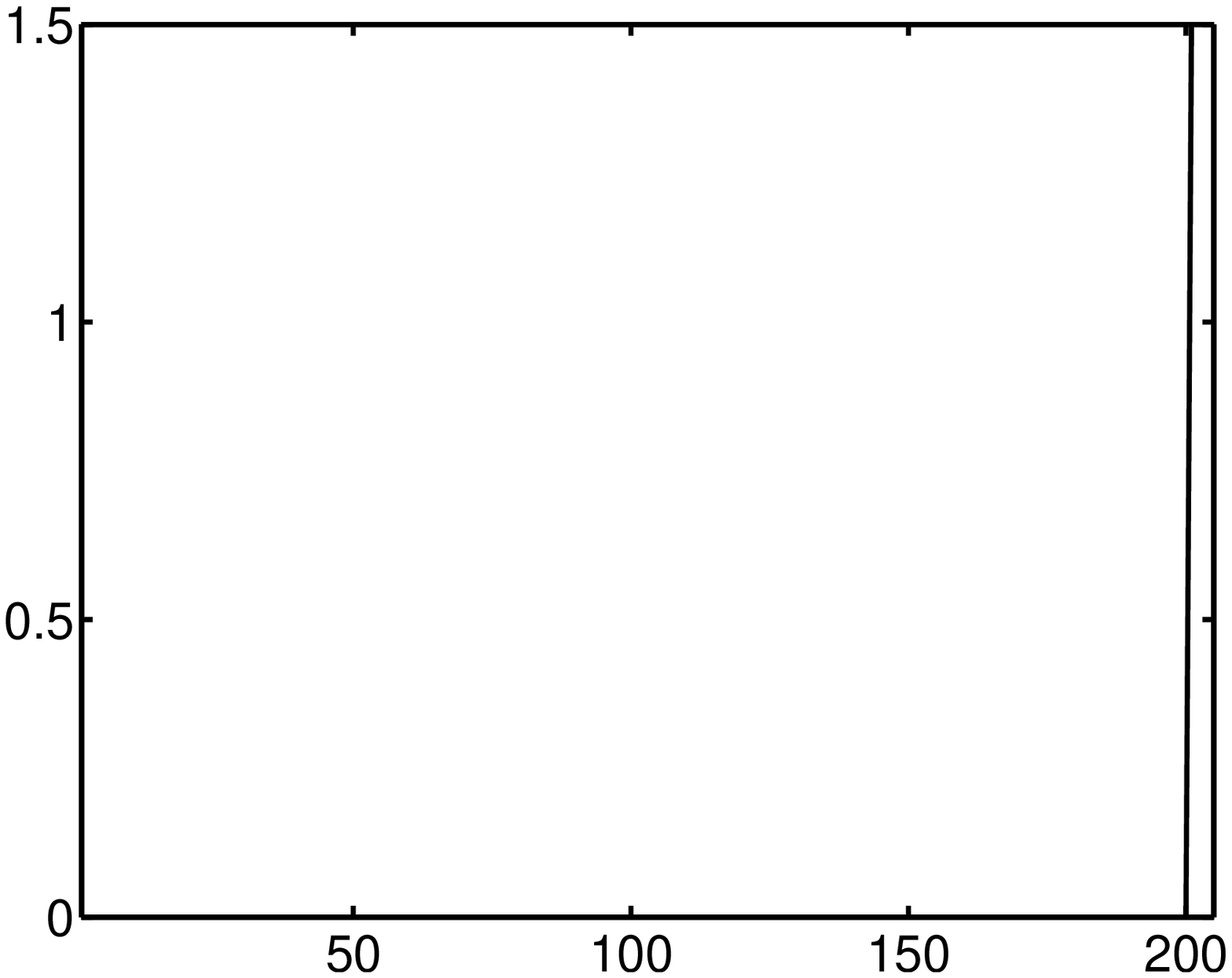}\label{fig:z_thresh_proj}}
\hfil
\subfigure[Optimal projection. The projection $\ell_2$-norm error is $2.5$.]{\includegraphics[width=1.5in]{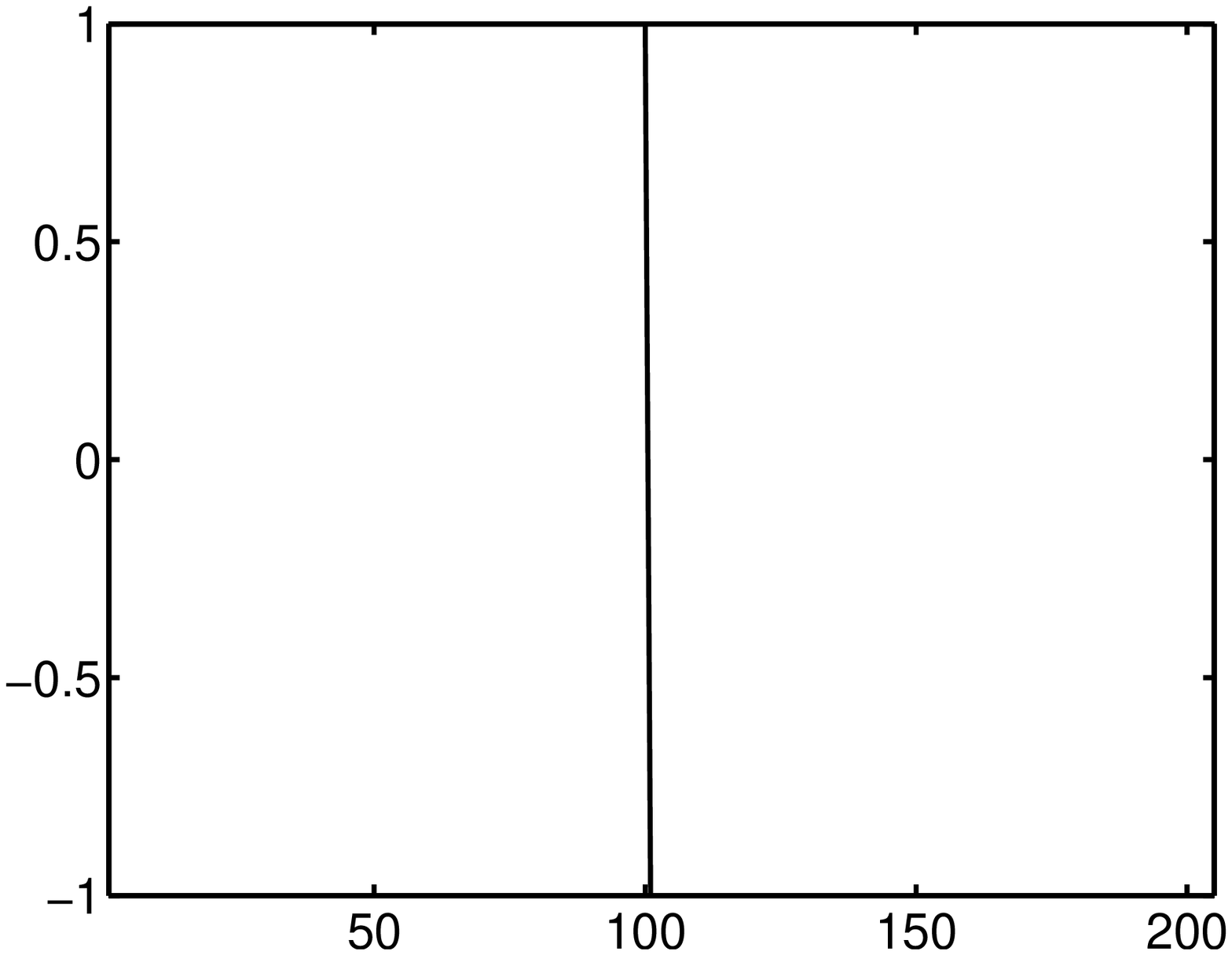}\label{fig:z_opt_proj}}}
\caption{Comparison between projection using thresholding cosupport selection and optimal cosupport selection. As it can be seen the thresholding
projection error is much larger than the optimal projection error by a factor much larger than $1$}
\label{fig:projection}
\end{figure*}

One example for an operator for which the thresholding is sub-optimal is the 1D-finite difference operator $\OM_{\text{1D-DIF}}$.
This operator is defined as:
\begin{equation}\label{1D Diff}
\OM_{\text{1D-DIF}} = \begin{pmatrix}
-1     & 1 &\cdots &~ &~ & ~ \\
\vdots &-1 & 1     & ~&~ & ~ \\
~      &~  &~      &\ddots  &~ & ~ \\
~      &~  &~      &~ &-1&1
\end{pmatrix}
\end{equation}
In this case, given a signal $\vect{z}$, applying $\OM_{\text{1D-DIF}}$ on it, result with a vector of coefficients that represents the differences in the signal. The thresholding selection method will select the indices of the $\cosp$ smallest elements in $\OM_{}\vect{z}$ as the cosupport $\Lambda_\vect{z}$.
For example, for the signal $\vect{z}\in \RR{201}$ in Fig~\ref{fig:z_sig} that contains $100$ times one, $100$ times minus one and $1.5$ as the last element, the thresholding
will select the cosupport to be the first $199$ coefficients in $\OM_{\text{1D-DIF}}\vect{z}$ that appears in Fig~\ref{fig:omega_z_sig}
and thus the projected vector will be the one in Fig~\ref{fig:z_thresh_proj}.
Its error in the $\ell_2$-norm sense is $\sqrt{200}$. However, selecting the cosupport to be the first $99$ elements and last $100$ elements result with the projected vector in Fig.~\ref{fig:z_opt_proj},
which has a smaller projection error ($2.5$). Thus, it is clear that the thresholding is sub-optimal for $\OM_{\text{1D-DIF}}$.
In a similar way it is also sub-optimal for the 2D-finite difference operator $\OM_{\text{2D-DIF}}$ that returns the
vertical and horizontal differences of a two dimensional signal.
Though not optimal, the use of thresholding with this operator is illustrated in Section~\ref{sec:exp} demonstrating that
also when a good projection is not at hand, good reconstruction is still possible.

\subsection{Optimal Analysis Projection Operators}

As mentioned above, in general it would appear that determining the optimal projection is
computationally difficult with the only general solution being to fully enumerate the
projections onto all possible cosupports.
Here we highlight two cases where it is relatively easy (polynomial complexity)
to calculate the optimal cosparse projection.

\subsubsection{Case 1: 1D finite difference}

For the 1D finite difference operator the analysis operator is not
redundant ($\pdim = \sdim-1$) but neither is it invertible.
As we have seen, a simple thresholding does not provide us with the optimal cosparse projection.
Thus, in order to determine the best $\cosp$-cosparse approximation for a given vector $\vect{z}$
we take another route and note that
we are looking for the closest (in the $\ell_2$-norm sense to $\vect{z}$) piecewise constant vector with $\pdim-\cosp$
change-points. This problem has been solved previously in the signal
processing literature using dynamic
programming (DP), see for example: \cite{Han04Optimal}. Thus for this operator
it is possible to calculate the best cosparse representation in ${\cal
  O}(\sdim^2)$ operations. The existence of a DP
solution follows from the ordered localized nature of the finite difference
operator. To the best of our knowledge, there is no known extension to 2D finite difference.

\subsubsection{Case 2: Fused Lasso Operator}
A redundant operator related to the 1D finite difference operator is
the so-called fused Lasso operator,
usually used with the analysis $\ell_1$-minimization \cite{Tibshirani05Sparsity}.
This usually takes the form:

\begin{equation}\label{fusion operator}
\OM_{\text{FUS}} = \begin{pmatrix} \OM_{\text{1D-DIF}} \\ \epsilon \matr{I}
\end{pmatrix}.
\end{equation}

Like $\OM_{\text{1D-DIF}}$ this operator works locally and therefore we can
expect to derive a DP solution to the approximation problem. This is
presented below.

\begin{rem} Note that
in terms of the cosparsity model the
$\epsilon$ parameter plays no role. This is in contrast to the
traditional convex optimization solutions where the value of
$\epsilon$ is pivotal \cite{Vaiter12Robust}. It is possible to mimic the $\epsilon$
dependence within the cosparsity framework by considering a
generalized fused Lasso operator of the form:
\begin{equation}\label{fusion operator_epsilon}
\OM_{\epsilon {\text{FUS}}} = \begin{pmatrix} \OM_{\text{1D-DIF}}\\ \OM_{\text{1D-DIF}}\\
  \vdots \\
  \OM_{\text{1D-DIF}} \\ \matr{I}
\end{pmatrix}.
\end{equation}
where the number of repetitions of the $\OM_{\text{1D-DIF}}$ operator (and
possibly the $\matr{I}$ operator) can be selected to mimic a weight
on the number of nonzero coefficients of each type. For simplicity
 we only consider the case indicated by
\eqref{fusion operator}
\end{rem}

\subsubsection{A recursive solution to the optimal projector for $\OM_{\text{FUS}}$}
Rather than working directly with the operator $\OM_{\text{FUS}}$ we make
use of the following observation.
An $\cosp$-cosparse vector $\vect{v}$ (or $k$-sparse vector) for $\OM_{\text{FUS}}$
is a piecewise
constant vector with $k_1$ change points and $k_2$ non-zero entries
such that $k_1 + k_2 = k = p-\cosp$, where $p=2d-1$. To understand better the relation between
$k_1$ and $k_2$, notice that $k_1=0$ implies equality of all entries,
so $k_2 = 0$ or $d$, hence $\cosp = p$ or $d-1$.
Conversely, considering $d\le\cosp<p$ or $0\le\cosp<d-1$ implies $k_1 \neq 0$.
It also implies that there is at least one nonzero value, hence $k_2 \neq 0$.

Thus, an $\cosp$-cosparse vector $\vect{v}$ for $\OM_{\text{FUS}}$ can be parameterized in
terms of a set of change points, $\{n_i\}_{i=0:k_1+1}$, and a set of
constants, $\{\mu_i\}_{i=1:k_1+1}$, such that:
\begin{equation}
\vect{v}_j = \mu_i, n_{i-1} < j \leq n_i
\end{equation}
with the convention that $n_0 = 0$ and $n_{k_1+1} = d$, unless stated otherwise. We will also
make use of the indicator vector, $\s$, defined as:
\begin{equation} \s_i =
\begin{cases}
0& \text{if $\mu_i = 0$},\\
1& \text{otherwise}
\end{cases} \text{for $1\le i \le k_1+1$}.
\end{equation}
Using this alternative
parametrization we can write the minimum distance between a vector $\vect{z}$
and the set of $k$-sparse fused Lasso coefficients as:
\begin{equation}
\begin{split}
F_k(\vect{z}) &= \min_{1\leq k_1 \leq k} \min_{\substack{
\{n_i\}_{i=1:k_1}\\
\{\mu_i\}_{i=1:k_1+1}\\
n_{k_1}<d}}
\sum_{i=1}^{k_1+1} \sum_{j = n_{i-1}+1}^{n_i}(\vect{z}_j-\mu_i)^2,\\
&\mbox{~subject to~} \sum_{i = 1}^{k_1+1} \s_i (n_i-n_{i-1}) = k-k_1
\end{split}
\end{equation}
Although this looks a formidable optimization task we now show that
it can be computed recursively through a standard DP strategy,
modifying the arguments in \cite{Han04Optimal}.

Let us define the optimal cost, $I_k(L,\omega,k_1)$, for the vector $[\vect{z}_1, \ldots ,
\vect{z}_L]^T$ with $k_1$ change points and $\s_{k_1+1} = \omega$, as:
\begin{equation}
\begin{split}
I_k(L,\omega,k_1) &= \min_{\substack{
\{n_i\}_{i=1:k_1}\\
\{\s_i\}_{i=1:k_1+1}\\
n_{k_1} <L,~ n_{k_1+1} =L\\
\s_{k_1+1}=\omega}}
\sum_{i=1}^{k_1+1} \sum_{j = n_{i-1}+1}^{n_i}(\vect{z}_j-\mu_i)^2,\\
&\mbox{~subject to~} \sum_{i = 1}^{k_1+1} \s_i (n_i-n_{i-1}) = k-k_1\\
&\mbox{~and~} \mu_i = \frac{\s_i}{n_i-n_{i-1}}
\sum_{l=n_{i-1}+1}^{n_i} \vect{z}_l
\end{split}
\end{equation}
where we have set $\mu_i$ to the optimal sample
means. Notice that calculating $I_k(L,\omega,k_1)$ is easy for $k_1\le k \le 1$.
Thus, we calculate it recursively considering two separate
scenarios:
\begin{description}
\item[Case 1: $\omega = 0$] where the last block of coefficients are zero.
This gives:
\begin{equation}
\begin{split}
I_k(L,0,k_1) = &  \min_{n_{k_1<L}}\left( \sum_{j = n_{k_1}+1}^L (\vect{z}_j)^2 +
\min_{\substack{\{n_i\}_{i=1:k_1-1}\\
\{\s_i\}_{i=1:k_1-1}\\
n_{k_1-1}<n_{k_1}\\
\s_{k_1} = 1}} \sum_{i=1}^{k_1} \sum_{j =
n_{i-1}+1}^{n_i}(\vect{z}_j-\mu_i)^2 \right),\\
~&\mbox{~subject to~} \sum_{i = 1}^{k_1} s_i (n_i-n_{i-1}) = (k-1)-(k_1-1)\\
~&\mbox{~and~} \mu_i = \frac{\s_i}{n_i-n_{i-1}} \sum_{l=n_{i-1}+1}^{n_i} \vect{z}_l,
\end{split}
\end{equation}
(noting that if $\s_{k_1+1} = 0$ then $\s_{k_1} = 1$ since otherwise $n_{k_1}$ would not have been a change point). This simplifies
to the recursive formula:
\begin{equation}\label{eq:recusive1}
I_k(L,0,k_1) = \min_{n_{k_1}< L}\left( \sum_{j = n_{k_1}+1}^L (\vect{z}_j)^2 +
I_{k-1}(n_{k_1},1,k_1-1)\right)
\end{equation}

\item[Case 2: $\omega = 1$] when the final block of coefficients are
  non-zero we have:
\begin{equation}
\begin{split}
I_k(L,1,k_1) = &  \min_{
\substack{n_{k_1}<L \\n_{k_1+1}=L\\ \s_{k_1}} }
\left( \sum_{j = n_{k_1}+1}^L (\vect{z}_j-\mu_{k_1+1})^2 +
\min_{\substack{\{n_i\}_{i=1:k_1-1}\\
\{\s_i\}_{i=1:k_1-1}\\
n_{k_1-1}<n_{k_1}}}
\sum_{i=1}^{k_1} \sum_{j = n_{i-1}+1}^{n_i}(\vect{z}_j-\mu_i)^2 \right),\\
~&\mbox{~subject to~} \sum_{i = 1}^{k_1} \s_i (n_i-n_{i-1}) = (k-L+n_{k_1}-1)-(k_1-1)\\
~&\mbox{~and~} \mu_i = \frac{\s_i}{n_i-n_{i-1}} \sum_{l=n_{i-1}+1}^{n_i} \vect{z}_l.
\end{split}
\end{equation}
This simplifies to the recursive relationship:
\begin{equation}\label{eq:recursive2}
\begin{split}
I_k(L,1,k_1) &= \min_{\substack{n_{k_1}<L\\ \s_{k_1}} } \left(\sum_{j = n_{k_1}+1}^L (\vect{z}_j-\mu_{k_1+1})^2 +
I_{k-L+n_{k_1}-1}(n_{k_1},\s_{k_1},k_1-1) \right)\\
&\text{subject to~}  \mu_{k_1+1} = \sum_{l=n_{k_1}+1}^{L} \vect{z}_l/\big(L-n_{k_1}\big)
\end{split}
\end{equation}
\end{description}

Equations \eqref{eq:recusive1} and \eqref{eq:recursive2}
are sufficient to enable the calculation of the optimal projection in
polynomial time,starting with $k_1\le k\le 1$ and recursively evaluating
the costs for $k\ge k_1 \ge 1$. Finally, we have $F_k(\vect{z}) = \min_{k_1\le k, \omega \in \{0,1\}} I_k(d,\omega,k_1)$. The implementation details are left as an
exercise for the reader.

\section{New Analysis algorithms}
\label{sec:analysis_alg}

\subsection{Quick Review of the Greedy-Like Methods}
Before we turn to present the analysis versions of the greedy-like techniques we recall their synthesis versions.
These use a prior knowledge about the cardinality $k$ and actually
aim at approximating a variant of \eqref{eq:synthesisL0}
\begin{equation}
\label{eq:synthesisL0_k}
\argmin_{\alphabf} \norm{\y - \M \Dict \alphabf}_2^2 \subjectto \norm{\alphabf}_0 \le k.
\end{equation}
For simplicity we shall present the greedy-like pursuits for the case $\matr{D} = \matr{I}$.
In the general case $\matr{M}$ should be replaced with $\matr{M}\matr{D}$, $\vect{x}$ with $\alphabf$ and the reconstruction result should be $\hat\x = \matr{D}\hat\alphabf$.
In addition, in the algorithms' description we do not specify the stopping criterion.
Any standard stopping criterion, like residual's size or relative iteration change, can be used.
More details can be found in \cite{Needell09CoSaMP,Dai09Subspace}.

{\em IHT and HTP:} IHT \cite{Blumensath09Iterative} and HTP \cite{Foucart11Hard} are presented in Algorithm~\ref{alg:IHT_HTP}.
Each IHT iteration is composed of two basic steps.
The first is a gradient step, with a step size $\mu_t$,
in the direction of minimizing $\norm{\y - \M \x}_2^2$.
The step size can be either constant in all iterations ($\mu^t = \mu$) or changing \cite{Kyrillidis11Recipes}.
The result vector $\x_g$ is not guaranteed to be sparse and thus
the second step of IHT projects $\x_g$ to the closest $k$-sparse subspace
by keeping its largest $k$ elements.
The HTP takes a different strategy in the projection step.
Instead of using a simple projection to the closest $k$-sparse subspace, HTP
selects the vector in this subspace that minimizes $\norm{\y - \M \x}_2^2$ \cite{Foucart11Hard, Blumensath12Accelerated}.

\begin{algorithm}
\caption{Iterative hard thresholding (IHT) and hard thresholding pursuit (HTP)} \label{alg:IHT_HTP}
\begin{algorithmic}

\REQUIRE $k, \matr{M}, \vect{y}$ where $\vect{y} = \matr{M}\vect{x}
+ \vect{e}$, $k$ is the cardinality of $\vect{x}$ and $\vect{e}$ is
an additive noise.

\ENSURE $\hat{\vect{x}}_{\text{\tiny IHT}}$ or $\hat{\vect{x}}_{\text{\tiny HTP}}$: $k$-sparse approximation of
$\vect{x}$.

\STATE Initialize representation $\hat{\x}^0 = \vect{0}$
 and set $t = 0$.

\WHILE{halting criterion is not satisfied}

\STATE $t = t + 1$.

\STATE Perform a gradient step: $\vect{x}_g = \hat{\vect{x}}^{t-1} + \mu^t \M^*(\y - \M\hat{\x}^{t-1})$

\STATE Find a new support: $T^t = \supp(\x_g,k)$

\STATE Calculate a new representation: $\hat{\vect{x}}_{\text{\tiny IHT}}^t = (\vect{x}_g)_{T^t}$ for IHT, and
$\hat{\vect{x}}_{\text{\tiny HTP}}^t = \matr{M}^{\dag}_{T^t}\vect{y}$ for HTP.

\ENDWHILE

\STATE Form the final solution $\hat{\vect{x}}_{\text{\tiny IHT}} = \hat{\vect{x}}_{\text{\tiny IHT}}^t$ for IHT and $\hat{\vect{x}}_{\text{\tiny HTP}} = \hat{\vect{x}}_{\text{\tiny HTP}}^t$ for HTP.

\end{algorithmic}
\end{algorithm}

{\em CoSaMP and SP:} CoSaMP \cite{Needell09CoSaMP} and SP \cite{Dai09Subspace} are presented in Algorithm~\ref{alg:CoSaMP_SP}. The difference between these two
techniques is similar to the difference between IHT and HTP. Unlike IHT and HTP,
the estimate for the support of $\x$ in each CoSaMP and SP iteration is computed
by observing the residual $\vect{y}_{\text{resid}}^t = \vect{y} - \M\x^t$.
In each iteration, CoSaMP and SP extract new support indices from the residual by taking the indices of the largest
elements in $\matr{M}^*\vect{y}_{\text{resid}}^t$.  They add the new indices to the estimated support set from the previous iteration
creating a new estimated support $\tilde{T}^t$ with cardinality larger than $k$.
Having the updated support, in a similar way to the projection in HTP, an objective aware
projection is performed resulting with an estimate $\vect{w}$ for $\x$ that is supported on $\tilde{T}^t$.
Since we know that $\x$ is $k$-sparse we want to project $\vect{w}$ to a $k$-sparse subspace.
CoSaMP does it by simple hard thresholding like in IHT. SP does it by an objective aware projection
similar to HTP.

\begin{algorithm}[t]
\caption{Subspace Pursuit (SP) and CoSaMP} \label{alg:CoSaMP_SP}
\begin{algorithmic}

\REQUIRE $k, \matr{M}, \vect{y}$ where $\vect{y} = \matr{M}\vect{x}
+ \vect{e}$, $k$ is the cardinality of $\vect{x}$ and $\vect{e}$ is
an additive noise. $a = 1$ (SP), $a = 2$ (CoSaMP).

\ENSURE $\hat{\vect{x}}_{\text{\tiny CoSaMP}}$ or $\hat{\vect{x}}_{\text{\tiny SP}}$: $k$-sparse approximation of
$\vect{x}$.

\STATE Initialize the support $T^0 =\emptyset$, the residual $\vect{y}_{\text{resid}}^0 = \vect{y}$
 and set $t = 0$.

\WHILE{halting criterion is not satisfied}

\STATE $t = t + 1$.

\STATE Find new support elements: $T_\Delta =
\supp(\matr{M}^*\vect{y}^{t - 1}_{\text{resid}},a k)$.

\STATE Update the support: $\tilde{T}^t = T^{t -1} \cup
T_\Delta$.

\STATE Compute a temporary representation: $\vect{w} =
\matr{M}^{\dag}_{\tilde{T}^t}\vect{y}$.

\STATE Prune small entries: $T^t =
\supp(\vect{w},k)$.

\STATE Calculate a new representation: $\hat{\vect{x}}_{\text{\tiny CoSaMP}}^t = \vect{w}_{T^t}$ for CoSaMP, and
$\hat{\vect{x}}_{\text{\tiny SP}}^t = \matr{M}^{\dag}_{T^t}\vect{y}$ for SP.

\STATE Update the residual:
$\vect{y}_{\text{resid}}^t = \vect{y} - \matr{M}\hat{\vect{x}}_{\text{\tiny CoSaMP}}^t$
 for CoSaMP, and $\vect{y}_{\text{resid}}^t = \vect{y} -\matr{M}\hat{\vect{x}}_{\text{\tiny SP}}^t$ for SP.

\ENDWHILE
\STATE Form the final solution $\hat{\vect{x}}_{\text{\tiny CoSaMP}} = \hat{\vect{x}}_{\text{\tiny CoSaMP}}^t$ for CoSaMP and $\hat{\vect{x}}_{\text{\tiny SP}} = \hat{\vect{x}}_{\text{\tiny SP}}^t$ for SP.

\end{algorithmic}
\end{algorithm}

\subsection{Analysis greedy-like methods}

\begin{table*}
   \centering
  \begin{tabular}{||C{3.6cm}|C{3.5cm}|| C{4cm}| C{3.5cm}||}
  \hline
  {\bf Synthesis operation name} & {\bf Synthesis operation} & {\bf Analysis operation name }  & {\bf Analysis operation }    \\ \hline
  Support selection & Largest $k$ elements: $T=\supp(\cdot, k)$ & Cosupport selection  & Using a near optimal projection: $\Lambda = \hat\CF_\cosp(\cdot)$ \\ \hline
  Orthogonal Projection of $\vect{z}$ to a $k$-sparse subspace with support $T$ & $\vect{z}_T$  & Orthogonal projection of $\vect{z}$ to an $\cosp$-cosparse subspace with cosupport $\Lambda$  & $\Q_{\Lambda}\vect{z}$ \\ \hline
  Objective aware projection to a $k$-sparse subspace with support $T$ & $\M_T^\dag \y$ = $\argmin_{\vect{v}}\norm{\y - \M\vect{v}}_2^2$ s.t. $\vect{v}_{T^C} = 0$ & Objective aware projection to an $\cosp$-cosparse subspace with cosupport $\Lambda$ & $\argmin_{\vect{v}}\norm{\y - \M\vect{v}}_2^2$ s.t. $\OM_{\Lambda}\vect{v} = 0$  \\ \hline
  Support of $\vect{v}_1 + \vect{v}_2$ where $\supp(\vect{v}_1)= T_1$ and $\supp(\vect{v}_2)= T_2$ & $\supp(\vect{v}_1+\vect{v}_2)\subseteq T_1 \cup T_2$ & Cosupport of $\vect{v}_1 + \vect{v}_2$ where $\cosupp(\vect{v}_1)= \Lambda_1$ and $\cosupp(\vect{v}_2)= \Lambda_2$  &
  $\cosupp(\vect{v}_1+\vect{v}_2) \supseteq \Lambda_1 \cap \Lambda_2$ \\ \hline
  Maximal size of $T_1 \cup T_2$ where $\abs{T_1} \le k_1$ and $\abs{T_2} \le k_2$ & $\abs{T_1 \cup T_2} \le k_1+k_2$ & Minimal size of $\Lambda_1 \cap \Lambda_2$ where $\abs{\Lambda_1} \ge \cosp_1$ and $\abs{\Lambda_2} \ge \cosp_2$  & $\abs{\Lambda_1 \cap \Lambda_2} \ge \cosp_1 + \cosp_2 - p$ \\
  \hline
\end{tabular}
  \caption{Parallel synthesis and analysis operations}
  \label{tbl:Synthesis_Analysis_parallels}
\end{table*}

Given the synthesis greedy-like pursuits, we would like to define their analysis counterparts.
For this task we need to 'translate' each synthesis operation into an analysis one.
This gives us a general recipe for converting algorithms between the two schemes.
The parallel lines between the schemes are presented in Table~\ref{tbl:Synthesis_Analysis_parallels}.
Those become more intuitive and clear when we keep in mind that while the synthesis approach focuses on the non-zeros,
the analysis concentrates on the zeros.

For clarity we dwell a bit more on the equivalences.
For the cosupport selection, as mentioned in Section~\ref{sec:near_opt_proj},
computing the optimal cosupport is a combinatorial problem and thus
the approximation $\hat\CF_\cosp$ is used.
Having a selected cosupport $\Lambda$, the projection to its
corresponding cosparse subspace becomes trivial, given by $\Q_\Lambda$.

Given two vectors $\vect{v}_1 \in \A_{\cosp_1}$ and $\vect{v}_2 \in \A_{\cosp_2}$ such that
$\Lambda_1 = \cosupp(\OM\vect{v}_1)$ and $\Lambda_2 = \cosupp(\OM\vect{v}_2)$, we know that
$\abs{\Lambda_1} \ge \cosp_1$ and $\abs{\Lambda_2} \ge \cosp_2$.
Denoting $T_1 = \supp(\OM\vect{v}_1)$ and $T_2 = \supp(\OM\vect{v}_2)$ it is clear that $\supp(\OM(\vect{v}_1 + \vect{v}_1)) \subseteq T_1 \cup T_2$.
Noticing that $\supp(\cdot) =\cosupp(\cdot)^C$ it is clear that
$\abs{T_1}\le p-\cosp_1$, $\abs{T_2}\le p-\cosp_2$ and $\cosupp(\OM(\vect{v}_1 + \vect{v}_2)) \supseteq (T_1 \cup T_2)^C =T_1^C \cap T_2^C = \Lambda_1 \cap \Lambda_2$.
From the last equality we can also deduce that $\abs{\Lambda_1 \cap \Lambda_2} = p -\abs{T_1 \cup T_2} \ge p - (p-\cosp_1) - (p-\cosp_2) = \cosp_1 + \cosp_2 - p$.

With the above observations we can develop the analysis versions of the greedy-like algorithms.
As in the synthesis case, we do not specify a stopping criterion.
Any stopping criterion used for the synthesis versions can be used also for the analysis ones.

\begin{algorithm}
\caption{Analysis Iterative hard thresholding (AIHT) and analysis hard thresholding pursuit (AHTP)} \label{alg:Analysis_IHT_HTP}
\begin{algorithmic}

\REQUIRE $\cosp, \matr{M}, \matr{\Omega}, \vect{y}$ where $\vect{y} = \matr{M}\vect{x}
+ \vect{e}$, $\cosp$ is the cosparsity of $\vect{x}$ under $\matr{\Omega}$ and $\vect{e}$ is
the additive noise.

\ENSURE $\hat{\vect{x}}_{\text{\tiny AIHT}}$ or $\hat{\vect{x}}_{\text{\tiny AHTP}}$: $\cosp$-cosparse approximation of
$\vect{x}$.

\STATE Initialize estimate $\hat{\x}^0 = \vect{0}$ and set $t = 0$.

\WHILE{halting criterion is not satisfied}

\STATE $t = t + 1$.

\STATE Perform a gradient step: $\vect{x}_g = \hat{\vect{x}}^{t-1} + \mu^t \M^*(\y - \M\hat{\x}^{t-1})$

\STATE Find a new cosupport: $\hat\Lambda^t = \hat{\CF}_\cosp(\x_g)$

\STATE Calculate a new estimate: $\hat{\vect{x}}_{\text{\tiny AIHT}}^t = \Q_{\hat\Lambda^t}\vect{x}_g$ for AIHT, and
$\hat{\vect{x}}_{\text{\tiny AHTP}}^t = \argmin_{\tilde{\vect{x}}}\norm{\y - \M\tilde{\vect{x}}}_2^2$ s.t. $\OM_{\hat\Lambda^t}\tilde{\vect{x}} = 0$ for AHTP.

\ENDWHILE

\STATE Form the final solution $\hat{\vect{x}}_{\text{\tiny AIHT}} = \hat{\vect{x}}_{\text{\tiny AIHT}}^t$ for AIHT and $\hat{\vect{x}}_{\text{\tiny AHTP}} = \hat{\vect{x}}_{\text{\tiny AHTP}}^t$ for AHTP.

\end{algorithmic}
\end{algorithm}

{\em AIHT and AHTP:}
Analysis IHT (AIHT) and analysis HTP (AHTP) are presented in Algorithm~\ref{alg:Analysis_IHT_HTP}.
As in the synthesis case, the choice of the gradient stepsize $\mu^t$ is crucial: If $\mu^t$'s are chosen too small,
the algorithm gets stuck at a wrong solution and if too large, the algorithm diverges.
We consider two options for $\mu^t$.

In the first we choose $\mu^t = \mu$ for some constant $\mu$ for all iterations.
A theoretical discussion on how to choose $\mu$ properly is given in Section~\ref{sec:AIHT_AHTP_guarantees}.

The second option is to select a different $\mu$ in each iteration.
One way for doing it is to choose an `optimal' stepsize $\mu^t$ by solving
the following problem
\begin{equation}
\label{eq:muOptimalM}
\mu^t := \argmin_\mu \norm{\y - \M\hat\x^t }_2^2.
\end{equation}
Since $\hat\Lambda^t = \hat\CF_\cosp(\hat\x^{t-1} + \mu^t\M^*(\y-\M\hat\x^{t-1}))$
and $\hat\x^t = \Q_{\hat\Lambda^t}(\x_g)$,
the above requires a line search over different values of $\mu$ and
along the search $\hat\Lambda^t$ might change several times.
A simpler way is an adaptive step
size selection as proposed in \cite{Kyrillidis11Recipes} for IHT.
In a heuristical way we limit the search to the cosupport
$\tilde{\Lambda} = \hat{\CF}_\cosp(\M^*(\y-\M\hat\x^{t-1})) \cap \hat{\Lambda}^{t-1}$.
This is the intersection of the cosupport of $\hat\x^{t-1}$ with the $\cosp$-cosparse cosupport of
the estimated closest $\cosp$-cosparse subspace to $\M^*(\y-\M\hat\x^{t-1})$.
Since $\hat\x^{t-1} = \Q_{\tilde\Lambda} \hat\x^{t-1}$, finding $\mu$ turns to be
\begin{equation}
\label{eq:muOptimalM_approx}
\mu^t := \argmin_\mu \norm{\y - \M(\hat\x^{t-1} + \mu\Q_{\tilde{\Lambda}}\M^*(\y-\M\hat\x^{t-1}))}_2^2.
\end{equation}
This procedure of selecting $\mu^t$ does not require a line search and it has a simple closed form solution.

To summarize, there are three main options for the step size selection:
\begin{itemize}
\item Constant step-size selection -- uses a constant step size $\mu^t = \mu$ in all iterations.
\item Optimal changing step-size selection -- uses different values for $\mu^t$ in each iterations by minimizing $\norm{\vect{y} -\M\hat\x^t}_2$.
\item Adaptive changing step-size selection -- uses \eqref{eq:muOptimalM_approx}.
\end{itemize}

\begin{algorithm}[t]
\caption{Analysis Subspace Pursuit (ASP) and Analysis CoSaMP (ACoSaMP)} \label{alg:Analysis_CoSaMP_SP}
\begin{algorithmic}[l]

\REQUIRE $\cosp, \matr{M}, \matr{\Omega}, \vect{y}, a$ where $\vect{y} = \matr{M}\vect{x}
+ \vect{e}$, $\cosp$ is the cosparsity of $\vect{x}$ under $\matr{\Omega}$ and $\vect{e}$ is
the additive noise.

\ENSURE $\hat{\vect{x}}_{\text{\tiny ACoSaMP}}$ or $\hat{\vect{x}}_{\text{\tiny ASP}}$: $\cosp$-cosparse approximation of
$\vect{x}$.

\STATE Initialize the cosupport $\Lambda^0 =\{i,1\le i \le p\}$, the residual $\vect{y}_{\text{resid}}^0 = \vect{y}$ and set $t = 0$.

\WHILE{halting criterion is not satisfied}

\STATE $t = t + 1$.

\STATE Find new cosupport elements: $\Lambda_\Delta=\hat{\CF}_{a\cosp}(\matr{M}^*\vect{y}^{t - 1}_{\text{resid}})$.

\STATE Update the cosupport: $\tilde{\Lambda}^t = \hat\Lambda^{t -1} \cap
\Lambda_\Delta$.

\STATE Compute a temporary estimate: $\vect{w} = \argmin_{\tilde{\vect{x}}} \norm{\vect{y} - \matr{M}\tilde{\vect{x}}}_2^2$ s.t. $\matr{\Omega}_{\tilde{\Lambda}^t}\tilde{\vect{x}} = 0$.

\STATE Enlarge the cosupport: $\hat\Lambda^t =
\hat{\CF}_{\cosp}(\vect{w})$.

\STATE \begin{flushleft}Calculate a new estimate: $\hat{\vect{x}}_{\text{\tiny ACoSaMP}}^t = \Q_{\hat\Lambda^t}\vect{w}$ for ACoSaMP, and
$\hat{\vect{x}}_{\text{\tiny ASP}}^t = \argmin_{\tilde{\vect{x}}} \norm{\vect{y} - \matr{M}\tilde{\vect{x}}}_2^2$ s.t. $\matr{\Omega}_{\hat\Lambda^t}\tilde{\vect{x}} = 0$
 for ASP.\end{flushleft}

\STATE \begin{flushleft}Update the residual:
$\vect{y}_{\text{resid}}^t = \vect{y} - \matr{M}\hat{\vect{x}}_{\text{\tiny ACoSaMP}}^t$
 for ACoSaMP, and $\vect{y}_{\text{resid}}^t = \vect{y} -\matr{M}\hat{\vect{x}}_{\text{\tiny ASP}}^t$ for ASP.\end{flushleft}

\ENDWHILE
\STATE Form the final solution $\hat{\vect{x}}_{\text{\tiny ACoSaMP}} = \hat{\vect{x}}_{\text{\tiny ACoSaMP}}^t$ for ACoSaMP and $\hat{\vect{x}}_{\text{\tiny ASP}} = \hat{\vect{x}}_{\text{\tiny ASP}}^t$ for ASP.
\end{algorithmic}
\end{algorithm}


{\em ACoSaMP and ASP:}
analysis CoSaMP (ACoSaMP) and analysis SP (ASP)
are presented in Algorithm~\ref{alg:Analysis_CoSaMP_SP}.
The stages are parallel to those of the synthesis CoSaMP and SP.
We dwell a bit more on the meaning of the parameter $a$ in the algorithms.
This parameter determines the size of the new cosupport $\Lambda_\Delta$ in each iteration.
$a=1$ means that the size is $\cosp$ and according to Table~\ref{tbl:Synthesis_Analysis_parallels}
it is equivalent to $a=1$ in the synthesis as done in SP in which we select new $k$ indices for the support in each iteration.
In synthesis CoSaMP we use $a=2$ and select $2k$ new elements. $2k$ is the maximal support size of two added $k$-sparse vectors.
The corresponding minimal size in the analysis case is $2\cosp-p$ according to Table~\ref{tbl:Synthesis_Analysis_parallels}.
For this setting we need to choose $a = \frac{2\cosp-p}{\cosp}$.

\subsection{The Unitary Case}
For $\matr{\Omega} = \matr{I}$
the synthesis and the analysis greedy-like algorithms become equivalent.
This is easy to see since in this case we have $p=d$, $k = d-\cosp$,  $\Lambda = T^C$,
$\Q_\Lambda\x = \x_T$ and $T_1 \cup T_2 = \Lambda_1 \cap \Lambda_2$ for $\Lambda_1 = T_1^C$ and $\Lambda_2 = T_2^C$.
In addition, $\hat{\CF}_{\cosp} = \CF_{\cosp}^*$ finds the closest $\cosp$-cosparse subspace
by simply taking the smallest $\cosp$ elements. Using similar arguments, also in the case where $\OM$ is a unitary matrix
the analysis methods coincide with the synthesis ones. In order to get exactly the same algorithms
$\M$ is replaced with $\M\OM^*$ in the synthesis techniques and the output is multiplied by $\OM^*$.

Based on this observation, we can deduce that the guarantees of the synthesis greedy-like methods apply
also for the analysis ones in a trivial way. Thus,
it is tempting to assume that the last should have similar guarantees
based on the $\OM$-RIP. In the next section we develop such claims.

\subsection{Relaxed Versions for High Dimensional Problems}
\label{sec:relaxed_alg}
Before moving to the next section we mention a variation of the analysis greedy-like techniques.
In AHTP, ACoSaMP and ASP  we need to solve the constrained minimization problem
$\min_{\tilde{\vect{x}}}\norm{\vect{y}-\matr{M}\tilde{\vect{x}}}_2^2$ s.t. $\norm{\OM_{\Lambda}\tilde{\vect{x}}}_2^2 =0$.
For high dimensional signals this problem is hard to solve and we suggest to replace it with
minimizing $\norm{\vect{y}-\matr{M}\tilde{\vect{x}}}_2^2 + \lambda\norm{\OM_{\Lambda}\tilde{\vect{x}}}_2^2$,
where $\lambda$ is a relaxation constant.
This results in a relaxed version of the algorithms.
We refer hereafter to these versions as relaxed AHTP (RAHTP) relaxed ASP (RASP) and relaxed ACoSaMP (RACoSaMP).


\section{Algorithms Guarantees}
\label{sec:guarantees}

In this section we provide theoretical guarantees for the reconstruction performance of the analysis greedy-like methods.
For AIHT and AHTP we study both the constant step-size and the optimal step-size selections.
For ACoSaMP and ASP the analysis is made for $a = \frac{2\cosp-p}{\cosp}$,
but we believe that it can be extended also to other values of $a$, such as $a=1$.
The performance guarantees we provide are summarized in the following two theorems.
The first theorem, for AIHT and AHTP, is a simplified version of Theorem~\ref{thm:AIHT_AHTP_theorem_noisy_iter}
and the second theorem, for ASP and ACoSaMP, is a combination of Corollaries~\ref{cor:ACoSaMP_bound} and \ref{cor:ASP_bound},
all of which appear hereafter along with their proofs.
Before presenting the theorems we recall the problem we aim at solving:
\begin{defn}[Problem $\cal{P}$]
Consider a measurement vector $\y \in \RR{\mdim}$
such that $\y=\M\x + \e$ where $\x\in \RR{\sdim}$ is $\cosp$-cosparse,
$\M\in \RR{\mdim\times \sdim}$ is a degradation operator and $\e\in \RR{\mdim}$ is a bounded additive noise.
The largest singular value of $\M$ is $\sigma_\M$ and its $\OM$-RIP constant is $\delta_\cosp$.
The analysis operator $\OM \in \RR{\pdim \times \sdim}$ is given and fixed.
A procedure $\hat{\CF}_\cosp$ for finding a cosupport that implies
a near optimal projection with a constant $C_\cosp$ is assumed to be at hand.
Our task is to recover $\x$ from $\y$. The recovery result is denoted by $\hat\x$.
\end{defn}


\begin{thm}[Stable Recovery of AIHT and AHTP]
\label{thm:AIHT_AHTP_general_bound}
Consider the problem $\cal P$ and apply
either AIHT or AHTP with a certain constant step-size or
an optimal changing step-size, obtaining $\hat \x^t$ after $t$ iterations. If
\begin{eqnarray}
\label{eq:AIHT_AHTP_general_bound_cond}
\frac{(C_\cosp - 1) \sigma_{\M}^2}{C_\cosp} < 1
\end{eqnarray}
and $$\deltaB< \delta_1(C_\cosp,\sigma^2_\M),$$
where $\delta_1(C_\cosp,\sigma^2_\M)$ is a constant guaranteed to be greater than zero whenever
\eqref{eq:AIHT_AHTP_general_bound_cond} is satisfied and $C_\cosp$ is the near-optimal projection constant for cosparsity $\cosp$ (Definition~\ref{def:C_optimal_proj}),
then after a finite number of iterations $t^*$
\begin{eqnarray}
\norm{{\x} - {\hat\x}^{t^*}}_2 \le  c_1\norm{\vect{e}}_2,
\end{eqnarray}
implying that these algorithms lead to a stable recovery.
The constant $c_{1}$ is a function of $\deltaB$, $C_\cosp$ and $\sigma_\M^2$,
and the constant step-size used is dependent on $\delta_1(C_\cosp,\sigma^2_\M)$.
\end{thm}

\begin{thm}[Stable Recovery of ASP and ACoSaMP]
\label{thm:ACoSaMP_ASP_general_bound}
Consider the problem $\cal P$ and apply
either ACoSaMP or ASP with $a = \frac{2\cosp - \pdim}{\cosp}$,
obtaining $\hat \x^t$ after $t$ iterations. If
 \begin{eqnarray}
\label{eq:ACoSaMP_ASP_general_bound_cond}
\frac{(C_{\hat\CF}^2-1)\sigma_{\matr{M}}^2 }{C_{\hat\CF}^2} < 1,
\end{eqnarray}
and
\begin{eqnarray*}
\deltaD< \delta_2(C_{\hat\CF},\sigma^2_\M),
\end{eqnarray*}
where $C_{\hat\CF} = \max(C_\cosp, C_{2\cosp-\pdim})$ and
$\delta_2(C_{\hat\CF},\sigma^2_\M)$ is a constant guaranteed to be greater than zero whenever
\eqref{eq:ACoSaMP_ASP_general_bound_cond} is satisfied,
then after a finite number of iterations $t^*$
\begin{eqnarray}
 \norm{\x - \hat{\vect{x}}^{t^*} }_2 \le c_2\norm{\vect{e}}_2,
\end{eqnarray}
implying that these algorithms lead to a stable recovery.
The constant $c_{2}$ is a function of $\deltaD$, $C_\cosp$, $C_{2\cosp-\pdim}$ and $\sigma_\M^2$.
\end{thm}

Before we proceed to the proofs, let us comment on the constants in the above theorems.
Their values can be calculated using Theorem~\ref{thm:AIHT_AHTP_theorem_noisy_iter}, and Corollaries~\ref{cor:ACoSaMP_bound} and \ref{cor:ASP_bound}.
In the case where $\OM$ is a unitary matrix, \eqref{eq:AIHT_AHTP_general_bound_cond} and \eqref{eq:ACoSaMP_ASP_general_bound_cond} are trivially satisfied
since $C_{\cosp}=C_{2\cosp - \pdim} = 1$.
In this case the $\OM$-RIP conditions become $\deltaB < \delta_1(1,\sigma^2_\M) =  1/3$
for AIHT and AHTP, and $\deltaD < \delta_2(1,\sigma^2_\M) = 0.0156$ for ACoSaMP and ASP.
In terms of synthesis RIP for $\M\OM^{*}$, the condition $\deltaB < 1/3$ parallels $\delta_{2k}(\M\OM^*) < 1/3$
and similarly $\deltaD < 0.0156$ parallels $\delta_{4k}(\M\OM^*) < 0.0156$.
Note that the condition we pose for AIHT and AHTP in this case
is the same as the one presented for synthesis IHT with a constant step size \cite{Garg09Gradient}.
Better reference constants were achieved in the synthesis case for all four algorithms
and thus we believe that there is still room for improvement of the reference constants in the analysis context.

In the non-unitary case, the value of $\sigma_\M$ plays a vital role, though we believe that this is just an artifact of our proof technique.
For a random Gaussian matrix whose entries are i.i.d with a zero-mean and a variance $\frac{1}{\mdim}$,
$\sigma_\M$ behaves like $\frac{d}{m}\left(1+\sqrt{\frac{\sdim}{\mdim}}\right)$.
This is true also for other types of distributions for which the fourth moment is known to be bounded \cite{Rudelson10Non}.
For example, for $\sdim/\mdim = 1.5$ we have found empirically that $\sigma_\M^2 \simeq 5$.
In this case we need $C_\cosp \le \frac{5}{4}$ for \eqref{eq:AIHT_AHTP_general_bound_cond} to hold
and $C_{\hat\CF} \le  1.118$ for \eqref{eq:ACoSaMP_ASP_general_bound_cond} to hold,
and both are quite demanding on the quality of the near-optimal projection.
For $C_\cosp = C_{\hat\CF} = 1.05$ we have the conditions $\deltaB \le 0.289 $ for AIHT and AHTP, and $\deltaD \le 0.0049$ for ACoSaMP and ASP;
and for $C_\cosp = C_{\hat\CF} = 1.1$ we have $\deltaB \le 0.24 $ for AIHT and AHTP, and $\deltaD \le 0.00032$ for ACoSaMP and ASP.

As in the synthesis case, the $\OM$-RIP requirements for the theoretical bounds of AIHT and AHTP are better than those for ACoSaMP and ASP.
In addition, in the migration from the synthesis to the analysis
we lost more precision in the bounds for ACoSaMP and ASP than in those of AIHT and AHTP.
In particular, even in the case where $\OM$ is the identity
we do not coincide with any of the synthesis parallel RIP reference constants.
We should also remember that the synthesis bound for SP is in terms of $\delta_{3k}$ and not $\delta_{4k}$ \cite{Dai09Subspace}.
Thus, we expect that it will be possible to give a condition for ASP in terms of $\deltaC$
with better reference constants.
However, our main interest in this work is to show the existence of such bounds,
and in Section~\ref{sec:thm_conds} we dwell more on their meaning.

We should note that here and elsewhere we can replace the conditions on $\deltaB$ and $\deltaD$
in the theorems to conditions on $\delta_{2r-p}^{\text{corank}}$ and $\delta_{4r-3p}^{\text{corank}}$
and the proofs will be almost the same\footnote{At a first glance one would think that the conditions should be in terms of $\delta_{2r-d}^{\text{corank}}$ and $\delta_{4r-3d}^{\text{corank}}$.
However, given two cosparse vectors with coranks $r_1$ and $r_2$ the best estimation we can have for the corank of their sum is $r_1+r_2-\pdim$.}.
In this case we will be analyzing a version of the algorithms which is driven by the corank instead of the cosparsity.
This would mean we need the near-optimal projection to be in terms of the corank.
In the case where $\matr{\Omega}$ is in a general position, there is no difference between the cosparsity $\cosp$ and the corank $r$.
However, when we have linear dependencies in $\OM$ the two measures differ and an $\cosp$-cosparse vector
is not necessarily a vector with a corank $r$.

As we will see hereafter, our recovery conditions require $\deltaB$ and $\deltaD$ to be as small as possible
and for this we need $2\cosp -\pdim$ and $4\cosp-3\pdim$ to be as large as possible.
Thus, we need $\cosp$ to be as close as possible to $\pdim$ and for highly redundant $\OM$
this cannot be achieved without having linear dependencies in $\OM$.
Apart from the theoretical advantage of linear dependencies in $\OM$,
we also show empirically that an analysis dictionary with linear dependencies
has better recovery rate than analysis dictionary in a general position of the same dimension.
Thus, we deduce that linear dependencies in $\OM$ lead to better bounds and restoration performance.

Though linear dependencies allow $\cosp$ to be larger than $d$ and be in the order of $\pdim$, the value of the
corank is always bounded by $d$ and cannot be expected to be large enough for highly redundant analysis dictionaries.
In addition, we will see hereafter that the number of measurements $m$ required by the $\OM$-RIP
is strongly dependent on $\cosp$ and less effected by the value of $r$.
From the computational point of view we note also that using corank requires its computation in each iteration
which increases the overall complexity of the algorithms.
Thus, it is more reasonable to have conditions on $\deltaB$ and $\deltaD$ than on $\delta_{2r-p}^{\text{corank}}$ and
$\delta_{4r-3p}^{\text{corank}}$, and our study will be focused on the cosparsity based algorithms.

\subsection{AIHT and AHTP Guarantees}
\label{sec:AIHT_AHTP_guarantees}

A uniform guarantee for AIHT in the case that an optimal projection is given, is presented in \cite{Blumensath09Sampling}.
The work in \cite{Blumensath09Sampling} dealt with a general union of subspaces, $\A$, and assumed that
$\M$ is bi-Lipschitz on the considered union of subspaces.
In our case $\A = \A_\cosp$ and the bi-Lipschitz constants of $\M$ are the largest $B_{L}$ and smallest $B_{U}$
where $0<B_{L} \le B_{U}$ such that for all $\cosp$-cosparse vectors $\vect{v}_1,\vect{v}_2$:
\begin{eqnarray}
\label{eq:ORIP}
B_{L} \norm{\vect{v}_1 + \vect{v}_2}_2^2 \le \norm{\M(\vect{v}_1 + \vect{v}_2)}_2^2 \le B_{U} \norm{\vect{v}_1+\vect{v}_2}_2^2.
\end{eqnarray}
Under this assumption, one can apply Theorem~2 from \cite{Blumensath09Sampling}
to the idealized AIHT that has access to an optimal projection and uses a constant step size $\mu^t = \mu$.
Relying on Table~\ref{tbl:Synthesis_Analysis_parallels} we present this theorem
and replace $B_{L}$ and $B_{U}$ with $1-\deltaB$ and $1+\deltaB$ respectively.
\begin{thm}[Theorem~2 in \cite{Blumensath09Sampling}]\label{thm:AIHT_optimal_case_theorem}
Consider the problem $\cal P$ with $C_\cosp = 1$ and apply
AIHT with a constant step size $\mu$.
If $1+\deltaB \le \frac{1}{\mu} < 1.5(1-\deltaB)$
then after a finite number of iterations $t^*$
\begin{eqnarray}
\norm{{\x}-{\hat\x}^{t^*}}_{2} \le c_3\norm{\mathbf{e}}_{2},
\end{eqnarray}
implying that AIHT leads to a stable recovery.
The constant $c_3$ is a function of $\deltaB$ and $\mu$.
\end{thm}

In this work we extend the above in several ways: First, we refer to the case where optimal projection is not known,
and show that the same flavor guarantees apply for a near-optimal projection\footnote{Remark that we even improve the condition of the idealized case in \cite{Blumensath09Sampling} to be $\deltaB \le \frac{1}{3}$ instead of $\deltaB \le \frac{1}{5}$.}.
The price we seemingly have to pay is that $\sigma_\M$ enters the game.
Second, we derive similar results for the AHTP method.
Finally, we also consider the optimal step size and show that the same performance guarantees hold true in that case.


\begin{thm}
\label{thm:AIHT_AHTP_theorem_noisy_iter}
Consider the problem $\cal P$ and
apply either AIHT or AHTP with a constant step size $\mu$ or an optimal changing step size.
For a positive constant $\eta > 0$,
let
\[
b_1 := \frac{\eta}{1+\eta} \quad\text{and}\quad b_2 := \frac{(C_\ell - 1)\sigma_\M^2 b_1^2}{C_\ell(1-\deltaB)}.
\]
Suppose $\frac{b_2}{b_1^2} = \frac{(C_\ell - 1)\sigma_\M^2}{C_\ell(1-\deltaB)} < 1$,
$1+\deltaB \le \frac{1}{\mu} < \left(1+\sqrt{1-\frac{b_2}{b_1^2}}\right)b_1(1-\deltaB)$ and $\frac{1}{\mu} \le \sigma_\M^2$.
Then for
\begin{eqnarray}
\label{eq:AIHT_AHTP_iter_num}
t \ge t^* \triangleq \frac{\log\left(\frac{\eta\norm{\e}_2^2}{\norm{\y}_2^2}\right)}{\log\left((1+\frac{1}{\eta})^2(\frac{1}{\mu(1-\deltaB)} - 1)C_\cosp + (C_\cosp -1)(\mu\sigma_\M^2 -1) + \frac{C_\cosp}{\eta^2}\right)},
\end{eqnarray}
\begin{eqnarray}
\label{eq:AIHT_AHTP_iter_bound}
\norm{\x-\hat\x^t}_2^2 \le \frac{(1+\eta)^2}{1-\deltaB}\norm{\e}_2^2,
\end{eqnarray}
implying that AIHT and AHTP lead to a stable recovery.
Note that for an optimal changing step-size we set $\mu = \frac{1}{1+\deltaB}$ in $t^*$
and the theorem conditions turn to be $\frac{b_2}{b_1^2} < 1$ and $1+\deltaB <(1+\sqrt{1-\frac{b_2}{b_1^2}})b_1(1-\deltaB)$.
\end{thm}

This theorem is the parallel to Theorems~2.1 in \cite{Garg09Gradient} for IHT.
A few remarks are in order for the nature of the theorem, especially in regards to the constant $\eta$. One can view that $\eta$ gives a trade-off between satisfying the theorem conditions and the amplification of the noise. In particular, one may consider that the above theorem proves the convergence result for the noiseless case by taking $\eta$ to infinity; one can imagine solving the problem $\cal P$ where $\e \to 0$, and applying the theorem with appropriately chosen $\eta$ which approaches infinity. It is indeed possible to show that the iterate solutions of AIHT and AHTP converges to $\x$ when there is no noise. However, we will not give a separate proof since the basic idea of the arguments is the same for both cases.

As to the minimal number of iterations $t^*$ given in \eqref{eq:AIHT_AHTP_iter_num}, one may ask whether it can be negative. In order to answer this question it should be noted that according to the conditions of the Theorem the term inside the log in the denominator \eqref{eq:AIHT_AHTP_iter_num} is always greater than zero. Thus, $t^*$ will be negative only if $\norm{\vect{y}}_2^2 < \eta\norm{\vect{e}}_2^2$. Indeed, in this case $0$ iterations suffice for having the bound in \eqref{eq:AIHT_AHTP_iter_bound}.

The last remark is on the step-size selection. The advantage of the optimal changing step-size over the constant step-size is that we get the guarantee
of the optimal constant step-size $\mu = \frac{1}{1+\deltaB}$ without computing it. This is important since in practice we cannot evaluate the value of $\deltaB$.
However, the disadvantage of using the optimal changing step-size is its additional complexity for the algorithm.
Thus, one option is to approximate the optimal selection rule by replacing it with an adaptive one,
for which we do not have a theoretical guarantee. Another option is to set $\mu = 6/5$ which meets the theorem conditions for small enough $\deltaB$,
in the case where an optimal projection is at hand.

We will prove the theorem by proving two key lemmas first. The proof technique is based on ideas
from \cite{Garg09Gradient} and \cite{Blumensath09Sampling}.
Recall that the two iterative algorithms try to reduce the objective
$\norm{\y-\M\hat\x^t}_2^2$ over iterations $t$.
Thus, the progress of the algorithms can be indirectly measured by how much
the objective $\norm{\y-\M\hat\x^t}_2^2$ is reduced at each iteration $t$.
The two lemmas that we present capture this idea.
The first lemma is similar to Lemma~3 in \cite{Blumensath09Sampling}
and relates $\norm{\y-\M\hat\x^t}_2^2$ to $\norm{\y-\M\hat\x^{t-1}}_2^2$ and similar quantities at iteration $t-1$.
We remark that the constraint $\frac{1}{\mu} \le \sigma_\M^2$ in Theorem~\ref{thm:AIHT_AHTP_theorem_noisy_iter} may not be necessary and is added only for having a simpler derivation of the results in this theorem.
Furthermore, this is a very mild condition compared to $\frac{1}{\mu} < \left(1+\sqrt{1-\frac{b_2}{b_1^2}}\right)b_1(1-\deltaB)$
and can only limit the range of values that can be used with the constant step size
versions of the algorithms.

\begin{lem}
\label{lem:AIHT_lemma}
Consider the problem $\cal P$ and
apply either AIHT or AHTP with a constant step size $\mu$ satisfying
$\frac{1}{\mu} \ge 1+\delta_{2\cosp - \pdim}$ or an optimal step size.
Then, at the $t$-th iteration, the following holds:
\begin{eqnarray}
\label{eq:AIHT_lemma}
&&\norm{\y-\M\hat\x^t}_2^2 - \norm{\y-\M\hat\x^{t-1}}_2^2
\le C_\ell\left(\norm{\y-\M\x}_2^2 - \norm{\y-\M\hat\x^{t-1}}_2^2\right)\\
&& \nonumber ~~~~    + C_\ell \left(\frac{1}{\mu(1-\deltaB)}-1\right)\norm{\M(\x-\hat\x^{t-1})}_2^2
    + (C_\ell-1)\mu\sigma_\M^2 \norm{\y-\M\hat\x^{t-1}}_2^2.
\end{eqnarray}
For the optimal step size the bound is achieved with the value $\mu = \frac{1}{1+\deltaB}$.
\end{lem}

The proof of the above lemma appears in \ref{sec:AIHT_lemma_proof}.
The second lemma is built on the result of Lemma~\ref{lem:AIHT_lemma}.
It shows that once the objective $\norm{\y-\M\hat\x^{t-1}}_2^2$
at iteration $t-1$ is small enough, then we are guaranteed
to have small $\norm{\y-\M\hat\x^{t}}_2^2$ as well.
Given the presence of noise, this is quite natural; one cannot expect
it to approach $0$ but may expect it not to become worse.
Moreover, the lemma also shows that if $\norm{\y-\M\hat\x^{t-1}}_2^2$ is not small,
then the objective in iteration $t$ is necessarily reduced by a constant factor.

\begin{lem}
\label{lem:AIHT_AHTP_lemma_noisy}
Suppose that the same conditions of Theorem~\ref{thm:AIHT_AHTP_theorem_noisy_iter} hold true. If $\norm{\y-\M\hat\x^{t-1}}_2^2 \le \eta^2\norm{\e}_2^2$, then $\norm{\y-\M\hat\x^t}_2^2 \le \eta^2\norm{\e}_2^2$.
Furthermore, if $\norm{\y-\M\hat\x^{t-1}}_2^2 > \eta^2\norm{\e}_2^2$, then
\begin{eqnarray}
\label{eq:AIHT_AHTP_lemma_noisy}
\norm{\y-\M\hat\x^t}_2^2 \le c_4 \norm{\y-\M\x^{t-1}}_2^2
\end{eqnarray}
where
\[
c_4 := \left(1+\frac{1}{\eta}\right)^2\left(\frac{1}{\mu(1-\deltaB)} - 1\right)C_\cosp + (C_\ell -1)(\mu\sigma_\M^2 -1) + \frac{C_\cosp}{\eta^2} < 1.
\]
\end{lem}

Having the two lemmas above, the proof of the theorem is straightforward.

{\em Proof:}[Proof of Theorem~\ref{thm:AIHT_AHTP_theorem_noisy_iter}]
When we initialize $\hat\x^0 = \mathbf{0}$, we have $\norm{\y-\M\hat\x^0}_2^2 = \norm{\y}_2^2$.
Assuming that $\norm{\y}_2 > \eta\norm{\e}_2$ and applying Lemma~\ref{lem:AIHT_AHTP_lemma_noisy} repeatedly, we obtain
\[
\norm{\y-\M\hat\x^t}_2^2 \le \max(c_4^t\norm{\y}_2^2, \eta^2\norm{\e}_2^2).
\]
Since $c_4^t\norm{\y}_2^2 \le \eta^2\norm{\e}_2^2$ for $t \ge t^*$, we have simply
\begin{eqnarray}
\label{eq:AIHT_AHTP_theorem_noisy_iter_step1}
\norm{\y-\M\hat\x^t}_2^2 \le \eta^2\norm{\e}_2^2
\end{eqnarray}
for $t \ge t^*$. If $\norm{\y-\M\hat\x^0}_2 = \norm{\y}_2 \le \eta\norm{\e}_2$ then according to Lemma~\ref{lem:AIHT_AHTP_lemma_noisy},
$\eqref{eq:AIHT_AHTP_theorem_noisy_iter_step1}$ holds for every $t > 0$.
Finally, we observe
\begin{eqnarray}
\label{eq:AIHT_AHTP_theorem_noisy_iter_step2}
\norm{\x-\hat\x^t}_2^2 \le \frac{1}{1-\deltaB}\norm{\M(\x-\hat\x^t)}_2^2
\end{eqnarray}
and, by the triangle inequality,
\begin{eqnarray}
\label{eq:AIHT_AHTP_theorem_noisy_iter_step3}
\norm{\M(\x-\hat\x^t)}_2 \le \norm{\y-\M\hat\x^t}_2 + \norm{\e}_2.
\end{eqnarray}
By plugging \eqref{eq:AIHT_AHTP_theorem_noisy_iter_step1} into \eqref{eq:AIHT_AHTP_theorem_noisy_iter_step3}
and then the resulted inequality into \eqref{eq:AIHT_AHTP_theorem_noisy_iter_step2}, the result of the Theorem follows.
\hfill $\Box$
\bigskip

As we have seen, the above AIHT and AHTP results hold for the cases of
using a constant or an optimal changing step size.
The advantage of using an optimal one is that
we do not need to find $\mu$ that satisfies the conditions of the theorem --
the knowledge that such a $\mu$ exists is enough.
However, its disadvantage is the additional computational complexity it introduces.
In Section~\ref{sec:analysis_alg} we have introduced a third option of using an approximated adaptive step size.
In the next section we shall demonstrate this option in simulations,
showing that it leads to the same reconstruction result as the optimal selection method.
Note, however, that our theoretical guarantees do not cover this case.

\subsection{ACoSaMP Guarantees}

Having the results for AIHT and AHTP we turn to ACoSaMP and ASP.
We start with a theorem for ACoSaMP. Its proof is based on the proof for CoSaMP in \cite{foucart10Sparse}.

\begin{thm}
\label{thm:ACoSaMP_iter_bound}
Consider the problem $\cal P$ and apply ACoSaMP with $a = \frac{2\cosp - \pdim}{\cosp}$.
Let $C_{\hat\CF} = \max({C_\cosp},C_{2\cosp-\pdim})$ and suppose
that there exists $\gamma >0$ such that
 \begin{eqnarray}
\label{eq:C_l_tilda_C_2lp_cond}
(1 + C_{\hat\CF})\left(1   -\bigg(\frac{C_{\hat\CF}}{(1+\gamma)^2}  -(C_{\hat\CF}-1)\sigma_{\matr{M}}^2 \bigg)\right) < 1.
\end{eqnarray}
Then, there exists $\delta_{\text{\tiny ACoSaMP}}(C_{\hat\CF},\sigma_{\matr{M}}^2,\gamma) > 0$
such that, whenever $\deltaD \le \delta_{\text{\tiny ACoSaMP}}(C_{\hat\CF},\sigma_{\matr{M}}^2,\gamma)$,
the $t$-th iteration of the algorithm satisfies
\begin{eqnarray}
\label{eq:ACoSaMP_iter_bound}
&& \hspace{-0.3in} \norm{\x-\hat\x^t}_2 \le \rho_1\rho_2\norm{\vect{x} - \hat\x^{t-1}}_2 + \left(\eta_1 + \rho_1\eta_2 \right)\norm{\vect{e}}_2,
\end{eqnarray}
where
$$\eta_1 \triangleq \frac{\sqrt{\frac{2+C_\cosp}{1+C_\cosp}+2\sqrt{C_\cosp}+C_\cosp}\sqrt{1+\deltaC}}{1-\deltaD},$$
$$\eta_2^2 \triangleq \bigg(\frac{1+\deltaC}{\gamma(1+\alpha)}
 +\frac{(1+\deltaB)C_{2\cosp-\pdim}}{\gamma(1+\alpha)(1+\gamma)}
 + \frac{(C_{2\cosp-\pdim}-1)(1+\gamma)\sigma_{\matr{M}}^2}{(1+\alpha)(1+\gamma)\gamma}\bigg),$$
 $$\rho_1^2 \triangleq \frac{1+2\deltaD\sqrt{C_\cosp} + C_\cosp}{1-\deltaD^2},$$
  $$\small {\rho_2^2 \triangleq 1   -\bigg(\sqrt{\deltaD} -
 \sqrt{\frac{C_{2\cosp-\pdim}}{(1+\gamma)^2}\left(1-\sqrt{\deltaB}\right)^2
 -(C_{2\cosp-\pdim}-1)(1+\deltaB)\sigma_{\matr{M}}^2} \bigg)^2}$$
and
$$ \small {\alpha = \frac{\sqrt{\delta_{4\cosp-3p}}}{\sqrt{\frac{C_{2\cosp-p}}{(1+\gamma)^2}\left(1-\sqrt{\delta_{2\cosp-p}}\right)^2 -(C_{2\cosp-p}-1)(1+\delta_{2\cosp-p})\sigma_{\matr{M}}^2} - \sqrt{\delta_{4\cosp-3p}}}}.$$
Moreover, $\rho_1^2\rho_2^2 < 1$, i.e., the iterates converges.
\end{thm}

The constant $\gamma$ plays a similar role to the constant $\eta$ of Theorem~\ref{thm:AIHT_AHTP_theorem_noisy_iter}. It gives a tradeoff between satisfying the theorem conditions and the noise amplification. However, as opposed to $\eta$, the conditions for the noiseless case are achieved when $\gamma$ tends to zero. An immediate corollary of the above theorem is the following.
\begin{cor}
\label{cor:ACoSaMP_bound}
Consider the problem $\cal P$ and
apply ACoSaMP with $a = \frac{2\cosp - \pdim}{\cosp}$.
If \eqref{eq:C_l_tilda_C_2lp_cond} holds and
$\deltaD< \delta_{\text{\tiny ACoSaMP}}(C_{\hat\CF},\sigma_{\matr{M}}^2,\gamma)$,
where $C_{\hat\CF}$ and $\gamma$ are as in Theorem~\ref{thm:ACoSaMP_iter_bound} and
$\delta_{\text{\tiny ACoSaMP}}(C_{\hat\CF},\sigma_{\matr{M}}^2,\gamma)$ is a constant guaranteed to be greater than zero whenever
\eqref{eq:ACoSaMP_ASP_general_bound_cond} is satisfied,
then for any $$t \ge  t^* = \ceil{\frac{\log(\norm{\vect{x}}_2/\norm{\vect{e}}_2)}{\log(1/\rho_1\rho_2)}},$$
\begin{eqnarray}
\label{eq:ACoSaMP_bound}
&& \hspace{-0.5in} \norm{\vect{x} - \hat\x^{t^*}_{\text{\tiny ACoSaMP}}}_2 \le
  \left(1 +  \frac{1-(\rho_1\rho_2)^{t^*}}{1-\rho_1\rho_2}\left(\eta_1 + \rho_1\eta_2 \right)\right)\norm{\vect{e}}_2,
\end{eqnarray}
implying that ACoSaMP leads to a stable recovery.
The constants $\eta_1$, $\eta_2$, $\rho_1$ and $\rho_2$ are the same as in Theorem~\ref{thm:ACoSaMP_iter_bound}.
\end{cor}
{\em Proof:}
By using \eqref{eq:ACoSaMP_iter_bound} and recursion we have that
\begin{eqnarray}
\label{eq:ACoSaMP_bound_step1}
&& \hspace{-0.3in} \norm{\vect{x} - \hat\x^{t^*}_{\text{\tiny ACoSaMP}}}_2 \le (\rho_1\rho_2)^{t^*}\norm{\vect{x} - \hat\x^{0}_{\text{\tiny ACoSaMP}}}_2
\\ \nonumber && \hspace{-0.3in}
  +(1+\rho_1\rho_2+(\rho_1\rho_2)^2+\dots (\rho_1\rho_2)^{t^*-1})\left(\eta_1 + \rho_1\eta_2 \right)\norm{\vect{e}}_2.
\end{eqnarray}
Since $\hat\x^{0}_{\text{\tiny ACoSaMP}}=0$, after $t^*$ iterations, one has
\begin{eqnarray}
\label{eq:ACoSaMP_bound_step2}
(\rho_1\rho_2)^{t^*}\norm{\vect{x} - \hat\x^{0}_{\text{\tiny ACoSaMP}}}_2 = (\rho_1\rho_2)^{t^*}\norm{\x}_2 \le \norm{\e}_2.
\end{eqnarray}
By using the equation of geometric series with \eqref{eq:ACoSaMP_bound_step1} and plugging \eqref{eq:ACoSaMP_bound_step2} into it,
we get \eqref{eq:ACoSaMP_bound}.
\hfill $\Box$ \bigskip

We turn now to prove the theorem. Instead of presenting the proof directly, we divide the proof into several lemmas.
The first lemma gives a bound for $\norm{\x - \vect{w}}_2$ as a function of $\norm{\vect{e}}_2$ and $\norm{\P_{\tilde{\Lambda}^t}(\x - \vect{w})}_2$.
\begin{lem}
\label{lem:ACoSaMP_xp_bound}
Consider the problem $\cal P$ and
apply ACoSaMP with $a = \frac{2\cosp - \pdim}{\cosp}$.
For each iteration we have
\begin{eqnarray}
\label{eq:ACoSaMP_xp_bound}
\norm{\x - \vect{w}}_2 &\le& \frac{1}{\sqrt{1-\deltaD^2}}\norm{\P_{\tilde{\Lambda}^t}(\vect{x} - \vect{w} )}_2+  \frac{\sqrt{1+\deltaC}}{1-\deltaD}\norm{\vect{e}}_2.
\end{eqnarray}
\end{lem}

The second lemma bounds $\norm{\x - \hat\x^{t}_{\text{\tiny ACoSaMP}}}_2$ in terms of $\norm{\P_{\tilde{\Lambda}^t}(\x- \hat\x^{t}_{\text{\tiny ACoSaMP}} )}_2$ and $\norm{\vect{e}}_2$ using the first lemma.
\begin{lem}
\label{lem:ACoSaMP_xt_bound1}
Consider the problem $\cal P$ and
apply ACoSaMP with $a = \frac{2\cosp - \pdim}{\cosp}$.
For each iteration we have
\begin{eqnarray}
&& \hspace{-0.5in} \norm{\x - \hat\x^t}_2 \le \rho_1\norm{\P_{\tilde{\Lambda}^t}(\x - \vect{w})}_2 +  \eta_1\norm{\vect{e}}_2,
\end{eqnarray}
where $\eta_1$ and $\rho_1$ are the same constants as in Theorem~\ref{thm:ACoSaMP_iter_bound}.
\end{lem}

The last lemma bounds $\norm{\P_{\tilde{\Lambda}^t}(\x - \vect{w})}_2$ with $\norm{\x - \hat\x^{t-1}_{\text{\tiny ACoSaMP}}}_2$ and $\norm{\vect{e}}_2$.
\begin{lem}
\label{lem:ACoSaMP_Pxp_bound}
Consider the problem $\cal P$ and
apply ACoSaMP with $a = \frac{2\cosp - \pdim}{\cosp}$.
if
\begin{eqnarray}
\label{eq:C_2lp_cond}
C_{2\cosp-p} < \frac{\sigma_{\matr{M}}^2(1+\gamma)^2}{\sigma_{\matr{M}}^2(1+\gamma)^2-1},
\end{eqnarray}
then there exists $\tilde{\delta}_{\text{ \tiny ACoSaMP}}(C_{2\cosp - \pdim}, \sigma_{\matr{M}}^2, \gamma) >0$ such that
for any $\deltaB < \tilde{\delta}_{\text{ \tiny ACoSaMP}}(C_{2\cosp - \pdim}, \sigma_{\matr{M}}^2, \gamma)$
\begin{eqnarray}
\label{eq:ACoSaMP_Pxp_bound}
&& \norm{\P_{\tilde{\Lambda}^t}(\x - \vect{w})}_2 \le
    \eta_2\norm{\vect{e}}_2
 + \rho_2\norm{\x - \hat\x^{t-1}}_2.
\end{eqnarray}
The constants $\eta_2$ and $\rho_2$ are as defined in Theorem~\ref{thm:ACoSaMP_iter_bound}.
\end{lem}

The proofs of Lemmas~\ref{lem:ACoSaMP_xp_bound}, \ref{lem:ACoSaMP_xt_bound1} and \ref{lem:ACoSaMP_Pxp_bound} appear in \ref{sec:ACoSaMP_xp_bound_proof}, \ref{sec:ACoSaMP_xt_bound1_proof} and \ref{sec:ACoSaMP_Pxp_bound_proof} respectively.
With the aid of the above three lemmas we turn to prove Theorem~\ref{thm:ACoSaMP_iter_bound}.

{\em Proof:}[Proof of Theorem~\ref{thm:ACoSaMP_iter_bound}]
Remark that since $1+C_{\hat{\CF}} > 1$ we have that \eqref{eq:C_l_tilda_C_2lp_cond} implies $\frac{C_{\hat\CF}}{(1+\gamma)^2}  -(C_{\hat\CF}-1)\sigma_{\matr{M}}^2 \ge 0$. Because of that
the condition in \eqref{eq:C_2lp_cond} in Lemma~\ref{lem:ACoSaMP_Pxp_bound} holds.
Substituting the inequality of Lemma~\ref{lem:ACoSaMP_Pxp_bound} into the inequality of Lemma~\ref{lem:ACoSaMP_xt_bound1} gives
\eqref{eq:ACoSaMP_iter_bound}. The iterates convergence if $\rho_1^2\rho_2^2 = \frac{1+2\deltaD\sqrt{C_\cosp} + C_\cosp}{1 -\deltaD^2}\rho_2^2 < 1$.
By noticing that $\rho_2^2<1$ it is enough to require $\frac{1 + C_\cosp}{1 -\deltaD^2}\rho_2^2 + \frac{2\deltaD\sqrt{C_\cosp}}{1 -\deltaD^2} < 1$. The last is equivalent to
\begin{eqnarray}
\label{eq:ACoSaMP_rho1rho2_cond_step1}
&& \hspace{-0.4in} (1 + C_\cosp)
   \left(1   -\left(\sqrt{\deltaD} -
 \sqrt{\frac{C_{2\cosp - \pdim}}{(1+\gamma)^2}\left(1-\sqrt{\deltaB}\right)^2
 -(C_{2\cosp - \pdim}-1)(1+\deltaB)\sigma_{\matr{M}}^2} \right)^2\right) \\ \nonumber &&
~~~~~~~~~~~~~~~~~~~~~~~~~~~~~~~~~~~~~~~~~~~~~~~~~~~~~~~~~~~~~~~~~~~~~~~~+2\deltaD\sqrt{C_\cosp} -1 +\deltaD^2  < 0.
\end{eqnarray}
It is easy to verify that $\zeta(C,\delta) \triangleq \frac{C}{(1+\gamma)^2}\left(1-\sqrt{\delta}\right)^2
 -(C-1)(1+\delta)\sigma_{\matr{M}}^2$ is a decreasing function of both $\delta$ and $C$ for $0 \le \delta \le 1$ and $C>1$.
Since $1 \le C_{2\cosp - \pdim}\le C_{\hat\CF}$, $\deltaB \le \deltaD$ and $\delta \ge 0$ we have that
$\zeta(C_{\hat\CF},\deltaD) \le \zeta(C_{2\cosp - \pdim},\deltaD) \le \zeta(C_{2\cosp - \pdim},\deltaB) \le \zeta(1,0) =\frac{1}{(1+\gamma)^2}\le 1$.
Thus we have that $-1 \le-(\sqrt{\deltaD} -\zeta(C_{2\cosp-\pdim},\deltaB))^2 \le -\deltaD +2\sqrt{\deltaD} - \zeta(C_{\hat\CF},\deltaD)$.
Combining this with the fact that $C_\cosp \le C_{\hat\CF}$ provides the following guarantee for $\rho_1^2\rho_2^2 < 1$,
\begin{eqnarray}
\label{eq:ACoSaMP_rho1rho2_cond_step2}
&& \hspace{-0.3in} (1+ C_{\hat\CF})\bigg(1 -\deltaD +2\sqrt{\deltaD}
\\ \nonumber && \hspace{-0.3in}
 -\frac{C_{\hat\CF}}{(1+\gamma)^2}\left(1-2\sqrt{\deltaD}+\deltaD\right)
 +(C_{\hat\CF}-1)(1+\deltaD)\sigma_{\matr{M}}^2 \bigg)+2\deltaD\sqrt{C_{\hat\CF}}  -1 + \deltaD^2 <0.
\end{eqnarray}
Let us now assume that $\deltaD \le \frac{1}{2}$. This necessarily means that $\delta_{\text{\tiny ACoSaMP}} \le \frac{1}{2}$ in the end.
This assumption implies $\deltaD^2\le \frac{1}{2}\deltaD$.
Using this and gathering coefficients, we now consider the condition
\begin{eqnarray}
\label{eq:ACoSaMP_rho1rho2_cond_step3}
&& \hspace{-0.3in} (1+ C_{\hat\CF})\left(1-\frac{C_{\hat\CF}}{(1+\gamma)^2}+(C_{\hat\CF}-1)\sigma_{\matr{M}}^2\right)-1 + 2(1+ C_{\hat\CF})\left(1+\frac{C_{\hat\CF}}{(1+\gamma)^2}\right)\sqrt{\deltaD}
\\ \nonumber && \hspace{-0.3in}
  +\left((1+ C_{\hat\CF})\left(-1-\frac{C_{\hat\CF}}{(1+\gamma)^2}+(C_{\hat\CF}-1)\sigma_{\matr{M}}^2\right)+ 2\sqrt{C_{\hat\CF}} +\frac{1}{2}\right)\deltaD <0.
\end{eqnarray}
The expression on the LHS is a quadratic function of $\sqrt{\deltaD}$.
Note that since \eqref{eq:C_l_tilda_C_2lp_cond} holds the constant term in the quadratic function is negative.
This guarantees the existence of a range of values $[0,\delta_{\text{\tiny ACoSaMP}}(C_{\hat\CF},\sigma_{\matr{M}}^2,\gamma)]$ for $\deltaD$ for which \eqref{eq:ACoSaMP_rho1rho2_cond_step3} holds,
where $\delta_{\text{\tiny ACoSaMP}}(C_{\hat\CF},\sigma_{\matr{M}}^2,\gamma)$ is the square of the positive solution of
the quadratic function. In case of two positive solutions we should take the smallest among them -- in this case
the coefficient of $\deltaD$ in \eqref{eq:ACoSaMP_rho1rho2_cond_step3} will be positive.

Looking back at the proof of the theorem, we observe that the value of the constant $\delta_{\text{\tiny ACoSaMP}}(C_{\hat\CF},\sigma_{\matr{M}}^2,\gamma)$ can
potentially be improved: at the beginning of the proof, we have used $\rho_2^2\le 1$.
By the end, we obtained $\rho_2^2\le \rho_1^{-2} \le 0.25$ since $\rho_1 > 2$. If we were to use this bound at the beginning, we would have obtained better
constant $\delta_{\text{\tiny ACoSaMP}}(C_{\hat\CF},\sigma_{\matr{M}}^2,\gamma)$.
\hfill $\Box$ \bigskip

\subsection{ASP Guarantees}

Having the result of ACoSaMP we turn to derive a similar result for ASP.
The technique for deriving a result for ASP based on the result of ACoSaMP is similar
to the one we used to derive a result for AHTP from the result of AIHT.

\begin{thm}
\label{thm:ASP_iter_bound}
Consider the problem $\cal P$ and
apply ASP with $a = \frac{2\cosp - \pdim}{\cosp}$.
If \eqref{eq:C_l_tilda_C_2lp_cond} holds
and $\deltaD \le
\delta_{\text{\tiny ASP}}(C_{\hat\CF},\sigma_{\matr{M}}^2,\gamma)$,
where  $C_{\hat\CF}$ and
 $\gamma$ are as in Theorem~\ref{thm:ACoSaMP_iter_bound}, and
$\delta_{\text{\tiny ASP}}(C_{\hat\CF},\sigma_{\matr{M}}^2,\gamma)$ is a constant guaranteed to be greater than zero whenever
\eqref{eq:C_l_tilda_C_2lp_cond} is satisfied,
then the $t$-th iteration of the algorithm satisfies
\begin{eqnarray}
\label{eq:ASP_iter_bound}
&& \hspace{-0.5in} \norm{\vect{x} - \hat\x^t_{\text{\tiny ASP}}}_2 \le \frac{1+\deltaB}{1-\deltaB
}\rho_1\rho_2\norm{\vect{x} - \hat{\vect{x}}^{t-1}_{\text{\tiny ASP} }}_2 +\left(\frac{1+\deltaB}{1-\deltaB}\left(\eta_1 + \rho_1\eta_2 \right)+ \frac{2}{1-\deltaB}\right)\norm{\vect{e}}_2.
\end{eqnarray}
and the iterates converges, i.e., $\rho_1^2\rho_2^2 < 1$.
The constants $\eta_1$, $\eta_2$, $\rho_1$ and $\rho_2$ are the same as in Theorem~\ref{thm:ACoSaMP_iter_bound}.
\end{thm}
{\em Proof:}
We first note that according to the selection rule of $\hat{\vect{x}}_{\text{\tiny ASP}}$ we have that
\begin{eqnarray}
\norm{\vect{y} - \matr{M}\hat{\vect{x}}^t_{\text{\tiny ASP}}}_2 \le \norm{\vect{y} - \matr{M}\matr{Q}_{\hat\Lambda^t}\vect{w}}_2.
\end{eqnarray}
Using the triangle inequality and the fact that $\vect{y} = \matr{M}\vect{x}+\e$ for both the LHS and the RHS we have
\begin{eqnarray*}
\norm{\matr{M}(\vect{x} - \hat{\vect{x}}^t_{\text{\tiny ASP}})}_2 -\norm{\vect{e}}_2 \le \norm{ \matr{M}(\vect{x} - \matr{Q}_{\hat\Lambda^t}\vect{w})}_2 + \norm{\vect{e}}_2.
\end{eqnarray*}
Using the $\matr{\Omega}$-RIP property of $\matr{M}$ with the fact that $\vect{x}$, $\hat{\vect{x}}_{\text{\tiny ASP}}$
and $\matr{Q}_{\hat\Lambda^t}\vect{w}$ are $\cosp$-cosparse we have
\begin{eqnarray*}
\norm{\vect{x} - \hat{\vect{x}}^t_{\text{\tiny ASP}}}_2 \le
\frac{1+\deltaB}{1-\deltaB}\norm{\vect{x} - \matr{Q}_{\hat\Lambda^t}\vect{w}}_2
+ \frac{2}{1-\deltaB}\norm{\vect{e}}_2.
\end{eqnarray*}
Noticing that $\matr{Q}_{\hat\Lambda^t}\vect{w}$ is the solution we get in one iteration of ACoSaMP with initialization of $\hat{\vect{x}}^{t-1}_{\text{\tiny ASP} }$, we can combine the above with the result of Theorem~\ref{thm:ACoSaMP_iter_bound} getting \eqref{eq:ASP_iter_bound}.
For $\frac{1+\deltaB}{1-\deltaB}\rho_1\rho_2<1$ to hold we need that
\begin{eqnarray}
\label{eq:ASP_rho1rho2_cond1}
&& \hspace{-0.3in} \frac{1+2\deltaD\sqrt{C_\cosp} + C_\cosp}{(1-\deltaD)^2}
   \left(1-\left( \left(\frac{\sqrt{\tilde{C}_{2\cosp-p}}}{1+\gamma} +1 \right)\sqrt{\deltaD} -\frac{\sqrt{\tilde{C}_{2\cosp-p}}}{1+\gamma}
 \right)^2\right) < 1.
\end{eqnarray}
Remark that the above differs from what we have for ACoSaMP only in the
denominator of the first element in the LHS. In ACoSaMP $1-\deltaD^2$ appears instead of $(1-\deltaD)^2$.
Thus, Using a similar process to the one in the proof of ACoSaMP we can show that \eqref{eq:ASP_rho1rho2_cond1} holds if the following holds
\begin{eqnarray}
\label{eq:ASP_rho1rho2_cond_step3}
&& \hspace{-0.3in} (1+ C_{\hat\CF})\left(1-\frac{C_{\hat\CF}}{(1+\gamma)^2}+(C_{\hat\CF}-1)\sigma_{\matr{M}}^2\right)-1 + 2(1+ C_{\hat\CF})\left(1+\frac{C_{\hat\CF}}{(1+\gamma)^2}\right)\sqrt{\deltaD}
\\ \nonumber && \hspace{-0.3in}
  +\left((1+ C_{\hat\CF})\left(-1-\frac{C_{\hat\CF}}{(1+\gamma)^2}+(C_{\hat\CF}-1)\sigma_{\matr{M}}^2\right)+ 2\sqrt{C_{\hat\CF}} +2\right)\deltaD <0.
\end{eqnarray}
Notice that the only difference of the above compared to \eqref{eq:ACoSaMP_rho1rho2_cond_step3}
is that we have $+2$ instead of $+0.5$ in the coefficient of $\deltaD$ and this is due to the difference we mentioned before in the denominator in \eqref{eq:ASP_rho1rho2_cond1}.
The LHS of \eqref{eq:ASP_rho1rho2_cond_step3} is a quadratic function of $\sqrt{\deltaD}$. As before, we notice that
if \eqref{eq:C_l_tilda_C_2lp_cond} holds then the constant term of the above is positive
and thus $\delta_{\text{\tiny ASP}}(C_{\hat\CF},\sigma_{\matr{M}}^2,\gamma) \ge 0$ exists
 and is the square of the positive solution of the quadratic function.
\hfill $\Box$ \bigskip

Having Theorem~\ref{thm:ASP_iter_bound} we can immediately have the following corollary
which is similar to the one we have for ACoSaMP.
The proof resembles the one of Corollary~\ref{cor:ACoSaMP_bound} and omitted.

\begin{cor}
\label{cor:ASP_bound}
Consider the problem $\cal P$ and
apply ASP with $a = \frac{2\cosp - \pdim}{\cosp}$.
If \eqref{eq:C_l_tilda_C_2lp_cond} holds
and $\deltaD \le
\delta_{\text{\tiny ASP}}(C_{\hat\CF},\sigma_{\matr{M}}^2,\gamma)$,
where  $C_{\hat\CF}$ and
 $\gamma$ are as in Theorem~\ref{thm:ACoSaMP_iter_bound}, and
$\delta_{\text{\tiny ASP}}(C_{\hat\CF},\sigma_{\matr{M}}^2,\gamma)$ is a constant guaranteed to be greater than zero whenever
\eqref{eq:ACoSaMP_ASP_general_bound_cond} is satisfied,
then for any
$$t\ge t^* = \ceil{\frac{\log(\norm{\vect{x}}_2/\norm{\vect{e}}_2)}{\log(1/\frac{1+\deltaB}{1-\deltaB}\rho_1\rho_2)}},$$
\begin{eqnarray}
\label{eq:ASP_bound}
&& \hspace{-0.5in} \norm{\vect{x}^{t}_{\text{\tiny ASP}} -\vect{x}}_2 \le
  \Bigg(1 +  \frac{1-\left(\frac{1+\deltaB}{1-\deltaB}\rho_1\rho_2\right)^{t}}{1-\frac{1+\deltaB}{1-\deltaB}\rho_1\rho_2}\cdot
   \left(\frac{1+\deltaB}{1-\deltaB}\left(\eta_1 + \rho_1\eta_2 \right)+ \frac{2}{1-\deltaB}\right)\Bigg)\norm{\vect{e}}_2.
\end{eqnarray}
implying that ASP leads to a stable recovery.
The constants $\eta_1$, $\eta_2$, $\rho_1$ and $\rho_2$ are the same as in Theorem~\ref{thm:ACoSaMP_iter_bound}.
\end{cor}

\subsection{Non-Exact Cosparse Case}

In the above guarantees we have assumed that the signal $\x$ is $\cosp$-cosparse. In many cases, it is not exactly $\cosp$-cosparse but only nearly so. Denote by $\x^\cosp = \Q_{\CF^*_\cosp(\x)}\x$ the best $\cosp$-cosparse approximation of $\x$, we have the following theorem that provides us with a guarantee also for this case. Similar result exists also in the synthesis case for the synthesis-$\ell_1$ minimization problem \cite{Candes08TheRIP}.

\begin{thm}\label{thm:non_exact_noisy_result}
Consider a variation of problem $\cal P$ where $\x$ is a general vector,
and apply either AIHT or AHTP both with either constant or changing step size; or ACoSaMP or ASP with $a = \frac{2\cosp - \pdim}{\cosp}$, and all are used with a zero initialization.
Under the same conditions of Theorems~\ref{thm:AIHT_AHTP_general_bound} and \ref{thm:ACoSaMP_ASP_general_bound} we
have for any $t\ge t^*$
\begin{eqnarray}
\label{eq:non_exact_bound}
\norm{\x - {\hat\x}}_2 \le \norm{\x - \x^\cosp}_2 + c\norm{\M(\x - \x^\cosp)}_2 + c\norm{\e}_2,
\end{eqnarray}
where $t^*$ and $c$ are the constants from Theorems~\ref{thm:AIHT_AHTP_general_bound} and \ref{thm:ACoSaMP_ASP_general_bound}.
\end{thm}
{\em Proof:}
First we notice that we can rewrite $\y = \M\x^\cosp + \M(\x - \x^\cosp) + \vect{e}$. Denoting $\vect{e}^\cosp = \M(\x - \x^\cosp) + \vect{e}$ we can use
Theorems~\ref{thm:AIHT_AHTP_general_bound} and \ref{thm:ACoSaMP_ASP_general_bound} to recover $\x^\cosp$ and have
\begin{eqnarray}
\norm{\x^\cosp - \hat\x}_2 \le c\norm{\e^\cosp}_2.
\end{eqnarray}
Using the triangle inequality for $\norm{\x - \hat\x}_2$ with the above gives
\begin{eqnarray}
\norm{\x - {\hat\x}}_2 \le \norm{\x - \x^\cosp}_2 +\norm{\x^\cosp - {\hat\x}}_2 \le \norm{\x - \x^\cosp}_2+ c\norm{\e^\cosp}_2.
\end{eqnarray}
Using again the triangle inequality for $\norm{\e^\cosp}_2 \le  \norm{\e}_2 + \norm{\M(\x - \x^\cosp)}_2$ provides us with the desired result.
\hfill $\Box$
\bigskip

\subsection{Theorem Conditions}
\label{sec:thm_conds}

Having the results of the theorems we ask ourselves whether their conditions are feasible.
As we have seen in Section~\ref{sec:omega_RIP}, the requirement on the $\OM$-RIP for many non-trivial matrices.
In addition, as we have seen in the introduction of this section we need $C_\cosp$ and $C_{2\cosp -\pdim}$ to be one or close to one for satisfying the conditions of the theorems.
Using the thresholding in \eqref{eq:thresh_cosupp_selection} for cosupport selection with a unitary $\OM$ satisfies the conditions in a trivial way since $C_\cosp = {C}_{2\cosp-p} = 1$. This case coincides with the synthesis model for which we already have theoretical guarantees. As shown in Section~\ref{sec:near_opt_proj}, optimal projection schemes
exist for $\OM_{\text{1D-DIF}}$ and $\OM_{\text{FUS}}$ which do not belong to the the synthesis framework.
For a general $\matr{\Omega}$, a general projection scheme is not known and if the thresholding method is used the constants in \eqref{eq:thresh_cosupp_selection} do not equal one
and are not even expected to be close to one \cite{Giryes11Iterative}.
It is interesting to ask whether there exists
an efficient general projection scheme that guarantees small constants for any given
operator $\matr{\Omega}$, or for specifically structured  $\OM$. We leave these questions as subject for future work.
Instead, we show empirically in the next section that a weaker projection scheme that
does not fulfill all the requirements of the theorems leads to a good reconstruction result.
This suggests that even in the absence of good near optimal projections we may still
use the algorithms practically.


\subsection{Comparison to Other Works}
\label{sec:comp_other}

Among the existing theoretical works that studied the performance of analysis algorithms
\cite{Nam12Cosparse,Vaiter12Robust,Peleg12Performance},
the result that resembles ours is the result for $\ell_1$-analysis
in \cite{Candes11Compressed}.
This work analyzed the $\ell_1$-analysis minimization problem with a synthesis perspective.
The analysis dictionary $\OM$ was replaced with the conjugate of a synthesis dictionary $\matr{D}$ which is assumed
to be a tight frame, resulting with the following minimization problem.
\begin{eqnarray}
\label{eq:l1_D_analysis}
\hat{\vect{x}}_{A-\ell_1} = \argmin_{\vect{z}} \norm{\matr{D}^*\vect{z}}_1 & s.t. & \norm{\vect{y} - \matr{M}\vect{z}}_2 \le \epsilon.
\end{eqnarray}
It was shown that if $\matr{M}$ has the $\matr{D}$-RIP \cite{Candes11Compressed,Blumensath09Sampling} with $\delta_{7k} < 0.6$, an extension of the synthesis RIP, then
\begin{eqnarray}
\label{eq:l1_candes_rec_error}
\norm{\hat{\vect{x}}_{A-\ell_1} - \vect{x}}_2 \le \tilde{C}_{\ell_1}\epsilon + \frac{\norm{\matr{D}^*\vect{x}-[\matr{D}^*\vect{x}]_k}_1}{\sqrt{k}}.
\end{eqnarray}
We say that a matrix $\matr{M}$ has a $\matr{D}$-RIP with a constant $\delta_k$
if for any signal $\vect{z}$ that has a $k$-sparse representation under $\matr{D}$
\begin{eqnarray}
(1-\delta_k)\norm{{\vect{z}}}_2^2 \le \norm{\matr{M}{\vect{z}}}_2^2 \le (1+\delta_k)\norm{{\vect{z}}}_2^2.
\end{eqnarray}
The authors in \cite{Candes11Compressed}
presented this result as a synthesis result
that allows linear dependencies in $\matr{D}$ at the cost of limiting
the family of signals to be those for which $\norm{\matr{D}^*\vect{x}-[\matr{D}^*\vect{x}]_k}_1$ is small.
However, having the analysis perspective,
we can realize that they provided a recovery guarantee
for $\ell_1$-analysis under the new analysis model for the case that $\OM$ is a tight frame.
An easy way to see it is to observe that for an $\cosp$-cosparse signal $\vect{x}$,
setting $k=\pdim-\cosp$, we have that $\norm{\OM\vect{x}-[\OM^*\vect{x}]_{\pdim-\cosp}}_1 =0$
and thus in the case $\epsilon = 0$ we get that $\eqref{eq:l1_candes_rec_error}$
guarantees the recovery of $\vect{x}$ by using \eqref{eq:l1_D_analysis} with $\matr{D}^* = \OM$.
Thus, though the result in \cite{Candes11Compressed} was presented as
a reconstruction guarantee for the synthesis model,
it is actually a guarantee for the analysis model.

Our main difference from \cite{Candes11Compressed} is that
the proof technique relies on the analysis model and not on the synthesis one
and that the results presented here are for general operators and not only for tight frames.
For instance, the operators $\OM_{\text{1D-DIF}}$ and $\OM_{\text{FUS}}$ for which the guarantees hold are not tight frames
where $\OM_{\text{1D-DIF}}$ is not even a frame.
However, the drawback of our approached compared to the work in \cite{Candes11Compressed} is that
it is still not known how to perform an optimal or a near optimal projection for a tight frame.

In the non-exact sparse case our results differ from the one in \eqref{eq:l1_candes_rec_error}
in the sense that it looks at the projection error and not at the values of $\OM\x$.
It would be interesting to see if there is a connection between the two and whether
one implies the other.

A recent work has studied the $\ell_1$-analysis minimization with the 2D-DIF operator,
also known as anisotropic two dimensional total-variation (2D-TV) \cite{Needell12Stable}.
It would be interesting to see whether similar results can be achieved for the greedy-like techniques proposed here with 2D-DIF.

\section{Experiments}
\label{sec:exp}

In this section we repeat some of the experiments performed in \cite{Nam12Cosparse}
for the noiseless case ($\vect{e} =0$) and some of the experiments performed in \cite{Nam11GAPN} for the noisy case\footnote{
A matlab package with code for the experiments performed in this paper is in preparation for an open source distribution.}.

\subsection{Targeted Cosparsity}
\label{sec:targ_cosp}
Just as in the synthesis counterpart of the proposed algorithms,
where a target sparsity level $\spar$ must be selected before running the algorithms,
we have to choose the targeted cosparsity level which will dictate the projection steps.
In the synthesis case it is known that it may be beneficial to over-estimate the sparsity $\spar$.
Similarly in the analysis framework the question arises:
In terms of recovery performance, does it help to under-estimate the cosparsity $\cosp$?
A tentative positive answer comes from the following heuristic:
Let $\tilde\Lambda$ be a subset of the cosupport $\Lambda$ of the signal $\x$
with $\tilde\cosp := |\tilde\Lambda| < \cosp = |\Lambda|$.
According to Proposition~3 in \cite{Nam12Cosparse}
\begin{eqnarray}
\label{eq:kappa_uniqueness}
\kappa_{\OM}(\tilde{\cosp}) \le \frac{m}{2}
\end{eqnarray}
is a sufficient condition to identify $\tilde\Lambda$
in order to recover $\x$ from the relations $\y = \M\x$ and $\OM_{\tilde\Lambda}\x = 0$.
$\kappa_{\OM}(\tilde{\cosp}) = \max_{\tilde{\Lambda} \in \L_{\tilde\cosp}} \dim (\aspace_{\tilde{\Lambda}})$ is a function of $\tilde\cosp$.
Therefore, we can replace $\cosp$ with the smallest $\tilde\cosp$ that
satisfies \eqref{eq:kappa_uniqueness} as the effective cosparsity in the algorithms.
Since it is easier to identify a smaller cosupport set it is better to run the algorithm
with the smallest possible value of $\tilde\cosp$, in the absence of noise.
In the presence of noise, larger values of $\cosp$ allows a better denoising.
Note, that in some cases the smallest possible value of $\tilde\cosp$
will be larger than the actual cosparsity of $\x$.
In this case we cannot replace $\cosp$ with $\tilde{\cosp}$.

We take two examples for selecting $\tilde{\cosp}$. The first is for $\OM$ which is in general position and the second is for
$\OM_{2D-DIF}$, the finite difference analysis operator that computes horizontal and vertical discrete derivatives
of an image which is strongly connected to the total variation (TV) norm minimization as noted before.
For $\OM$ that is in general position $\kappa_{\OM}(\tilde{\cosp}) = \max(\sdim -\cosp,0)$ \cite{Nam12Cosparse}.
In this case we choose
\begin{eqnarray}
\tilde{\cosp} = \min\left(\sdim -\frac{\mdim}{2}, \cosp\right).
\end{eqnarray}
For $\OM_{DIF}$ we have $\kappa_{\OM_{DIF}}(\tilde{\cosp}) \ge d - \frac{\cosp}{2} - \sqrt{\frac{\cosp}{2}}-1$ \cite{Nam12Cosparse}
and
\begin{eqnarray}
\tilde{\cosp} = \lceil\min((-1/\sqrt{2}+\sqrt{2d-m-1.5})^2, \cosp)\rceil.
\end{eqnarray}

Replacing $\cosp$ with $\tilde\cosp$ is more relevant to AIHT and AHTP than ACoSaMP and ASP
since in the last we intersect cosupport sets and thus the estimated cosupport set need to be large enough to avoid empty intersections.
Thus, for $\OM$ in general position we use the true cosparsity level for ACoSaMP and ASP. For $\OM_{DIF}$, where linear dependencies
occur, the corank does not equal the cosparsity and we use $\tilde{\cosp}$ instead of $\cosp$ since it will be favorable to run the algorithm targeting a cosparsity level in the middle. In this case $\cosp$ tends to be very large and it is more likely to have non-empty intersections .

\subsection{Phase Diagrams for Synthetic Signals in the Noiseless Case}

We begin with with synthetic signals in the noiseless case. We test the performance of AIHT with a constant step-size, AIHT with an adaptive changing step-size,
AHTP with a constant step-size, AHTP with an adaptive changing step-size,
ACoSaMP with $a=\frac{2\cosp-p}{\cosp}$, ACoSaMP with $a=1$,
ASP with $a=\frac{2\cosp-p}{\cosp}$ and ASP with $a=1$.
We compare the results to those of A-$\ell_1$-minimization \cite{elad07Analysis} and GAP \cite{Nam12Cosparse}.
We use a random matrix $\matr{M}$ and a random tight frame with $d=120$ and $p=144$,
where each entry in the matrices is drawn independently from the Gaussian distribution.

We draw a phase transition diagram \cite{Donoho09countingfaces} for each of the algorithms.
We test $20$ different possible values of $m$ and $20$ different values of $\cosp$ and for each pair repeat the experiment $50$ times.
In each experiment we check whether we have a perfect reconstruction. White cells in the diagram
denotes a perfect reconstruction in all the experiments of the pair and black cells denotes total failure in the reconstruction.
The values of $m$ and $\cosp$ are selected according to the following formula:
\begin{eqnarray}
m = \delta d && \cosp = d - \rho m,
\end{eqnarray}
where $\delta$, the sampling rate, is the x-axis of the phase diagram and
$\rho$, the ratio between the cosparsity of the signal and the number of measurements, is the y-axis.

\begin{figure*}[!t]
{\subfigure[AIHT, constant step-size]{\includegraphics[width=1.5in]{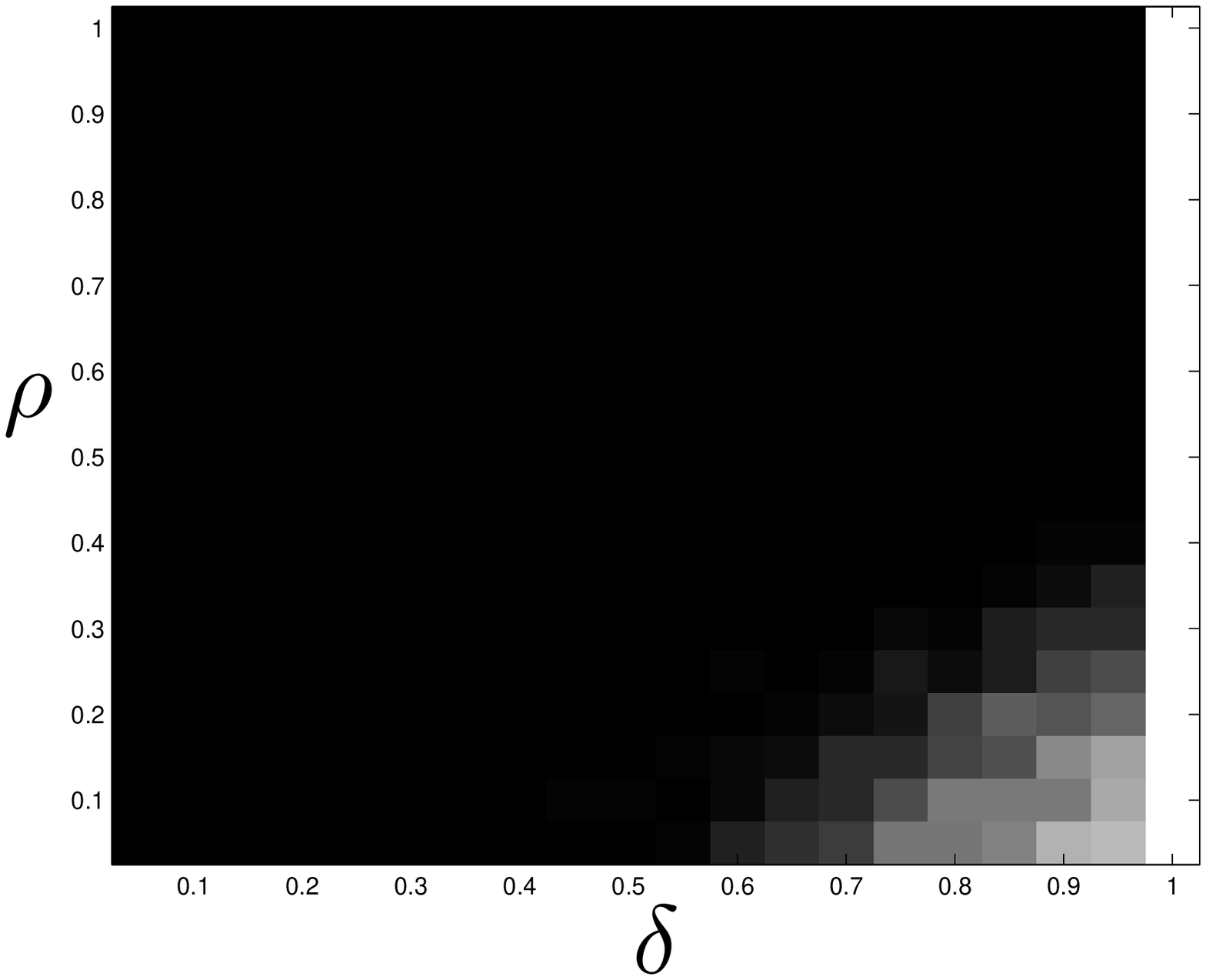}}%
\hfil
\subfigure[AIHT, adaptive step-size]{\includegraphics[width=1.5in]{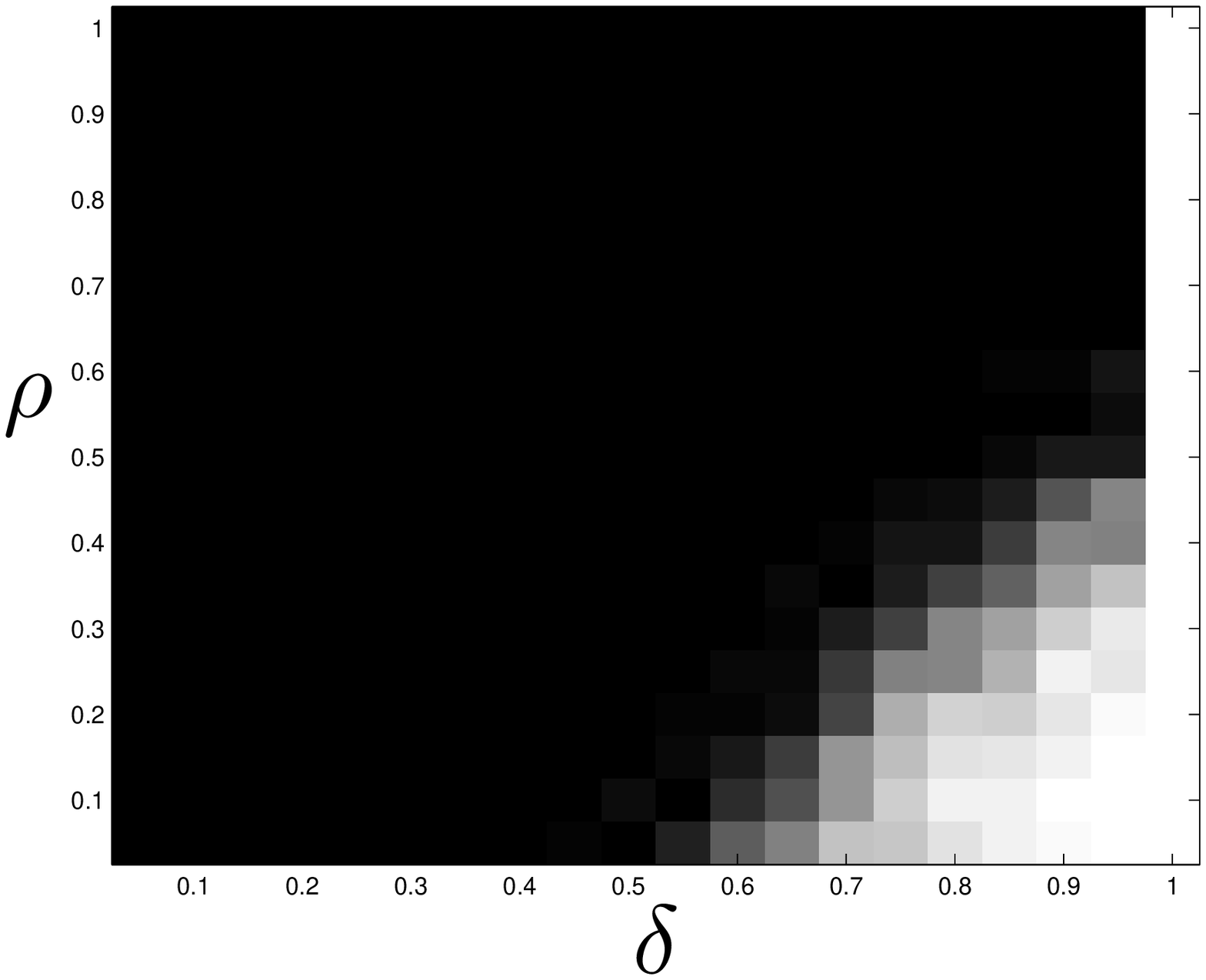}}}%
\hfil
\subfigure[AHTP, constant step-size]{\includegraphics[width=1.5in]{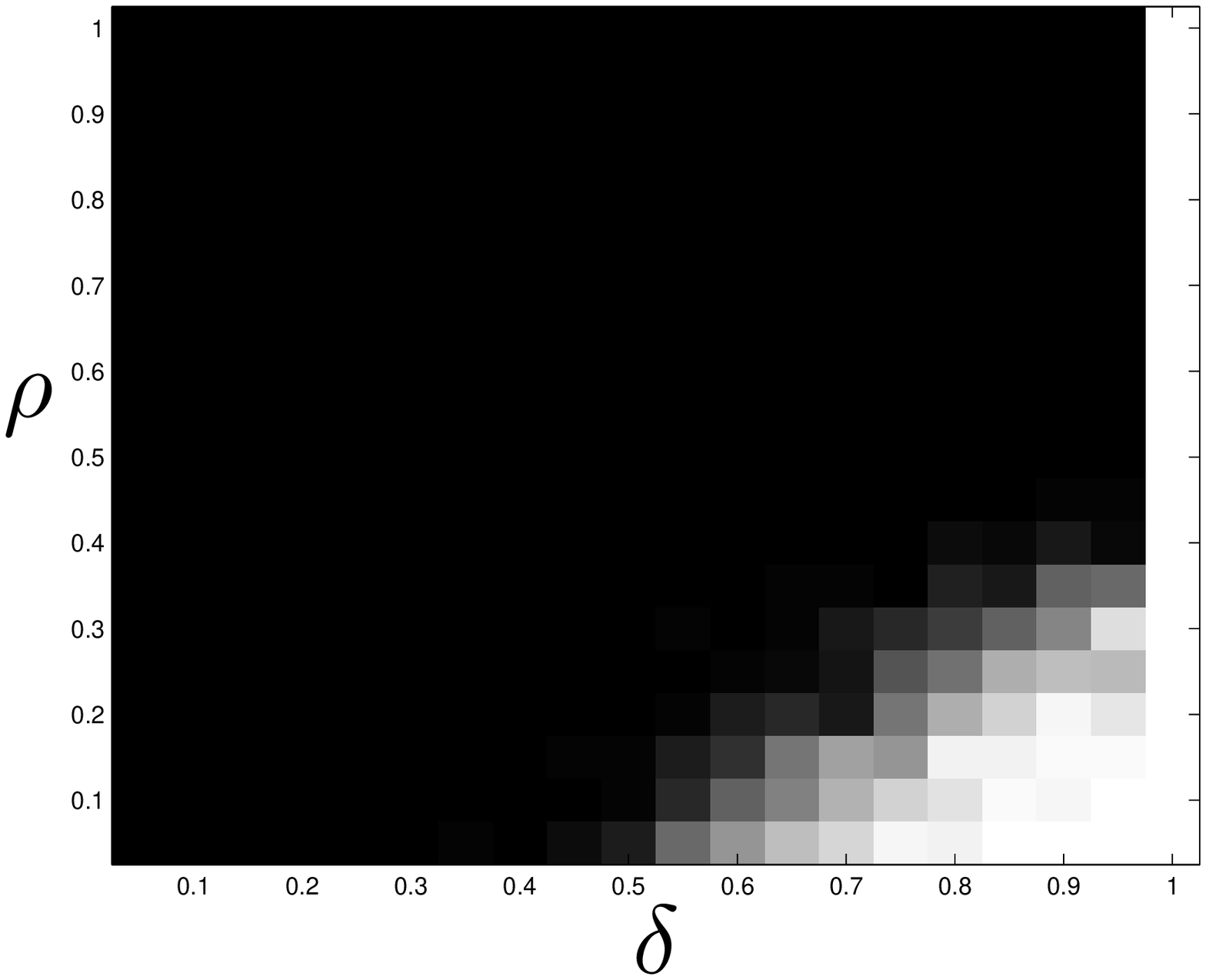}}%
\hfil
\subfigure[AHTP, adaptive step-size]{\includegraphics[width=1.5in]{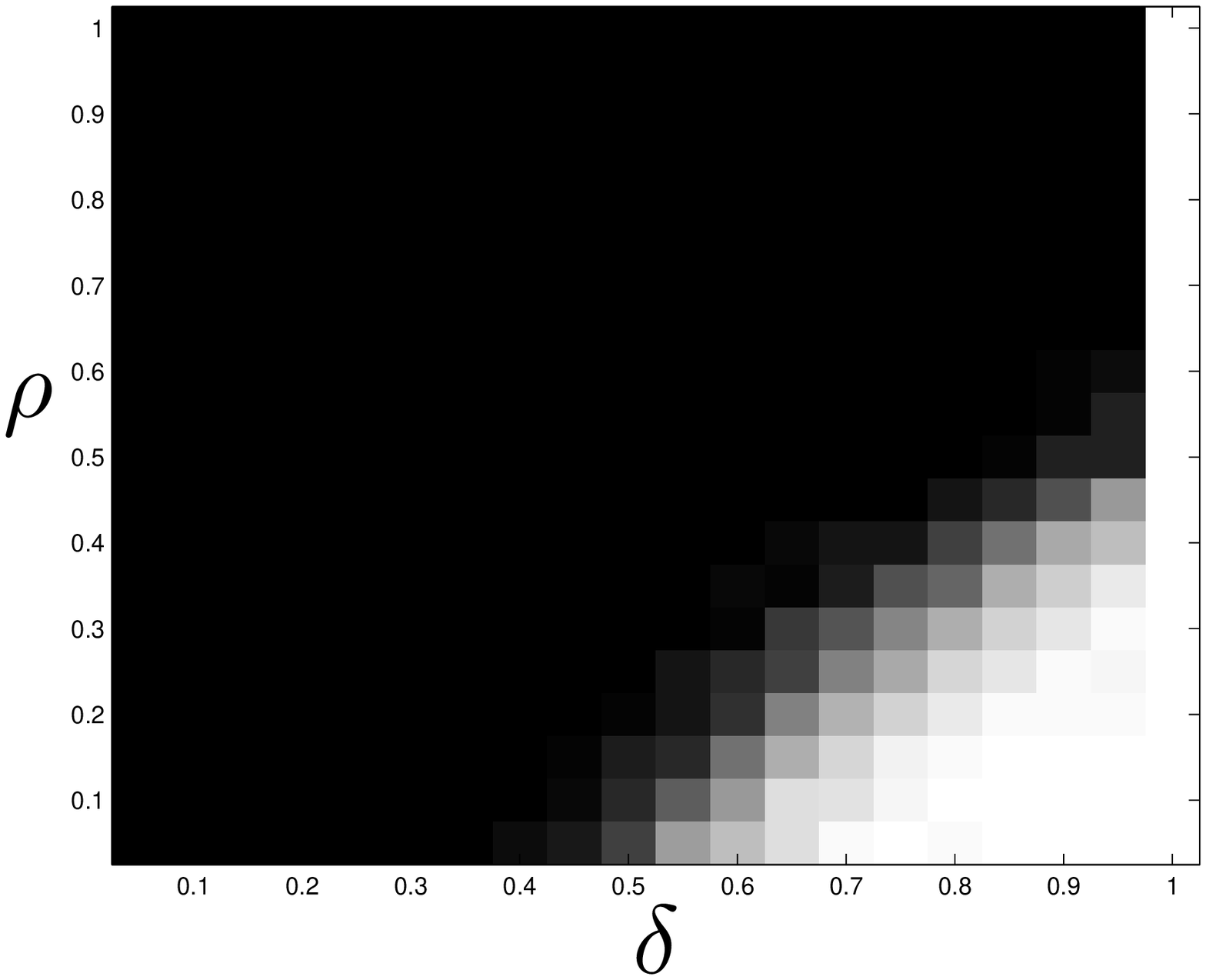}}%
\vfil
\subfigure[ACoSaMP, $a=\frac{2\cosp-p}{\cosp}$]{\includegraphics[width=1.5in]{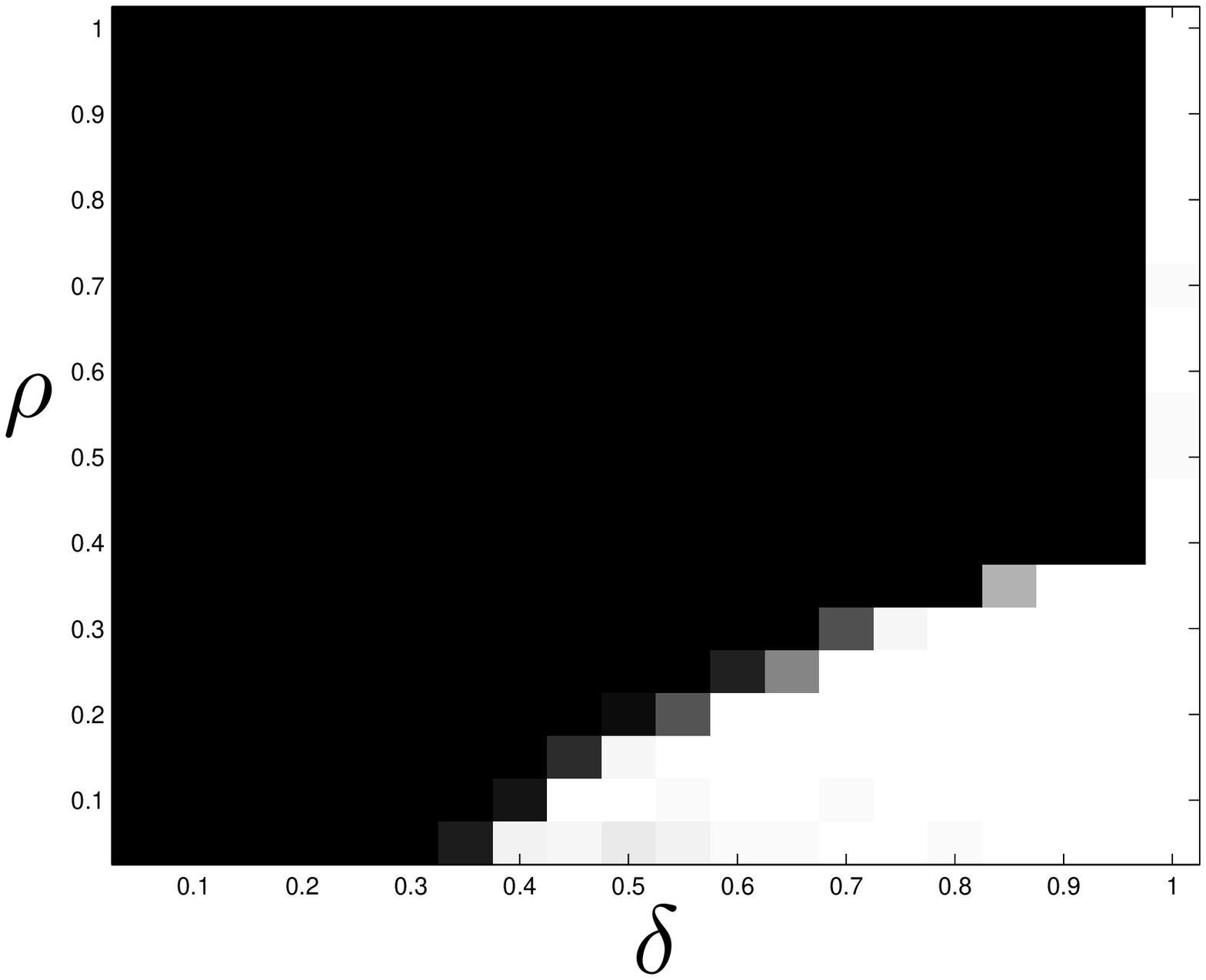}}%
\hfil
\subfigure[ACoSaMP, $a=1$]{\includegraphics[width=1.5in]{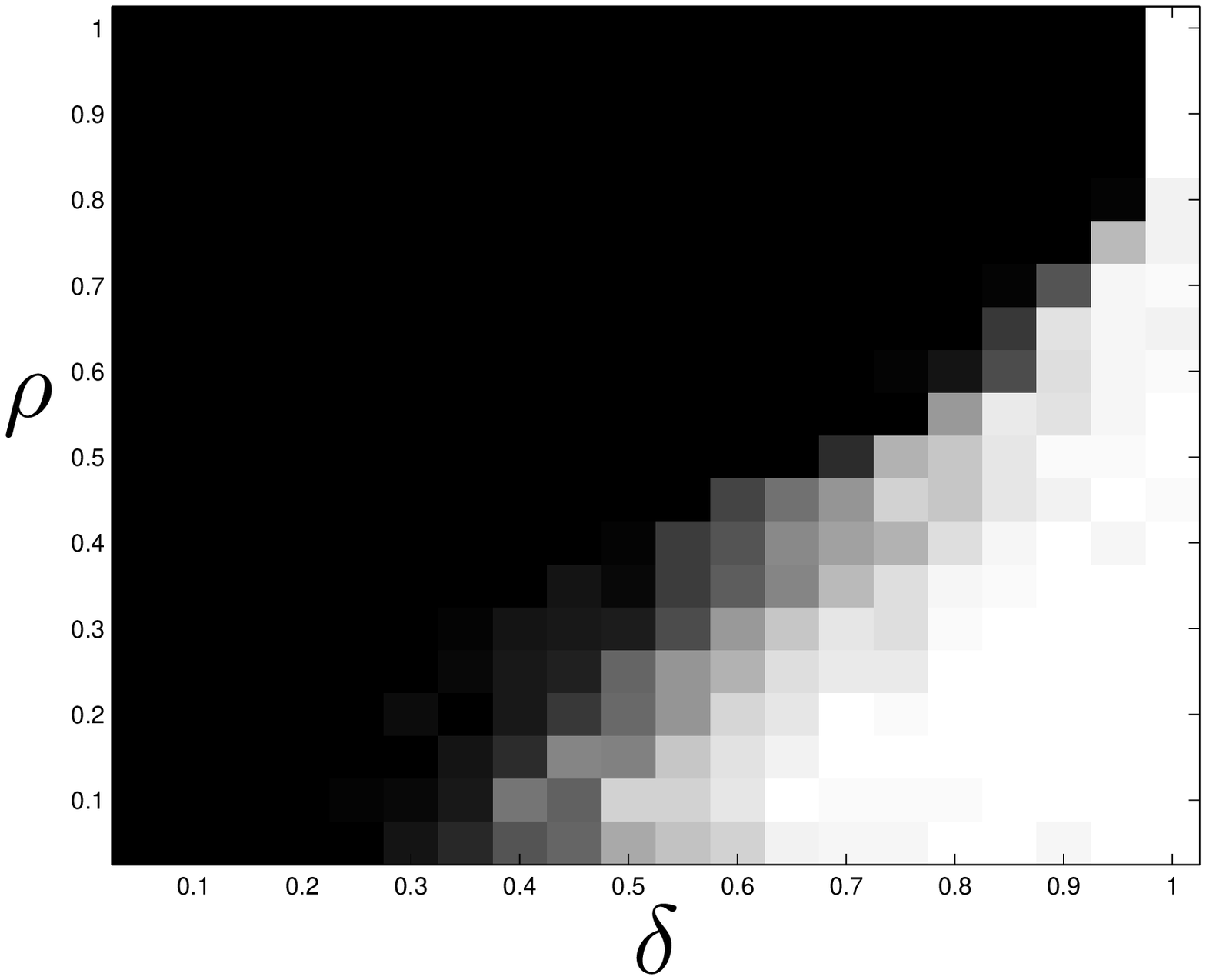}}%
\subfigure[ASP, $a=\frac{2\cosp-p}{\cosp}$]{\includegraphics[width=1.5in]{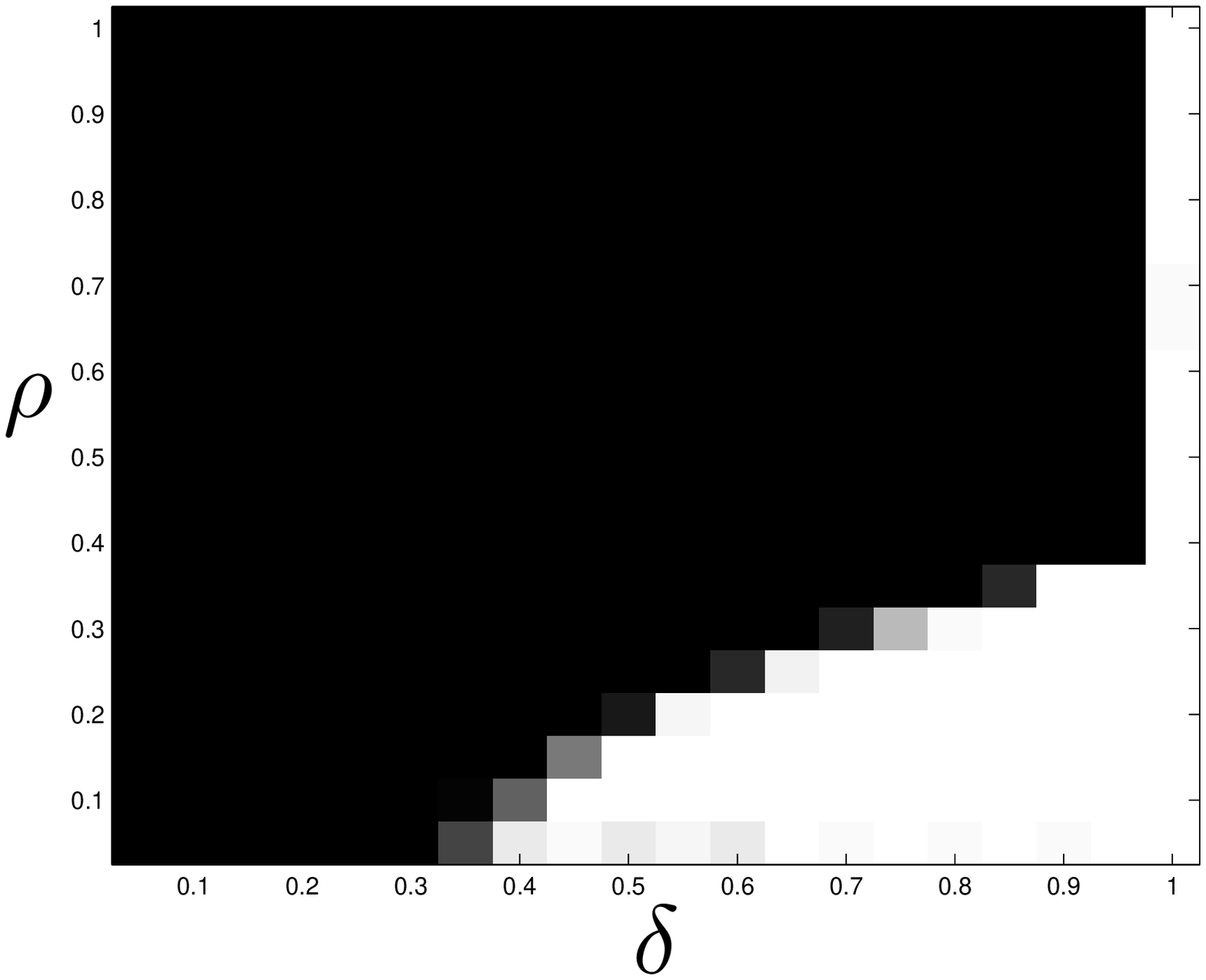}}%
\hfil
\subfigure[ASP, $a=1$]{\includegraphics[width=1.5in]{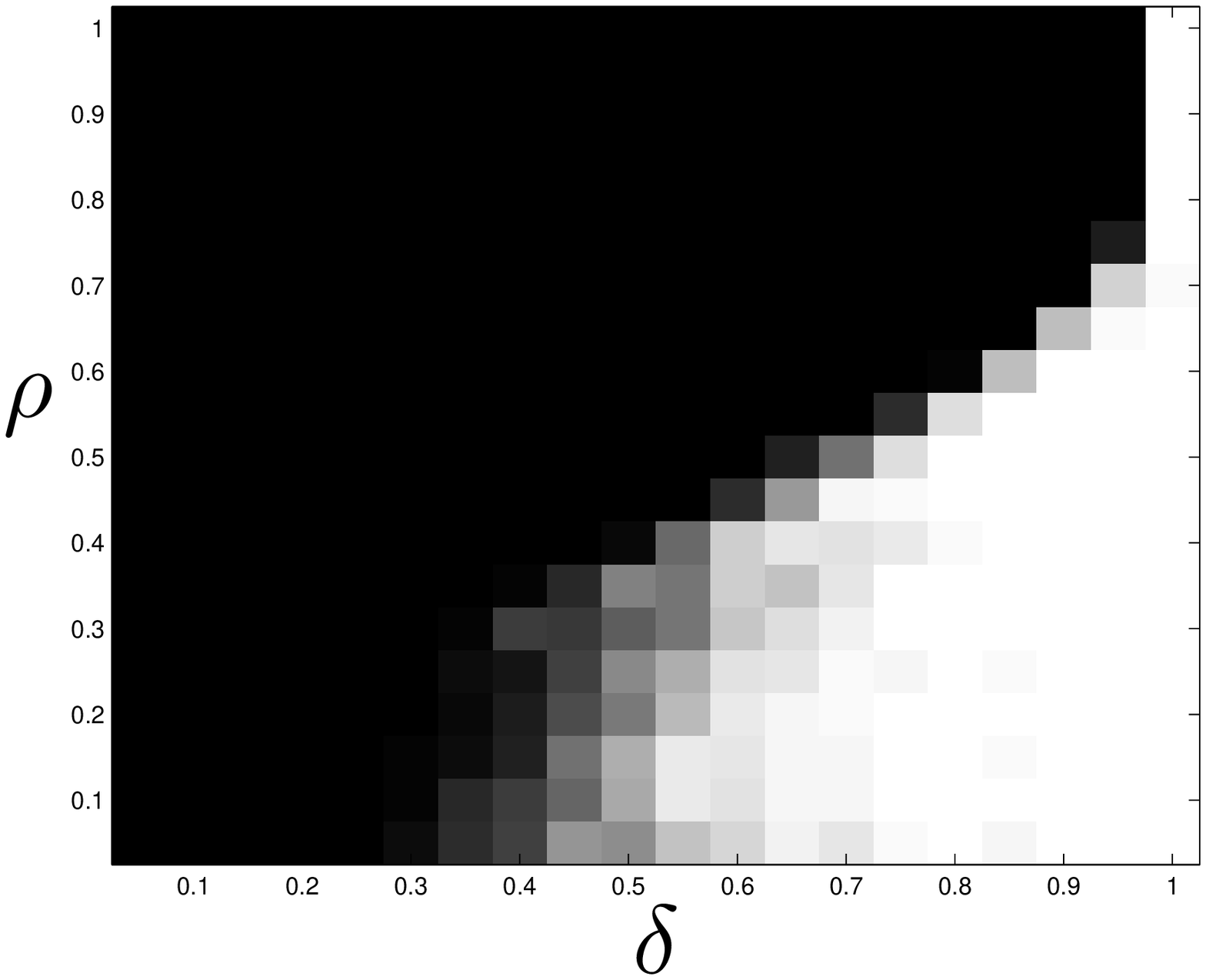}}%
\vfil
\subfigure[A-$\ell_1$-minimization]{\includegraphics[width=1.5in]{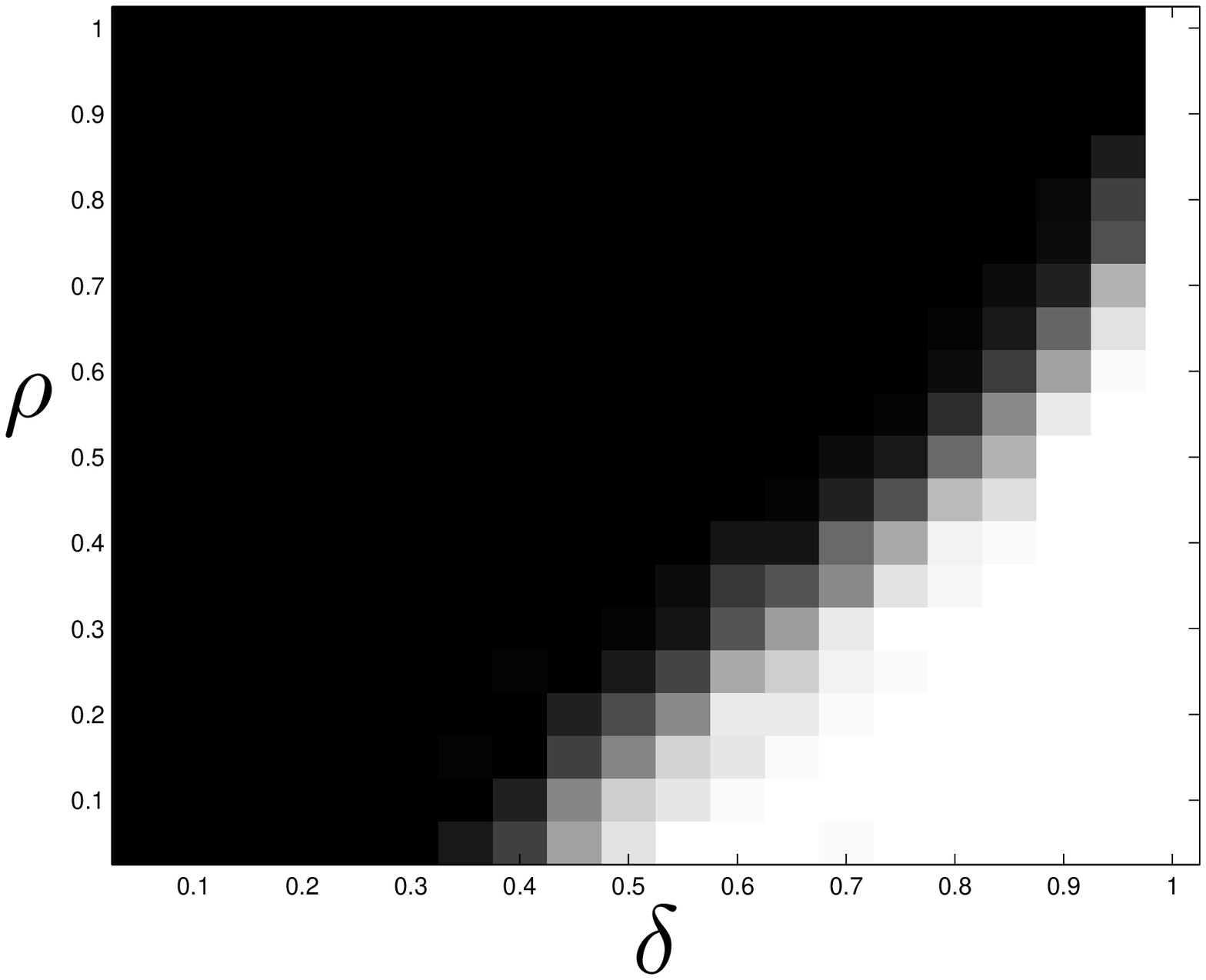}}%
\hfil
\subfigure[GAP]{\includegraphics[width=1.5in]{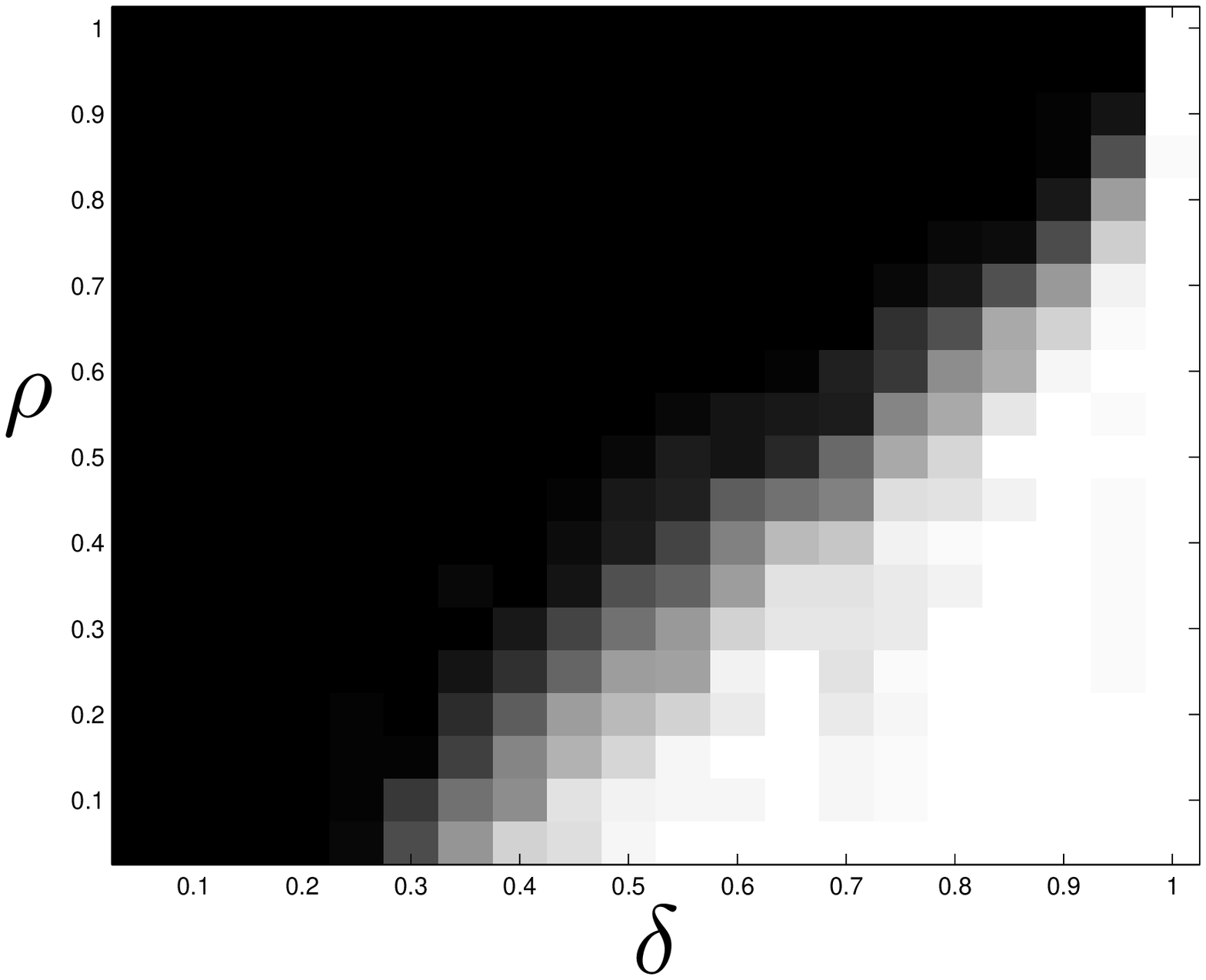}}%
\caption{Recovery rate for a random tight frame with $p=144$ and $d=120$. From left to right, up to bottom:
AIHT with a constant step-size, AIHT with an adaptive changing step-size,
AHTP with a constant step-size, AHTP with an adaptive changing step-size,
ACoSaMP with $a=\frac{2\cosp-p}{\cosp}$, ACoSaMP with $a=1$,
ASP with $a=\frac{2\cosp-p}{\cosp}$, ASP with $a=1$,
A-$\ell_1$-minimization and GAP.}
\label{fig:phaseDiagramAll1_2}
\end{figure*}

Figure~\ref{fig:phaseDiagramAll1_2} presents the reconstruction results of the algorithms.
It should be observed that AIHT and AHTP have better performance using the adaptive step-size than using the constant step-size.
The optimal step-size has similar reconstruction result like the adaptive one and thus not presented.
For ACoSaMP and ASP we observe that it is better to use $a=1$ instead of $a= \frac{2\cosp-p}{\cosp}$.
Compared to each other we see that ACoSaMP and ASP achieve better recovery than AHTP and AIHT. Between the last two, AHTP is better.
Though AIHT has inferior behavior, we should mention that with regards to running time AIHT is the most efficient. Afterwards we have AHTP
and then ACoSaMP and ASP.
Compared to $\ell_1$ and GAP we observe that ACoSaMP and ASP have competitive results.


\begin{figure}[t]
\begin{tabular}{ccccc}
\vspace{-0.15in} & {{\parbox{40mm}{\small AIHT, adaptive step-size}}} & {{\parbox{40mm}{\small AHTP, adaptive step-size}}} &
{{{\parbox{35mm}{\small ~~ACoSaMP, $a=1$}}}} & {{\parbox{35mm}{\small ~~~~ ASP, $a=1$}}} \\
\hspace{-0.15in} \begin{sideways}{{\parbox{35mm}{Random tight frame}}}\end{sideways} &\hspace{-0.3in}
\includegraphics[width=1.3in]{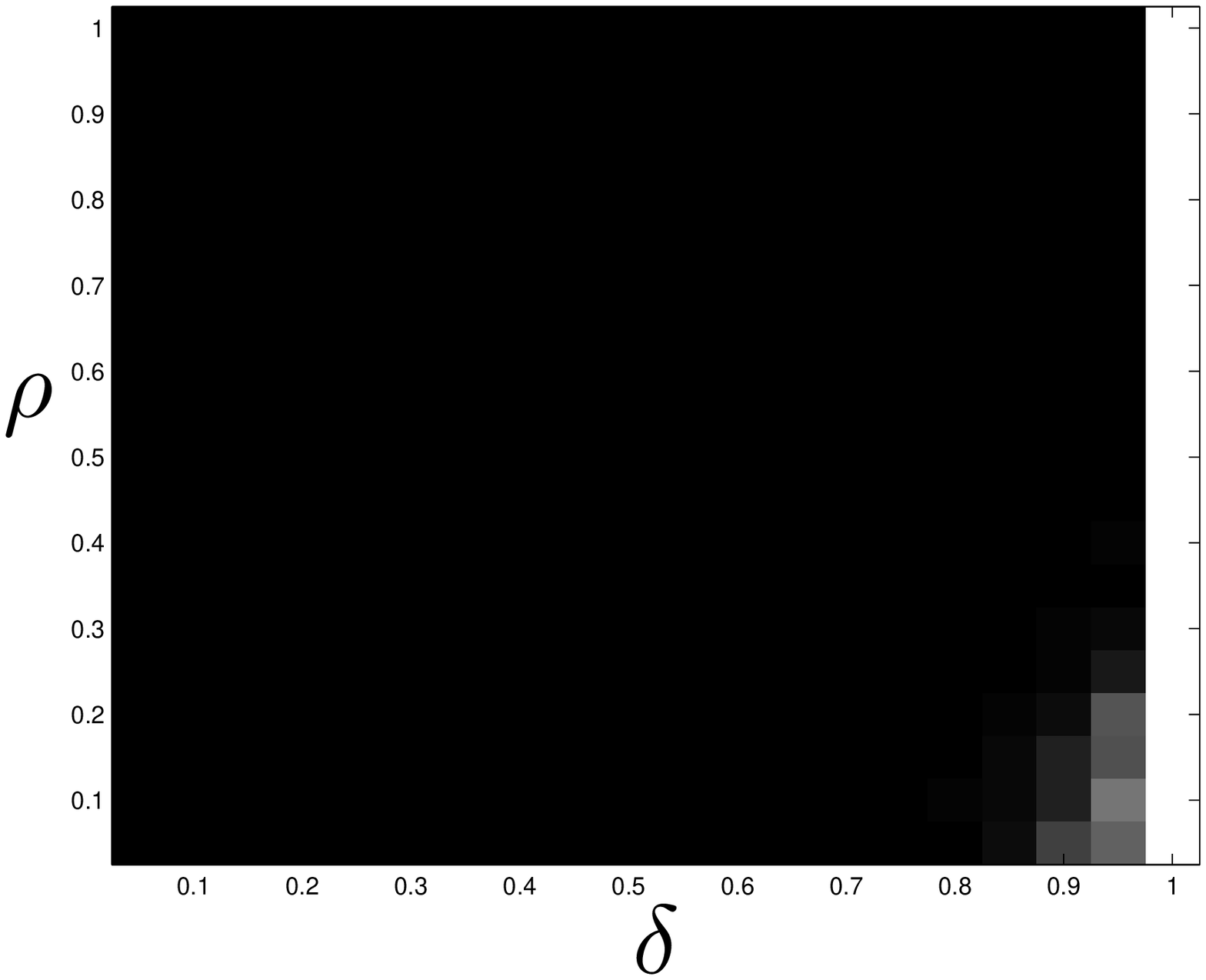} & \hspace{-0.35in} \includegraphics[width=1.3in]{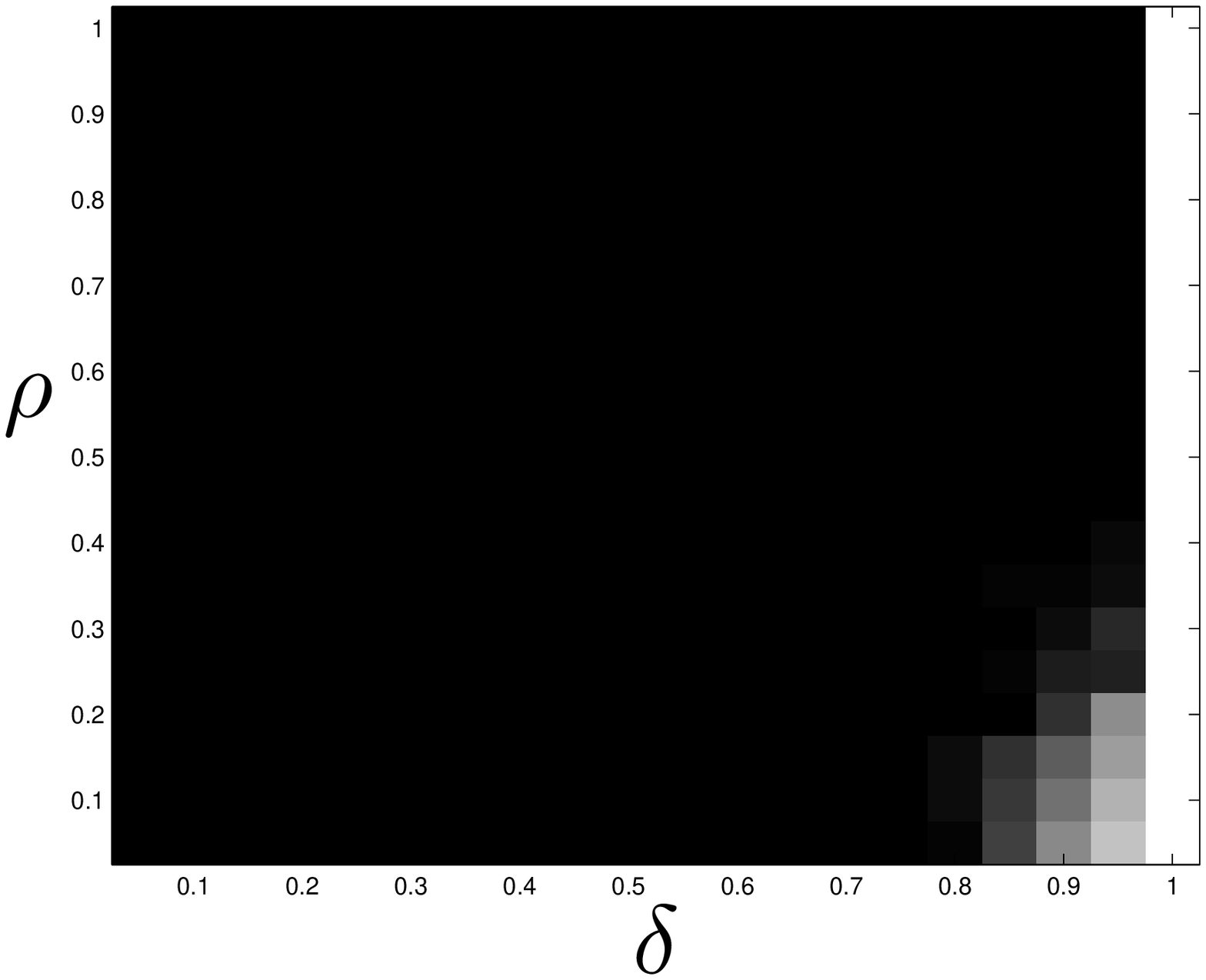} & \hspace{-0.35in}
\includegraphics[width=1.3in]{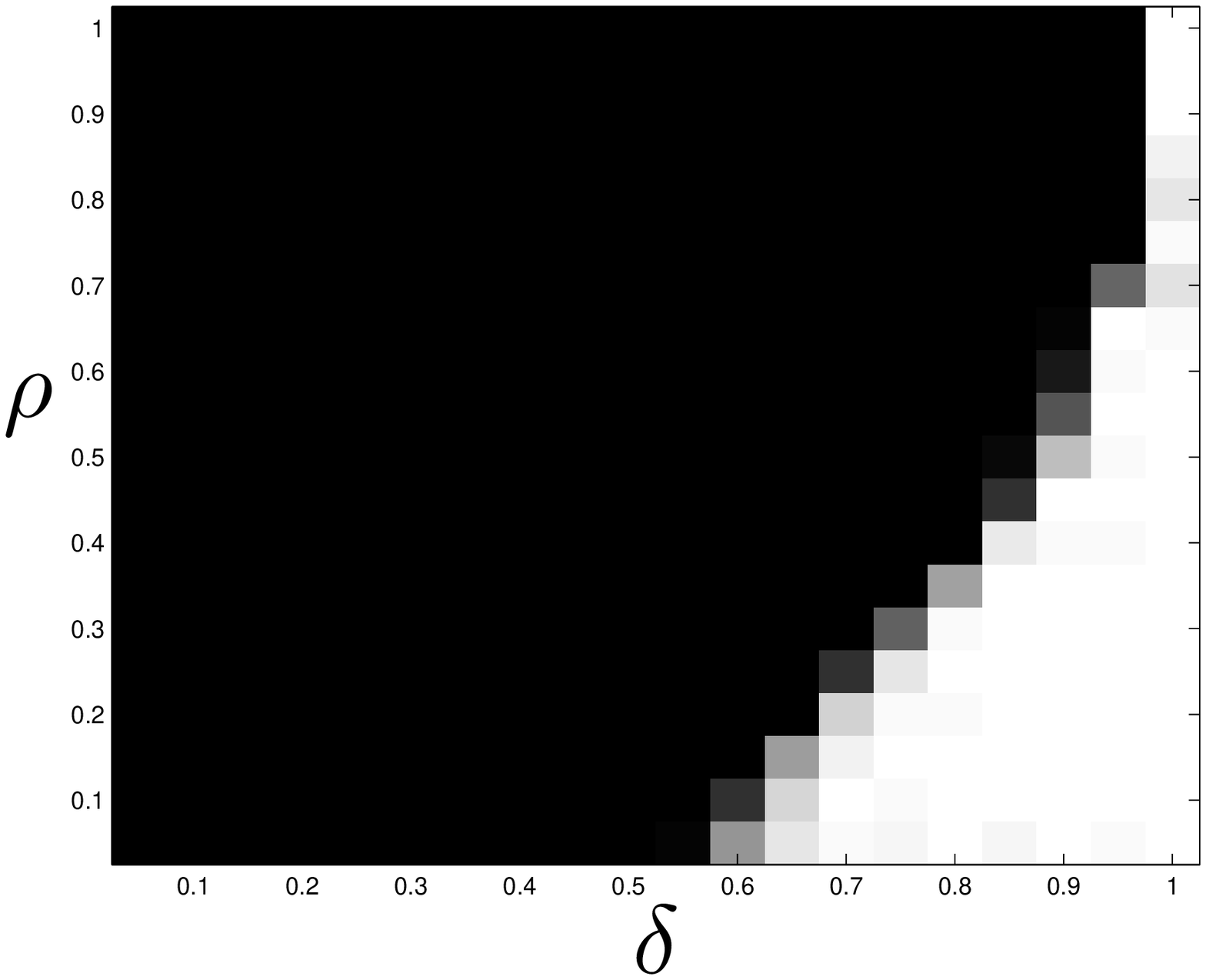} & \hspace{-0.35in} \includegraphics[width=1.3in]{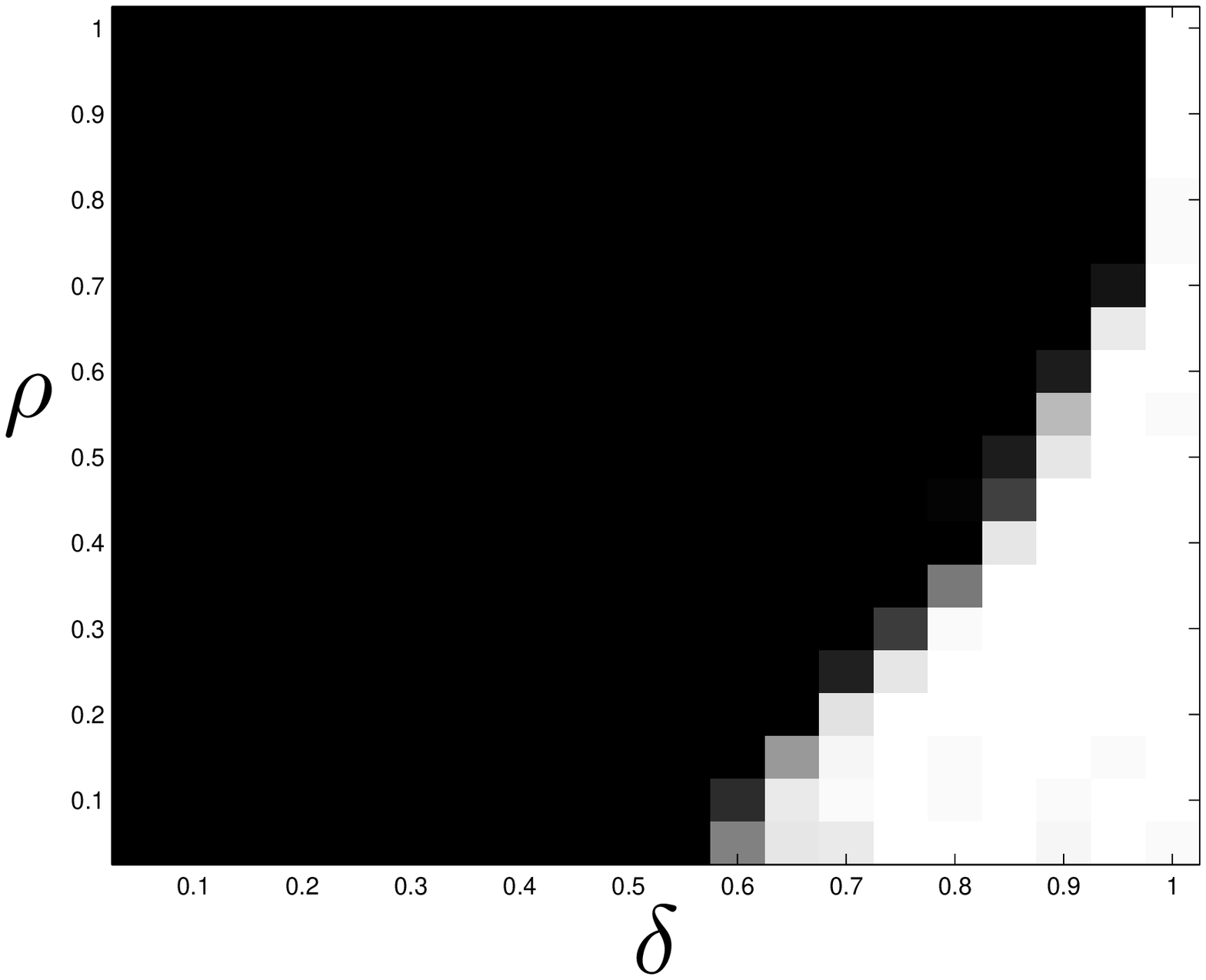}
\\ \hspace{-0.15in} \begin{sideways}{{\parbox{20mm}{~~~~~~~2D-DIF}}}\end{sideways} & \hspace{-0.3in}
\includegraphics[width=1.3in]{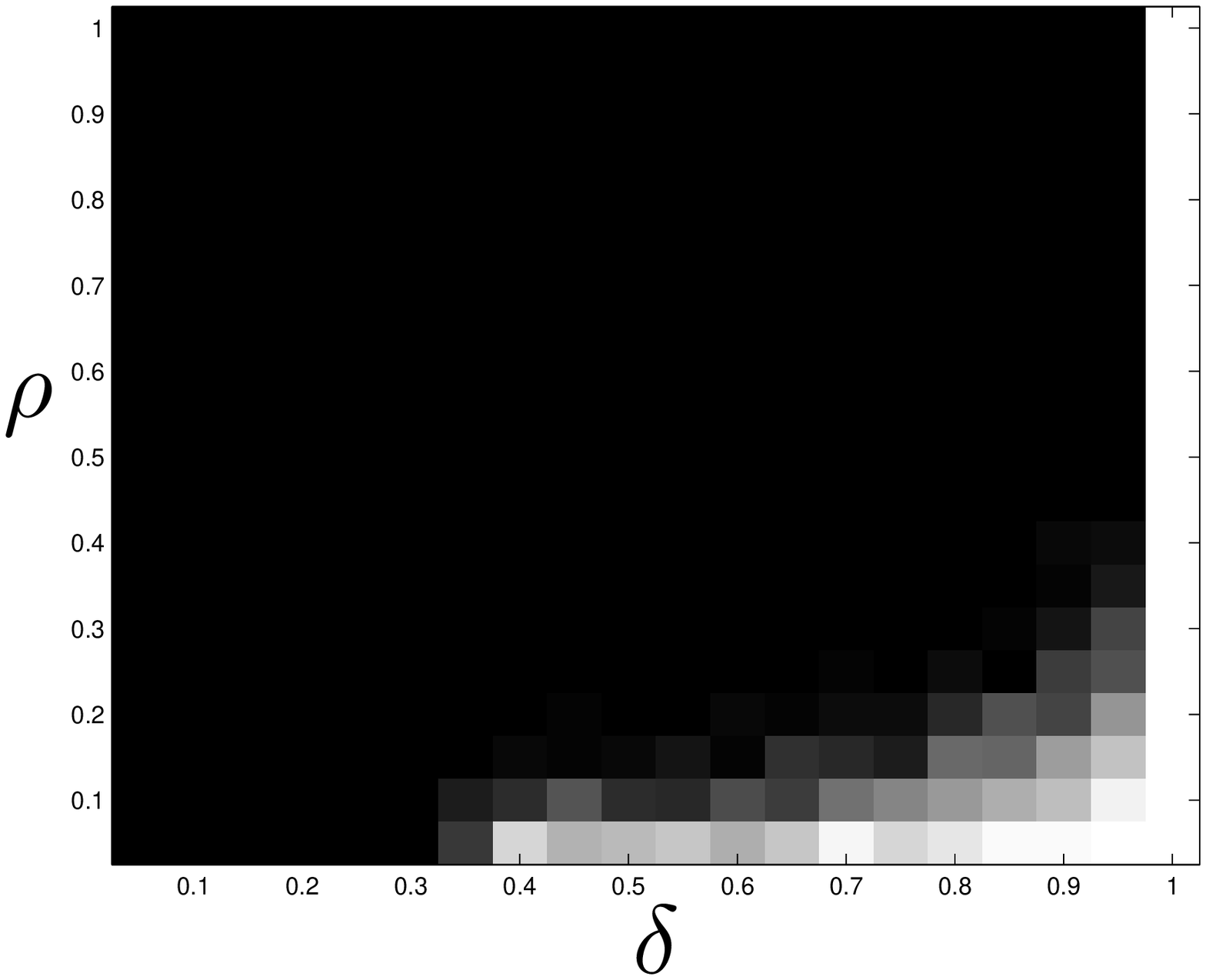} & \hspace{-0.35in} \includegraphics[width=1.3in]{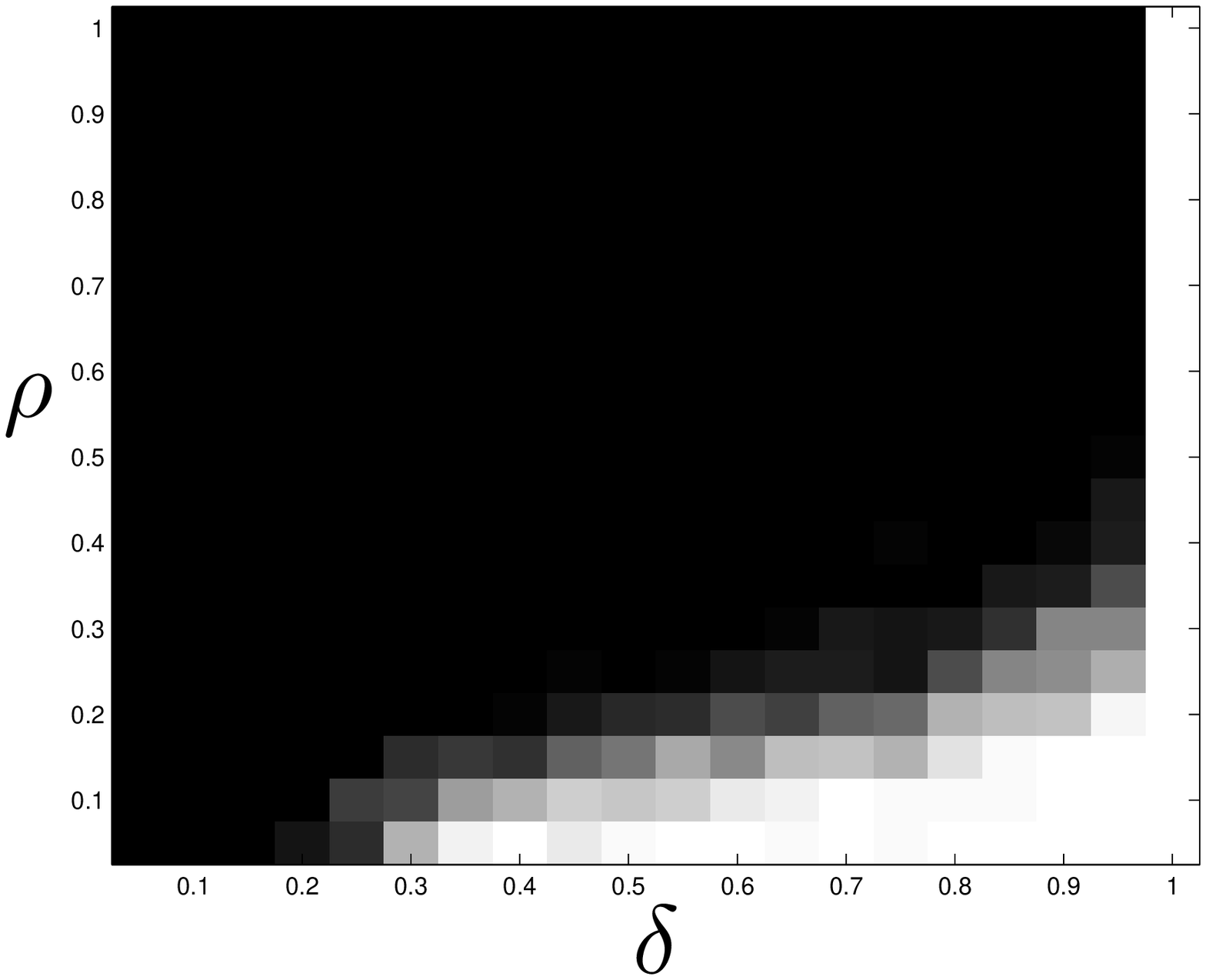} & \hspace{-0.35in}
\includegraphics[width=1.3in]{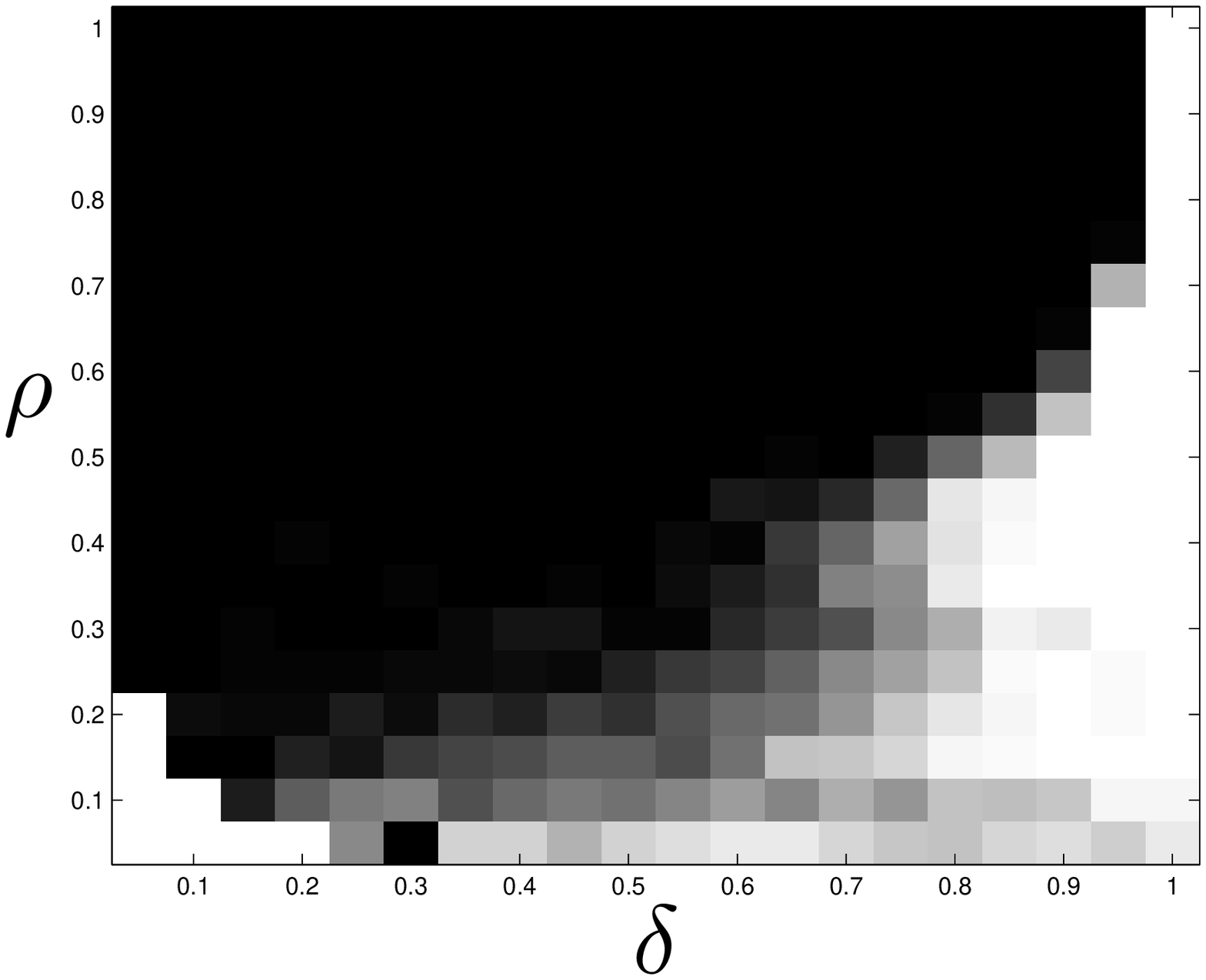} & \hspace{-0.35in} \includegraphics[width=1.3in]{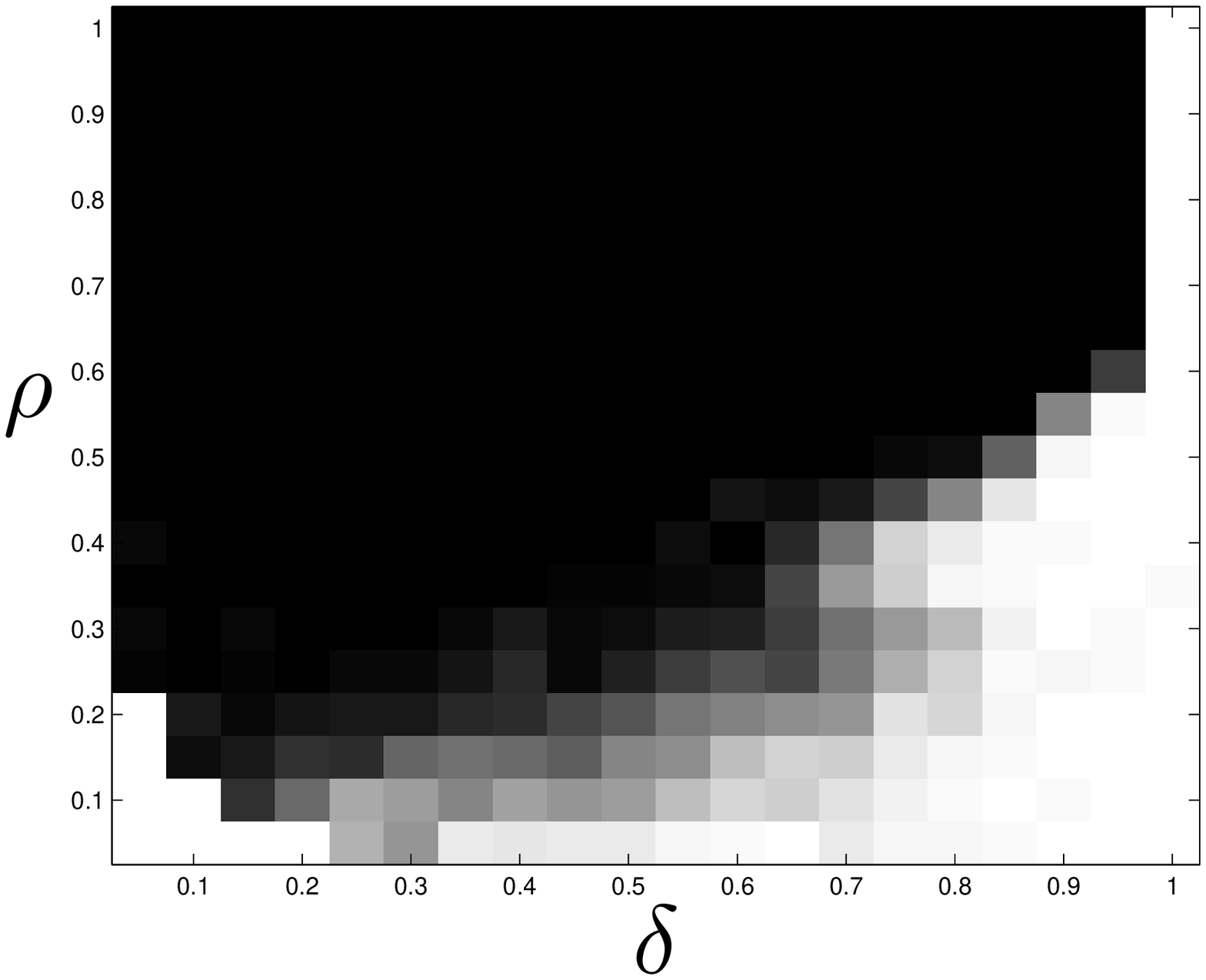}
\end{tabular}
\caption{Recovery rate for a random tight frame with $p=240$ and $d=120$ (up) and a finite difference operator (bottom).
From left to right:
AIHT and
AHTP with an adaptive changing step-size,
and ACoSaMP and ASP with $a=1$.}
\label{fig:phaseDiagramAll2}
\end{figure}

With the above observations, we turn to test operators with higher redundancy and see the effect of linear dependencies in them.
We test two operators. The first is a random tight frame as before but with redundancy factor of $2$. The second is the two dimensional finite difference operator $\OM_{\text{2D-DIF}}$.
In Fig.~\ref{fig:phaseDiagramAll2} we present the phase diagrams for both operators using
AIHT with an adaptive changing step-size,
AHTP with an adaptive changing step-size,
ACoSaMP with $a=1$,
and ASP with $a=1$.
As observed before, also in this case the ACoSaMP and ASP outperform AIHT and AHTP in both cases and AHTP outperform AIHT.
We mention again that the better performance comes at the cost of higher complexity.
In addition, as we expected, having redundancies in $\OM$ results with a better recovery.

\subsection{Reconstruction of High Dimensional Images in the Noisy Case}

\begin{figure}[!t]
\centering
{\subfigure[Phantom]{\includegraphics[width=1.5in]{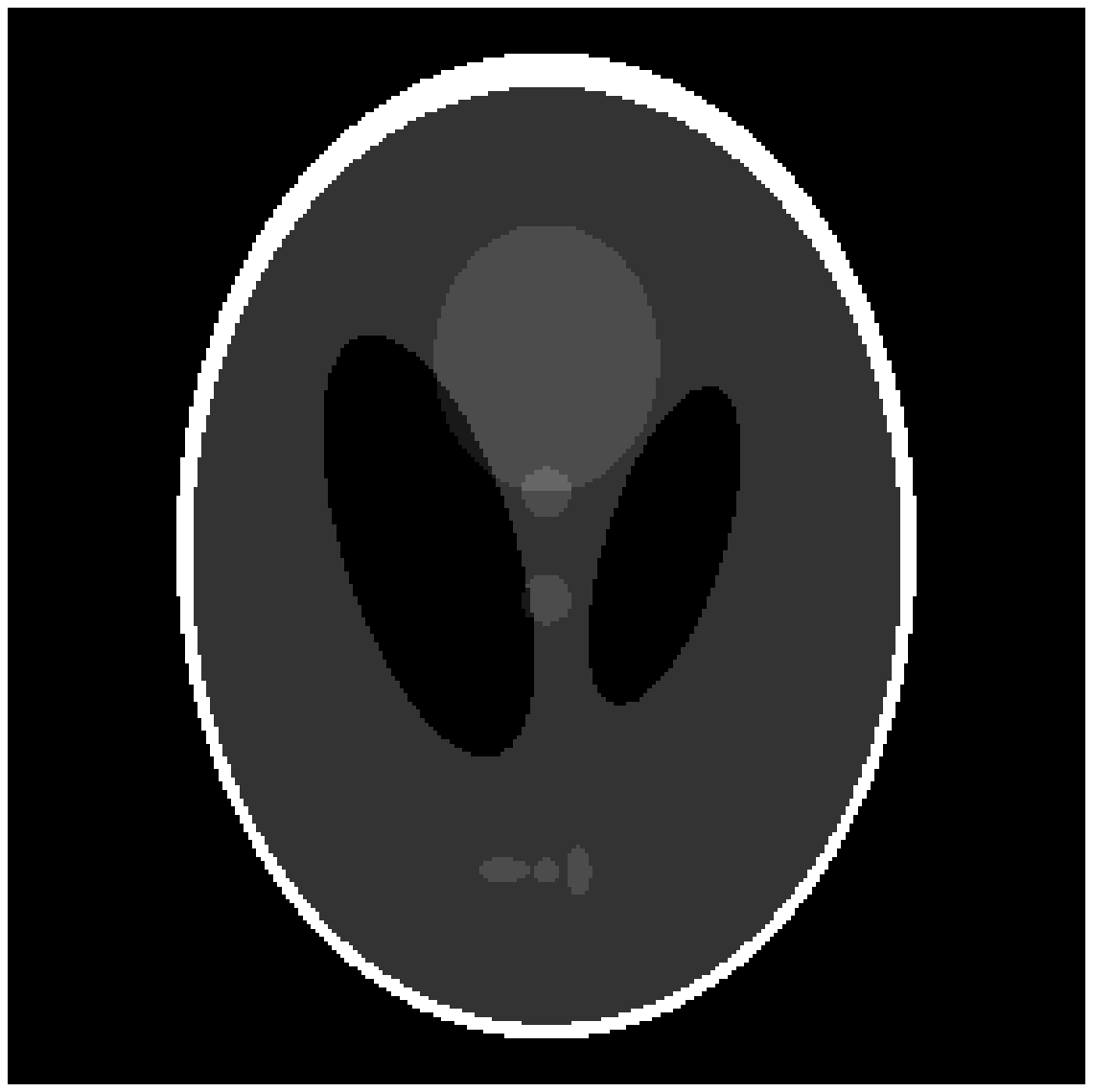} \label{fig:shepploganphantom}}%
\hfil
\subfigure[AIHT - noiseless]{\includegraphics[width=1.5in]{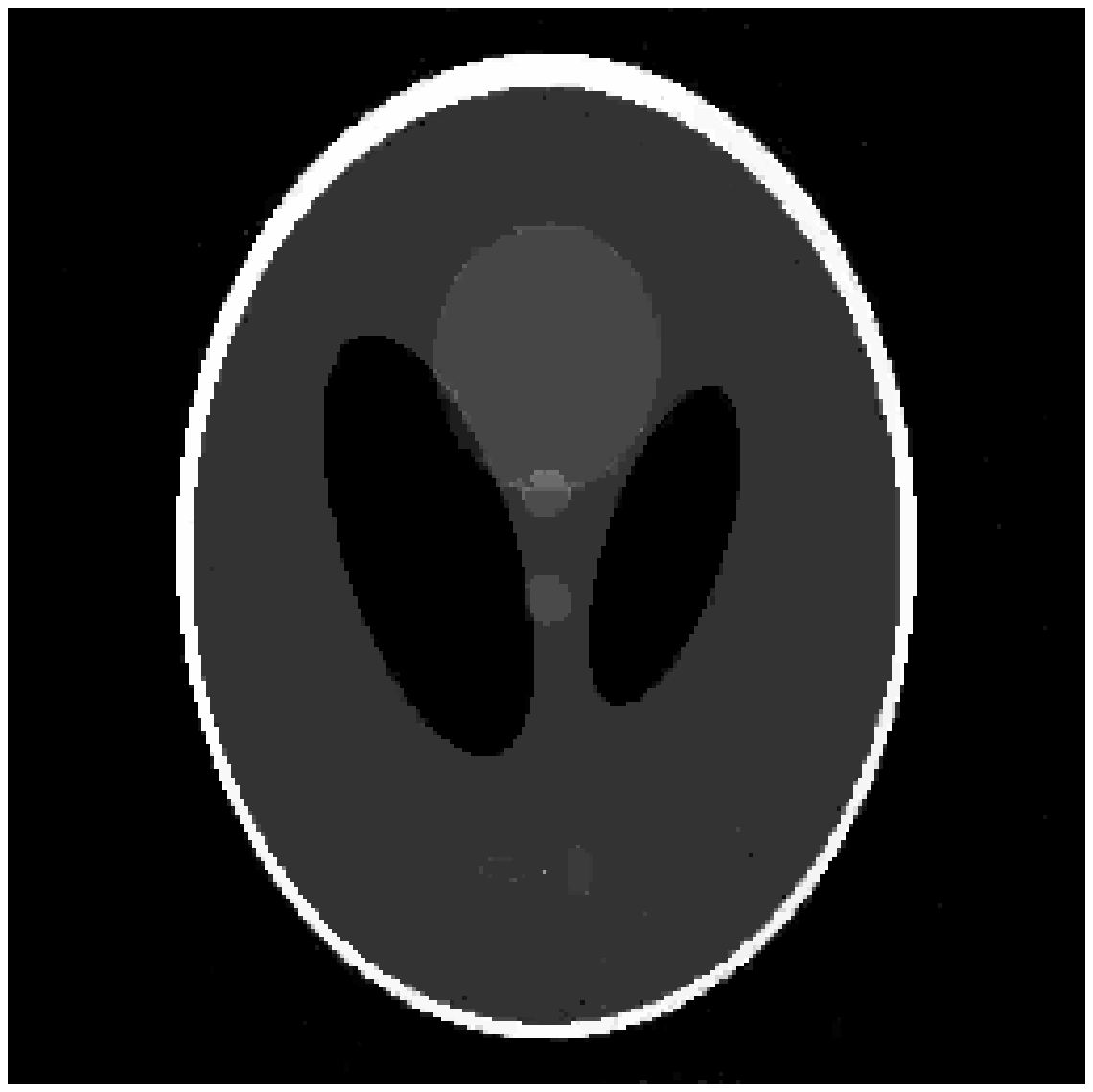} \label{fig:shepploganphantom_AIHT}}%
\hfil
\subfigure[Noisy Phantom]{\includegraphics[width=1.5in]{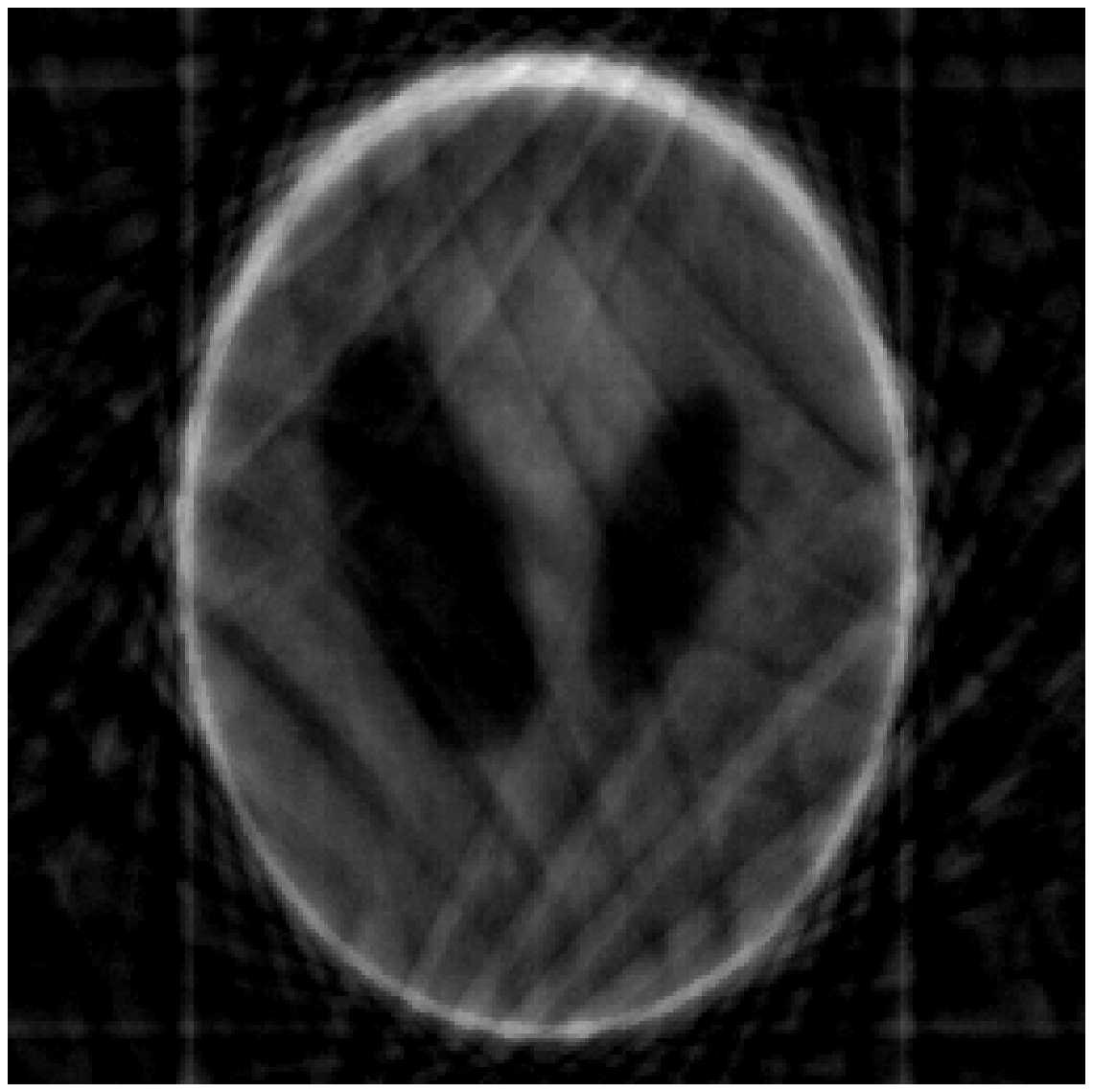} \label{fig:shepploganphantom_noisy}}%
\hfil
\subfigure[RASP - noisy]{\includegraphics[width=1.5in]{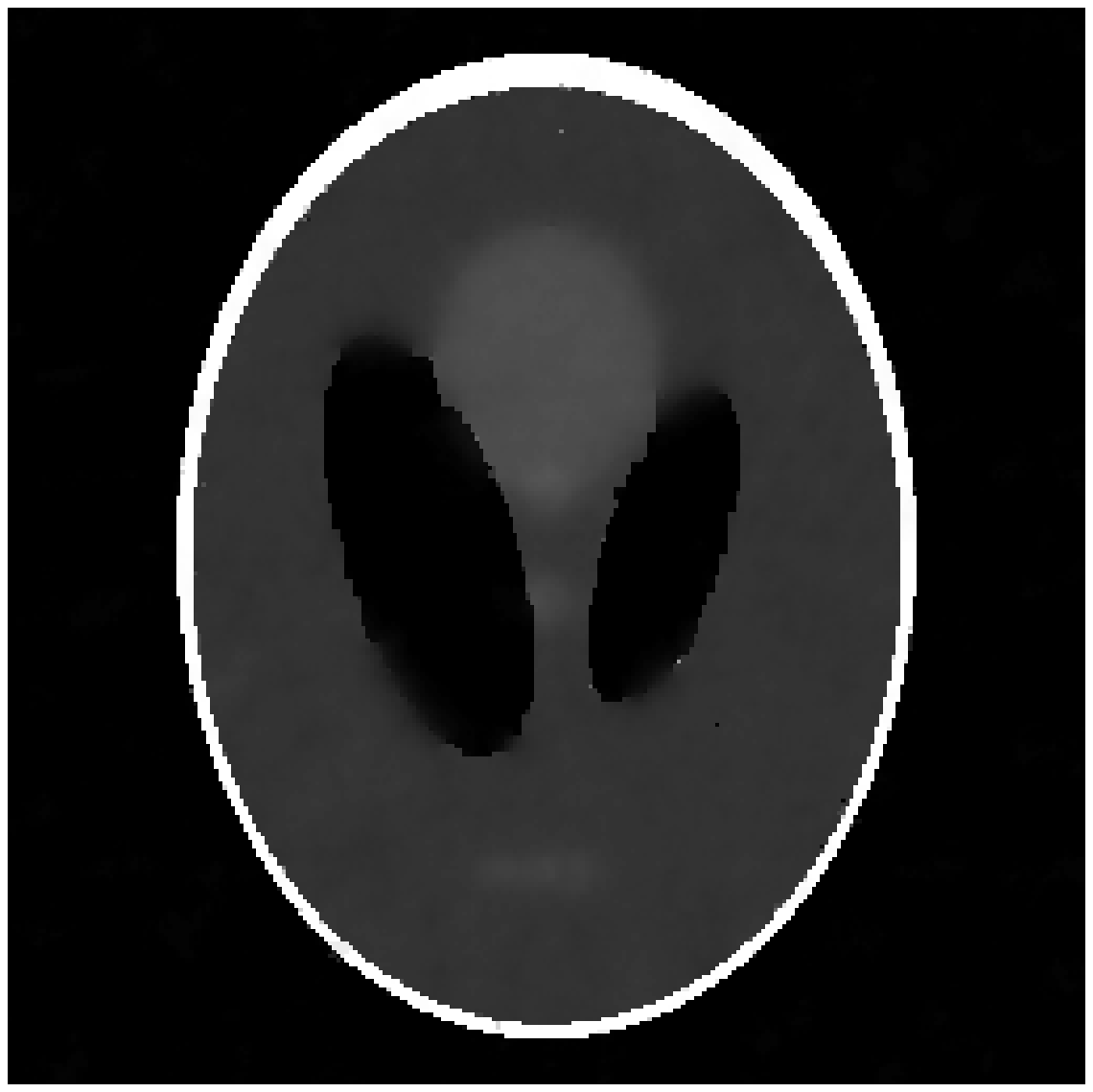} \label{fig:shepploganphantom_noisy_RASP}}}%
\caption{From left to right: Shepp Logan phantom image, AIHT reconstruction using $35$ radial lines,
 noisy image with SNR of $20$ and recovered image using RASP and only $22$ radial lines.
Note that for the noiseless case RASP and RACoSaMP get a perfect reconstruction using only $15$ radial lines.}
\label{fig:shepploganphantom_reconstruction}
\end{figure}

We turn now to test the methods for high dimensional signals. We use RASP and RACoSaMP (relaxed versions of ASP and ACoSaMP defined in Section~\ref{sec:relaxed_alg}) for
the reconstruction of the \emph{Shepp-Logan phantom} from few number of measurements. The sampling operator is a two dimensional Fourier transform that measures
only a certain number of radial lines from the Fourier transform.
The cosparse operator is $\OM_{\text{2D-DIF}}$ and the cosparsity used is the actual cosparsity of the signal under this operator ($\cosp = 128014$).
The phantom image is presented in Fig.~\ref{fig:shepploganphantom}. Using the RACoSaMP and RASP we get a perfect reconstruction using only $15$ radial lines, i.e., only $m=3782$ measurements out of $d=65536$ which is less then $6$ percent of the data in the original image. The algorithms requires less than $20$ iterations for having this perfect recovery.
For AIHT and RAHTP we achieve a reconstruction which is only close to the original image using $35$ radial lines. The reconstruction result of AIHT is presented in Fig~\ref{fig:shepploganphantom_AIHT}. The advantage of the AIHT, though it has an inferior performance, over the other methods is its running time. While the others need several minutes for each reconstruction, for the AIHT it takes only few seconds to achieve a visually reasonable result.

%

Exploring the noisy case, we perform a reconstruction using RASP of a noisy measurement of the phantom with $22$ radial lines and signal to noise ratio (SNR) of $20$.
Figure~\ref{fig:shepploganphantom_noisy} presents the noisy image, the result of applying inverse Fourier transform on the measurements,
and Fig.~\ref{fig:shepploganphantom_noisy_RASP} presents its reconstruction result.
Note that for the minimization process we solve conjugate gradients, in each iteration and take only the real part of the result and crop the values of the resulted image to be in the range of $[0,1]$. We get a peak SNR (PSNR) of $36dB$. We get similar results using RACoSaMP but using more radial lines (25).



\section{Discussion and Conclusion}
\label{sec:conc}

In this work we presented new pursuits for the cosparse analysis model.
A theoretical study of these algorithms was performed giving
guarantees for stable recovery under the assumptions of the $\OM$-RIP
and the existence of an optimal or a near optimal projection. We showed that optimal projections exists
for some non-trivial operators, i.e., operators that do not take us back to the synthesis case.
In addition, we showed experimentally that using
simpler kind of projections is possible in order to get good reconstruction results.
We demonstrated both in the theoretical and the empirical results that linear dependencies within
the analysis dictionary are favorable and enhance the recovery performance.

We are aware that there are still some open questions in this work and we leave them for future research.
This should deal with following:
\begin{itemize}
\item Our work assumed the existence of a procedure that finds a cosupport that implies a near optimal projection with a constant $C_\cosp$. Two examples for optimal cosupport slection schemes were given.
However, the existence of an optimal or a near optimal scheme for a general operator is still an open question.
The question is: for which types of $\OM$ and values of $C_\cosp$
we can find an efficient procedure that implies a near optimal projection.
\item As we have seen in the simulations, the thresholding procedure, though not near optimal with the theorems required constants,
provides good reconstruction results. A theoretical study of the analysis greedy-like techniques with this cosupport
selection scheme is required.
\item A family of analysis dictionaries that deserves a special attention is the family of tight frame operators.
In synthesis, there is a parallel between the guarantees
of $\ell_1$-synthesis and the greedy like algorithms.
The fact that a guarantee with a tight frame $\OM$ exists
for $\ell_1$-analysis encourage us to believe that similar guarantees exist also for the analysis greedy-like techniques.
\item In this paper, the noise $\e$ was considered to be adversarial. Random white Gaussian case
was considered for the synthesis case in \cite{Giryes12RIP} resulting with near-oracle performance guarantees.
It would be interesting to verify whether this is also the case for the analysis framework.
\end{itemize}


\appendix

\section{Proofs of Theorem~\ref{thm:analysis_RIP_cond} and Theorem~\ref{thm:analysis_RIP_cond_dependencies}}
\label{sec:analysis_RIP_proof}

{\em Theorem~\ref{thm:analysis_RIP_cond} (Theorem 3.3 in \cite{Blumensath09Sampling}):}
Let $\matr{M}\in \RR{m \times d}$ be a random matrix that satisfies that for any $\vect{z}\in \RR{d}$ and $0<\tilde{\epsilon} \le \frac{1}{3}$
\begin{eqnarray*}
P\left(\abs{\norm{\matr{M}\vect{z}}_2^2 - \norm{\vect{z}}_2^2} \ge \tilde{\epsilon}\norm{\vect{z}}_2^2\right) \le e^{-\frac{C_{\matr{M}}m\tilde{\epsilon}}{2}},
\end{eqnarray*}
where $C_{\matr{M}} >0$ is a constant.
For any value of $\epsilon_r >0$, if
\begin{eqnarray*}
m \ge \frac{32}{C_M\epsilon_r^2}\left( \log(\abs{\L_{r}^{\text{corank}}}) + (d-r)\log({9}/{\epsilon_r})+t\right),
\end{eqnarray*}
then $\delta_{r}^{\text{corank}} \le \epsilon_r$ with probability exceeding $1-e^{-t}$.

{\em Theorem~\ref{thm:analysis_RIP_cond_dependencies}:}
Under the same setup of Theorem~\ref{thm:analysis_RIP_cond}, for any $\epsilon_\cosp > 0$ if
\begin{eqnarray*}
m \ge \frac{32}{C_M\epsilon_\cosp^2}\left( (p-\cosp)\log\left(\frac{9p}{(p-\cosp)\epsilon_\cosp}\right)+t \right),
\end{eqnarray*}
then $\delta_\cosp \le \epsilon_\cosp$ with probability exceeding $1-e^{-t}$.

{\em Proof:}
Let $\tilde{\epsilon} = \epsilon_{r}/4$, $B^{\sdim-r} = \{\vect{z}\in \RR{d-r},\norm{\vect{z}}_2\le 1\}$ and
$\Psi$ an $\tilde{\epsilon}$-net for $B^{\sdim-r}$ with size $\abs{\Psi} \le \left(1+ \frac{2}{\tilde{\epsilon}} \right)^{\sdim-r}$ \cite{Mendelson08Uniform}.
For any subspace $\aspace_{\Lambda}^B = \aspace_{\Lambda} \cap B^{d-r}$ such that $\Lambda \in \L_{r}^{\text{corank}}$ we can build an orthogonal matrix $\matr{U}_{\matr{\Omega}} \in \RR{\sdim\times (\sdim - r)}$ such that $\aspace_{\Lambda}^B = \{\matr{U}_{\Lambda}\vect{z}, \vect{z}\in \RR{\sdim-r}, \norm{\vect{z}}_2\le 1 \} = \matr{U}_{\Lambda}B^{\sdim-r}$. It is easy to see that $\Psi_{\Lambda}=\matr{U}_\Lambda \Psi^{\sdim-r}$ is an $\tilde{\epsilon}$-net for $W_{\Lambda}^B$
and that $\Psi_{\A_{r}^{\text{corank}}} =\cup_{\Lambda \in \L_{r}^{\text{corank}}} \Psi_\Lambda$ is an $\tilde{\epsilon}$-net
for $\A_{r}^{\text{corank}}\cap B^{\sdim}$, where $\abs{\Psi_{\A_{r}^{\text{corank}}}} \le \abs{\L_{r}^{\text{corank}}}(1+\frac{2}{\tilde{\epsilon}})^{\sdim-r}$.

We could stop here and use directly Theorem~2.1 from \cite{Mendelson08Uniform} to get the desired result for Theorem~\ref{thm:analysis_RIP_cond}.
However, we present the remaining of the proof using a proof technique from \cite{Baraniuk08Simple,Rauhut08Compressed}.
Using union bound and the properties of $\matr{M}$ we have that with probability exceeding $1- \abs{\L_{r}^{\text{corank}}}(1+\frac{2}{\tilde{\epsilon}})^{\sdim-r}e^{-\frac{C_\M m\tilde{\epsilon}^2}{2}}$ every $\vect{v} \in \Psi_{\A_{r}^{\text{corank}}}$ satisfies
\begin{eqnarray}
\label{eq:epsilon_net_JL}
(1- \tilde{\epsilon})\norm{\vect{v}}_2^2 \le \norm{\matr{M}\vect{v}}_2^2 \le (1+\tilde{\epsilon})\norm{\vect{v}}_2^2.
\end{eqnarray}
According to the definition of $\delta_{r}^{\text{corank}}$ it holds that $\sqrt{1+\delta_{r}^{\text{corank}}} = \sup_{\vect{v}\in \A_{r}^{\text{corank}}\cap B^d}\norm{\matr{M}\vect{v}}_2$. Since $\A_{r}^{\text{corank}} \cap B^\sdim$ is a compact set there exists $\vect{v}_0 \in \A_{r}^{\text{corank}} \cap B^\sdim$ that achieves the supremum. Denoting by $\tilde{\vect{v}}$ its closest vector in $\Psi_{\A_{r}^{\text{corank}}}$ and using the definition of $\Psi_{\A_{r}^{\text{corank}}}$ we have $\norm{\vect{v}_0 - \tilde{\vect{v}}}_2 \le \tilde{\epsilon}$. This yields
\begin{eqnarray}
\label{eq:Omega_RIP_epsilon_stage}
\sqrt{1+\delta_{r}^{\text{corank}}} =  \norm{\matr{M}\vect{v}_0}_2 &\le& \norm{\matr{M}\tilde{\vect{v}}}_2 + \norm{\matr{M}(\vect{v}_0-\tilde{\vect{v}})}_2
\\ \nonumber &\le& \sqrt{1+\tilde{\epsilon}} + \norm{\matr{M}\frac{\vect{v}_0-\tilde{\vect{v}}}{\norm{\vect{v}_0-\tilde{\vect{v}}}}_2}_2\norm{\vect{v}_0-\tilde{\vect{v}}}_2
\le \sqrt{1+\tilde{\epsilon}} + \sqrt{1+\delta_{r}^{\text{corank}}}\tilde{\epsilon}.
\end{eqnarray}
The first inequality is due to the triangle inequality; the second one follows from \eqref{eq:epsilon_net_JL} and arithmetics;
and the last inequality follows from the definition of $\delta_{r}^{\text{corank}}$, the properties of $\tilde{\epsilon}$-net and the
fact that $\norm{\frac{\vect{v}_0-\tilde{\vect{v}}}{\norm{\vect{v}_0-\tilde{\vect{v}}}_2}}_2 =1$.
Reordering \eqref{eq:Omega_RIP_epsilon_stage} gives
\begin{eqnarray}
\label{eq:delta_l_epsilon_l_ineq}
1+\delta_{r}^{\text{corank}} \le \frac{1+\tilde{\epsilon}}{(1-\tilde{\epsilon})^2} \le 1 + 4\tilde{\epsilon} = 1+\epsilon_{r}.
\end{eqnarray}
where the inequality holds because $\epsilon_{r}\le0.5$ and $\tilde{\epsilon} = \frac{\epsilon_{r}}{4} \le \frac{1}{8}$.
Since we want \eqref{eq:delta_l_epsilon_l_ineq} to hold with probability greater than $1-e^{-t}$ it
remains to require $\abs{\L_{r}^{\text{corank}}}(1+\frac{8}{\epsilon_r})^{\sdim-r}e^{-\frac{C_M m\epsilon_r^2}{32}} \le e^{-t}$.
Using the fact that $(1+\frac{8}{\epsilon_r})\ge \frac{9}{\epsilon_r}$ and some arithmetics we get \eqref{eq:omega_RIP_m_size_cond}
and this completes the proof of the theorem.

We turn now to the proof of Theorem~\ref{thm:analysis_RIP_cond_dependencies}. Its proof is almost identical
to the previous proof but with the difference that instead of $r$, $\L_{r}^{\text{corank}}$ and $\delta_{r}^{\text{corank}}$ we look at $\cosp$, $\L_{\cosp}$ and $\delta_{\cosp}$.
In this case we do not know what is the dimension of the subspace that each cosupport implies. However, we can have a lower bound on it using $\pdim - \cosp$.
Therefore, we use $B^{\pdim- \cosp}$ instead of $B^{\sdim-r}$. This change provides us with a condition similar to \eqref{eq:omega_RIP_m_size_cond} but with $\pdim-\cosp$ in the second coefficient instead of $\sdim-r$. By using some arithmetics, noticing that the size of $\L_{\cosp}$ is $\pdim \choose \cosp$ and using Stirling's formula for upper bounding it we get \eqref{eq:omega_RIP_m_size_cond_dependencies} and this completes the proof.

\section{Proof of Lemma~\ref{lem:AIHT_lemma}}
\label{sec:AIHT_lemma_proof}

{\em Lemma~\ref{lem:AIHT_lemma}:}
Consider the problem $\cal P$ and
apply either AIHT or AHTP with a constant step size $\mu$ satisfying
$\frac{1}{\mu} \ge 1+\delta_{2\cosp - \pdim}$ or an optimal step size.
Then, at the $t$-th iteration, the following holds:
\begin{eqnarray}
\label{eq:AIHT_lemma}
&&\norm{\y-\M\hat\x^t}_2^2 - \norm{\y-\M\hat\x^{t-1}}_2^2
\le C_\ell\left(\norm{\y-\M\x}_2^2 - \norm{\y-\M\hat\x^{t-1}}_2^2\right)\\
&& \nonumber ~~~~    + C_\ell \left(\frac{1}{\mu(1-\deltaB)}-1\right)\norm{\M(\x-\hat\x^{t-1})}_2^2
    + (C_\ell-1)\mu\sigma_\M^2 \norm{\y-\M\hat\x^{t-1}}_2^2.
\end{eqnarray}
For the optimal step size the bound is achieved with the value $\mu = \frac{1}{1+\deltaB}$.

{\em Proof:}
We consider the AIHT algorithm first. We take similar steps to those taken in the proof of Lemma~3 in \cite{Blumensath09Sampling}.
Since $\frac{1}{\mu} \ge 1+\deltaB$, we have, from the $\OM$-RIP property of $\M$,
\[
\norm{\M(\hat\x^t-\hat\x^{t-1})}_2^2 \le \frac{1}{\mu}\norm{\hat\x^t-\hat\x^{t-1}}_2^2.
\]
Thus,
\begin{eqnarray*}
&&\norm{\y-\M\hat\x^t}_2^2 - \norm{\y-\M\hat\x^{t-1}}_2^2
= -2\inn{\M(\hat\x^t-\hat\x^{t-1})}{\y-\M\hat\x^{t-1}} + \norm{\M(\hat\x^t-\hat\x^{t-1})}_2^2 \\
&&~~~~~~~~~~\le -2\inn{\M(\hat\x^t-\hat\x^{t-1})}{\y-\M\hat\x^{t-1}} + \frac{1}{\mu}\norm{\hat\x^t-\hat\x^{t-1}}_2^2 \\
&&~~~~~~~~~~= -2\inn{\hat\x^t-\hat\x^{t-1}}{\M^*(\y-\M\hat\x^{t-1})} + \frac{1}{\mu}\norm{\hat\x^t-\hat\x^{t-1}}_2^2 \\
&&~~~~~~~~~~= -\mu\norm{\M^*(\y-\M\hat\x^{t-1})}_2^2 + \frac{1}{\mu}\norm{\hat\x^t-\hat\x^{t-1}-\mu \M^*(\y-\M\hat\x^{t-1})}_2^2.
\end{eqnarray*}
Note that by definition, $\hat\x^t = \Q_{\hat \CF_\cosp}\left(\hat\x^{t-1}+\mu \M^*(\y-\M\hat\x^{t-1})\right)$. Hence, by the $C_\cosp$-near optimality of the projection, we get
\begin{equation} \label{eq:aihtlemlast}
\norm{\y-\M\hat\x^t}_2^2 - \norm{\y-\M\hat\x^{t-1}}_2^2
\le -\mu\norm{\M^*(\y-\M\hat\x^{t-1})}_2^2 + \frac{C_\ell}{\mu}\norm{\x -\hat\x^{t-1}-\mu \M^*(\y-\M\hat\x^{t-1})}_2^2.
\end{equation}
Now note that
\begin{eqnarray*}
&& \norm{\x -\hat\x^{t-1}-\mu \M^*(\y-\M\hat\x^{t-1})}_2^2 \\
&&~~~~~= \norm{\x -\hat\x^{t-1}}_2^2 - 2\mu\inn{\M(\x-\hat\x^{t-1})}{\y-\M\hat\x^{t-1}} + \mu^2 \norm{\M^*(\y-\M\hat\x^{t-1})}_2^2 \\
&&~~~~~ \le \frac{1}{1-\deltaB} \norm{\M(\x-\hat\x^{t-1})}_2^2 - 2\mu\inn{\M(\x-\hat\x^{t-1})}{\y-\M\hat\x^{t-1}} + \mu^2 \norm{\M^*(\y-\M\hat\x^{t-1})}_2^2 \\
&&~~~~~ = \frac{1}{1-\deltaB} \norm{\M(\x-\hat\x^{t-1})}_2^2 + \mu\left(\norm{\y-\M\x}_2^2 - \norm{\y-\M\hat\x^{t-1}}_2^2 - \norm{\M(\x-\hat\x^{t-1})}_2^2\right) \\
   &&~~~~~~~~~ + \mu^2 \norm{\M^*(\y-\M\hat\x^{t-1})}_2^2.
\end{eqnarray*}
Putting this into \eqref{eq:aihtlemlast}, we obtain the desired result for the AIHT algorithm.

We can check that the same holds true for the AHTP algorithm as follows: suppose that $\hat\x^{t-1}_{\text{\tiny AHTP}}$ is the $(t-1)$-st estimate from the AHTP algorithm. If we now initialize the AIHT algorithm with this estimate and obtain the next estimate $\hat\x^t_{\tilde{\text{\tiny AIHT}}}$, then the inequality of the lemma holds true with $\hat\x^t_{\tilde{\text{\tiny AIHT}}}$ and $\hat\x^{t-1}_{\text{\tiny AHTP}}$ in place of $\hat\x^t$ and $\hat\x^{t-1}$ respectively. On the other hand, from the algorithm description, we know that the $t$-th estimate $\hat\x^t_{\text{\tiny AHTP}}$ of the AHTP satisfies
\[
\norm{\y-\M\hat\x^t_{\text{\tiny AHTP}}}_2^2 \le \norm{\y-\M\hat\x^t_{\tilde{\text{\tiny AIHT}}}}_2^2.
\]
This means that the result holds for the AHTP algorithm as well.

Using a similar argument for the optimal changing step size we note that it selects
the cosupport that minimizes $\norm{\M\x - \M{\hat\x}^{t}}_2^2$. Thus, for AIHT
and AHTP we have that $\norm{\M\x - \M{\hat\x_{\text{\tiny Opt}}}^{t}}_2^2 \le \norm{\M\x - \M{\hat\x}_{\mu}^{t}}_2^2$ for any value of
$\mu$, where $\hat\x_{\text{\tiny Opt}}^{t}$ and $\hat\x_{\mu}^{t}$ are the recovery results of AIHT or AHTP with an optimal changing step-size and a constant step-size $\mu$ respectively.
This yields that any theoretical result for a constant step-size selection with a constant $\mu$ holds true also to the optimal changing-step size selection.
In particular this is true also for $\mu = \frac{1}{1+\deltaB}$. This choice is justified in the proof of Lemma~\ref{lem:AIHT_AHTP_lemma_noisy}.
\hfill $\Box$
\bigskip

\section{Proof of Lemma~\ref{lem:AIHT_AHTP_lemma_noisy}}

{\em Lemma~\ref{lem:AIHT_AHTP_lemma_noisy}:}
Suppose that the same conditions of Theorem~\ref{thm:AIHT_AHTP_theorem_noisy_iter} hold true. If $\norm{\y-\M\hat\x^{t-1}}_2^2 \le \eta^2\norm{\e}_2^2$, then $\norm{\y-\M\hat\x^t}_2^2 \le \eta^2\norm{\e}_2^2$.
Furthermore, if $\norm{\y-\M\hat\x^{t-1}}_2^2 > \eta^2\norm{\e}_2^2$, then
\[
\norm{\y-\M\hat\x^t}_2^2 \le c_4 \norm{\y-\M\hat\x^{t-1}}_2^2
\]
where
\[
c_4 := \left(1+\frac{1}{\eta}\right)^2\left(\frac{1}{\mu(1-\deltaB)} - 1\right)C_\cosp + (C_\ell -1)(\mu\sigma_\M^2 -1) + \frac{C_\cosp}{\eta^2} < 1.
\]

{\em Proof:}
First, suppose that $\norm{\y-\M\hat\x^{t-1}}_2^2 > \eta^2\norm{\e}_2^2$.
From Lemma~\ref{lem:AIHT_lemma}, we have
\begin{eqnarray} \label{eq:aihtsecond}
 \norm{\y-\M\hat\x^t}_2^2
&\le& C_\cosp\norm{\y-\M\x}_2^2 + (C_\cosp - 1) (\mu\sigma_\M^2 - 1) \norm{\y-\M\hat\x^{t-1}}_2^2
 \\ \nonumber &&    + C_\cosp \left(\frac{1}{\mu(1-\deltaB)} - 1\right) \norm{\M(\x-\hat\x^{t-1})}_2^2.
\end{eqnarray}
Remark that all the coefficients in the above are positive because
$1 + \deltaB \le \frac{1}{\mu} \le \sigma_\M^2$ and $C_\cosp \ge 1$.
Since $\y-\M\x = \e$, we note
\[
\norm{\y-\M\x}_2^2 < \frac{1}{\eta^2}\norm{\y-\M\hat\x^{t-1}}_2^2
\]
and, by the triangle inequality,
\[
\norm{\M(\x-\hat\x^{t-1})}_2 \le \norm{\y-\M\x}_2 + \norm{\y-\M\hat\x^{t-1}}_2
< \left( 1 + \frac{1}{\eta} \right) \norm{\y-\M\hat\x^{t-1}}_2.
\]
Therefore, from \eqref{eq:aihtsecond},
\[
\norm{\y-\M\hat\x^t}_2^2 < c_4 \norm{\y-\M\hat\x^{t-1}}_2^2.
\]
This is the second part of the lemma.

Now, suppose that $\norm{\y-\M\hat\x^{t-1}}_2^2 \le \eta^2\norm{\e}_2^2$.
This time we have
\[
\norm{\M(\x-\hat\x^{t-1})}_2 \le \norm{\y-\M\x}_2 + \norm{\y-\M\hat\x^{t-1}}_2
\le (1+\eta) \norm{\e}_2.
\]
Applying this to \eqref{eq:aihtsecond}, we obtain
\begin{align*}
\norm{\y-\M\hat\x^t}_2^2 &\le
C_\cosp \norm{\e}_2^2 + (C_\ell - 1) (\mu\sigma_\M^2 - 1) \eta^2 \norm{\e}_2^2
+ C_\cosp \left(\frac{1}{\mu(1-\deltaB)} - 1\right) (1+\eta)^2 \norm{\e}_2^2 \\
&= \left(C_\cosp + (C_\cosp - 1) (\mu\sigma_M^2 - 1) \eta^2 + C_\cosp \left(\frac{1}{\mu(1-\deltaB)} - 1\right) (1+\eta)^2 \right) \norm{\e}_2^2  =c_4 \eta^2 \norm{\e}_2^2.
\end{align*}
Thus, the proof is complete as soon as we show $c_4 < 1$, or $c_4-1 < 0$.

To see $c_4 -1 < 0$, we first note that it is equivalent to--all the subscripts are dropped from here on for simplicity of notation--
\[
\frac{1}{\mu^2} - \frac{2(1-\delta)}{1+\frac{1}{\eta}} \frac{1}{\mu}
+ \frac{(C-1)\sigma^2 (1-\delta)}{C\left(1+\frac{1}{\eta}\right)^2} < 0,
\]
or
\[
\frac{1}{\mu^2} - 2 (1-\delta) b_1 \frac{1}{\mu} + (1-\delta)^2 b_2 < 0.
\]
Solving this quadratic equation in $\frac{1}{\mu}$, we want
\[
(1-\delta) \left(b_1 - \sqrt{b_1^2 - b_2}\right) < \frac{1}{\mu} <
(1-\delta) \left(b_1 + \sqrt{b_1^2 - b_2}\right).
\]
Such $\mu$ exists only when $\frac{b_2}{b_1^2} < 1$. Furthermore, we have already assumed
$1 + \delta \le \frac{1}{\mu}$ and we know $(1-\delta) \left(b_1 - \sqrt{b_1^2 - b_2}\right) < 1+\delta$, and hence the condition we require is
\[
1 + \delta \le \frac{1}{\mu} <
(1-\delta) \left(b_1 + \sqrt{b_1^2 - b_2}\right),
\]
which is what we desired to prove.

As we have seen in Lemma~\ref{lem:AIHT_lemma}, for changing optimal step-size selection, \eqref{eq:AIHT_AHTP_lemma_noisy}
holds for any value of $\mu$ that satisfies the above conditions.
Thus, in the bound of changing optimal step-size we put a value of $\mu$ that minimizes $c_4$.
This minimization result with $\frac{1}{\mu}=\sqrt{b_2}(1-\deltaB)$.
However, since we need $\frac{1}{\mu} \ge 1+\deltaB$ and have that $\sqrt{b_2}(1-\deltaB) < b_1(1-\deltaB) <1+\deltaB$ we set $\frac{1}{\mu} = 1+\deltaB$ in $c_4$ for the bound in optimal changing step-size case.
\hfill $\Box$
\bigskip

\section{Proof of Lemma~\ref{lem:ACoSaMP_xp_bound}}
\label{sec:ACoSaMP_xp_bound_proof}

{\em Lemma~\ref{lem:ACoSaMP_xp_bound}:}
Consider the problem $\cal P$ and
apply ACoSaMP with $a = \frac{2\cosp - \pdim}{\cosp}$.
For each iteration we have
\begin{eqnarray*}
\norm{\x - \vect{w}}_2 &\le& \frac{1}{\sqrt{1-\deltaD^2}}\norm{\P_{\tilde{\Lambda}^t}(\vect{x} - \vect{w} )}_2+  \frac{\sqrt{1+\deltaC}}{1-\deltaD}\norm{\vect{e}}_2.
\end{eqnarray*}

{\em Proof:}
Since $\vect{w}$ is the minimizer of $\norm{\vect{y} - \vect{M}\vect{v}}^2_2$ with the constraint $\matr{\Omega}_{\tilde{\Lambda}^t}\vect{v} =0$, then
\begin{eqnarray}
\langle \matr{M} \vect{w} - \vect{y}, \matr{M}\vect{u} \rangle =0,
\end{eqnarray}
for any vector $\vect{u}$ such that $\matr{\Omega}_{\tilde{\Lambda}^t}\vect{u} =0 $.
Substituting $\vect{y} = \matr{M}\vect{x} + \vect{e}$ and moving terms from the LHS to the RHS gives
\begin{eqnarray}
\label{eq:xp_x_property}
\langle \vect{w} - \vect{x}, \matr{M}^*\matr{M}\vect{u} \rangle = \langle \vect{e}, \matr{M}\vect{u} \rangle,
\end{eqnarray}
where $\vect{u}$ is a vector satisfying $\matr{\Omega}_{\tilde{\Lambda}^t}\vect{u} =0$.
Turning to look at $\norm{\Q_{\tilde{\Lambda}^t}(\x - \vect{w})}_2^2$ and using \eqref{eq:xp_x_property} with $\vect{u} = \Q_{\tilde{\Lambda}^t}(\x - \vect{w})$, we have
\begin{eqnarray}
\label{eq:Q_xp_x_norm}
&& \hspace{-0.3in} \norm{\Q_{\tilde{\Lambda}^t}(\vect{x} - \vect{w})}_2^2 = \langle \vect{x} - \vect{w}, \Q_{\tilde{\Lambda}^t}(\vect{x} - \vect{w}) \rangle \\
\nonumber && \hspace{-0.3in} = \langle \vect{x} - \vect{w}, (\matr{I} - \matr{M}^*\matr{M})\Q_{\tilde{\Lambda}^t}(\vect{x} - \vect{w}) \rangle - \langle \vect{e}, \matr{M}\Q_{\tilde{\Lambda}^t}(\vect{x} - \vect{w}) \rangle
\\ \nonumber && \hspace{-0.3in} \le \norm{ \vect{x} - \vect{w}}_2 \norm{\Q_{\Lambda \cap \tilde{\Lambda}^t} (\matr{I} - \matr{M}^*\matr{M})\Q_{\tilde{\Lambda}^t}}_2 \norm{\Q_{\tilde{\Lambda}^t}(\vect{x} - \vect{w})}_2  + \norm{ \vect{e}}_2\norm{ \matr{M}\Q_{\tilde{\Lambda}^t}(\vect{x} - \vect{w})}_2
\\ \nonumber && \hspace{-0.3in} \le \deltaD\norm{ \vect{x} - \vect{w}}_2 \norm{\Q_{\tilde{\Lambda}^t}(\vect{x} - \vect{w})}_2    + \norm{ \vect{e}}_2\sqrt{1+\deltaC}\norm{\Q_{\tilde{\Lambda}^t}(\vect{x} - \vect{w})}_2.
\end{eqnarray}
where the first inequality follows from the Cauchy-Schwartz inequality, the projection property that $\Q_{\tilde{\Lambda}^t} = \Q_{\tilde{\Lambda}^t}\Q_{\tilde{\Lambda}^t}$ and the fact that $\vect{x} - \vect{w} = \Q_{\Lambda \cap \tilde{\Lambda}^t  }(\vect{x} - \vect{w})$. The last inequality is due to the $\matr{\Omega}$-RIP properties, Corollary~\ref{cor:omega_RIP_norm_diff} and that according to Table~\ref{tbl:Synthesis_Analysis_parallels} $|{\tilde{\Lambda}^t}| \ge 3\cosp - 2\pdim$ and $|{\Lambda\cap\tilde{\Lambda}^t}| \ge 4\cosp - 3\pdim$.
After simplification of \eqref{eq:Q_xp_x_norm} by $\norm{\Q_{\tilde{\Lambda}^t}(\vect{x} - \vect{w})}_2$ we have
\begin{eqnarray}
\nonumber \norm{\Q_{\tilde{\Lambda}^t}(\vect{x} - \vect{w})}_2 \le \deltaD\norm{\vect{x} - \vect{w}}_2   + \sqrt{1+\deltaC}\norm{\vect{e}}_2.
\end{eqnarray}
Utilizing the last inequality with the fact that $\norm{\vect{x} - \vect{w}}_2^2 = \norm{\P_{\tilde{\Lambda}^t}(\vect{x} - \vect{w})}_2^2+ \norm{\Q_{\tilde{\Lambda}^t}(\vect{x} - \vect{w})}_2^2$ gives
\begin{eqnarray}
&& \norm{\vect{x} - \vect{w}}_2^2  \le \norm{\P_{\tilde{\Lambda}^t}(\vect{x} - \vect{w})}_2^2 + \left(\deltaD\norm{\vect{x} - \vect{w}}_2   + \sqrt{1+\deltaC}\norm{\vect{e}}_2 \right)^2.
\end{eqnarray}
By moving all terms to the LHS we get a quadratic function of $\norm{\x-\vect{w}}_2$.
Thus, $\norm{\x-\vect{w}}_2$ is bounded from above by the larger root of that function;
this with a few simple algebraic steps gives the inequality in \eqref{eq:ACoSaMP_xp_bound}.
\hfill $\Box$ \bigskip

\section{Proof of Lemma~\ref{lem:ACoSaMP_xt_bound1}}
\label{sec:ACoSaMP_xt_bound1_proof}

{\em Lemma~\ref{lem:ACoSaMP_xt_bound1}:}
Consider the problem $\cal P$ and
apply ACoSaMP with $a = \frac{2\cosp - \pdim}{\cosp}$.
For each iteration we have
\begin{eqnarray*}
&& \hspace{-0.5in} \norm{\x - \hat\x^t}_2 \le \rho_1\norm{\P_{\tilde{\Lambda}^t}(\x - \vect{w})}_2 +  \eta_1\norm{\vect{e}}_2,
\end{eqnarray*}
where $\eta_1$ and $\rho_1$ are the same constants as in Theorem~\ref{thm:ACoSaMP_iter_bound}.

{\em Proof:}
We start with the following observation
\begin{eqnarray}
\label{eq:xt_x_diff_norm}
&& \hspace{-0.3in} \norm{\vect{x} -\hat\x^t}_2^2 = \norm{\vect{x} -\vect{w} + \vect{w}-\hat\x^t}_2^2
 = \norm{\vect{x}-\vect{w}}_2^2 + \norm{\hat\x^t - \vect{w}}_2^2 +2(\vect{x} - \vect{w})^*(\vect{w} - \hat\x^t),
\end{eqnarray}
and turn to bound the second and last terms in the RHS.
For the second term, using the fact that $\hat\x^t = \Q_{\hat\CF_\cosp(\vect{w})}\vect{w}$ with \eqref{eq:C_optimal_ineq} gives
\begin{eqnarray}
\label{eq:xt_x_diff_norm_coeff1}
\norm{\hat\x^t - \vect{w}}_2^2 \le C_\cosp\norm{\x -\vect{w}}_2^2.
\end{eqnarray}
For bounding the last term, we look at its absolute value and use \eqref{eq:xp_x_property} with $\vect{u} = \vect{w} - \hat\x^t = \Q_{\tilde{\Lambda}^t}(\vect{w} - \hat\x^t)$. This leads to
\begin{eqnarray*}
&& \hspace{-0.3in}\abs{(\vect{x} - \vect{w})^*(\vect{w} - \hat\x^t)} =   \abs{(\vect{x} - \vect{w})^*(\matr{I} - \matr{M}^*\matr{M})(\vect{w} - \hat\x^t) - \vect{e}^*\matr{M}(\vect{w} - \hat\x^t)}.
\end{eqnarray*}
By using the triangle and Cauchy-Schwartz inequalities with the fact that $\vect{x}-\vect{w} = \Q_{\Lambda \cap \tilde{\Lambda}^t}(\vect{x}-\vect{w})$ and $\vect{w}-\hat\x^t = \Q_{\tilde{\Lambda}^t}(\vect{w}-\hat\x^t)$ we have
\begin{eqnarray}
\label{eq:xt_x_diff_norm_coeff3}
&& \hspace{-0.55in} \abs{(\vect{x} - \vect{w})^*(\vect{w} - \hat\x^t)}
 \le \norm{\vect{x} - \vect{w}}_2\norm{\Q_{\Lambda \cap \tilde{\Lambda}^t}(\matr{I} - \matr{M}^*\matr{M})\Q_{\tilde{\Lambda}^t}}_2\norm{\vect{w} - \hat\x^t}_2
 + \norm{\vect{e}}_2\norm{\matr{M}(\vect{w} - \hat\x^t)}_2
\\ \nonumber && \hspace{0.83in} \le \deltaD\norm{\vect{x}-\vect{w}}_2\norm{\vect{w} - \hat\x^t}_2
 +\sqrt{1+\deltaC}\norm{\vect{e}}_2\norm{\vect{w} - \hat\x^t}_2,
\end{eqnarray}
where the last inequality is due to the $\matr{\Omega}$-RIP definition and Corollary~\ref{cor:omega_RIP_norm_diff}.

By substituting \eqref{eq:xt_x_diff_norm_coeff1} and \eqref{eq:xt_x_diff_norm_coeff3} into \eqref{eq:xt_x_diff_norm} we have
\begin{eqnarray}
\label{eq:xt_x_diff_norm_all}
&& \hspace{-0.3in}\norm{\x -\hat\x^t}_2^2  \le (1+C_\cosp)\norm{\vect{x} -\vect{w}}_2^2
  + 2\deltaD\sqrt{C_\cosp}\norm{\vect{x}-\vect{w}}_2^2 +2\sqrt{1+\deltaC}\sqrt{C_\cosp}\norm{\vect{e}}_2\norm{\x - \vect{w}}_2
\\ \nonumber &&  \le \Big((1 + 2\deltaD\sqrt{C_\cosp} +C_\cosp)\norm{\vect{x} -\vect{w}}_2
 +2\sqrt{(1+\deltaC)C_\cosp}\norm{\vect{e}}_2\Big) \norm{\vect{x} - \vect{w}}_2 \\  \nonumber
&&  \le \frac{1 + 2\deltaD\sqrt{C_\cosp} +C_\cosp}{1-\deltaD^2}\norm{\P_{\tilde{\Lambda}^t}(\x -\vect{w} )}_2^2 \\   \nonumber &&+   \frac{2\sqrt{1+\deltaC}(1  +(1+\deltaD)\sqrt{C_\cosp}+C_\cosp)}{(1-\deltaD)\sqrt{1-\deltaD^2}}\norm{\P_{\tilde{\Lambda}^t}(\x-\vect{w} )}_2\norm{\vect{e}}_2  + \frac{(1+\deltaC)(1  +2\sqrt{C_\cosp}+C_\cosp)}{(1-\deltaD)^2}\norm{\vect{e}}_2^2
\\ \nonumber && \le \Bigg(\frac{\sqrt{1+2\deltaD\sqrt{C_\cosp} + C_\cosp}}{\sqrt{1-\deltaD^2}}
\norm{\P_{\tilde{\Lambda}^t}(\x - \vect{w} )}_2 +\frac{\sqrt{\frac{2+C_\cosp}{1+C_\cosp}+2\sqrt{C_\cosp}+C_\cosp}\sqrt{1+\deltaC}}{1-\deltaD}\norm{\vect{e}}_2 \Bigg)^2,
\end{eqnarray}
where for the second inequality we use the fact that $\deltaD \le 1$ combined with the inequality of Lemma~\ref{lem:ACoSaMP_xp_bound},
and for the last inequality we use the fact that $(1  +(1+\deltaD)\sqrt{C_\cosp}+C_\cosp)^2 \le (1 + 2\deltaD\sqrt{C_\cosp} +C_\cosp)(\frac{2+C_\cosp}{1+C_\cosp}+2\sqrt{C_\cosp}+C_\cosp)$ together with a few algebraic steps.
Taking square-root on both sides of \eqref{eq:xt_x_diff_norm_all} provides the desired result.
\hfill $\Box$ \bigskip

\section{Proof of Lemma~\ref{lem:ACoSaMP_Pxp_bound}}
\label{sec:ACoSaMP_Pxp_bound_proof}

{\em Lemma~~\ref{lem:ACoSaMP_Pxp_bound}:}
Consider the problem $\cal P$ and
apply ACoSaMP with $a = \frac{2\cosp - \pdim}{\cosp}$. if
\begin{eqnarray*}
C_{2\cosp-p} < \frac{\sigma_{\matr{M}}^2(1+\gamma)^2}{\sigma_{\matr{M}}^2(1+\gamma)^2-1},
\end{eqnarray*}
then there exists $\tilde{\delta}_{\text{ \tiny ACoSaMP}}(C_{2\cosp - \pdim}, \sigma_{\matr{M}}^2, \gamma) >0$ such that
for any $\deltaB < \tilde{\delta}_{\text{ \tiny ACoSaMP}}(C_{2\cosp - \pdim}, \sigma_{\matr{M}}^2, \gamma)$
\begin{eqnarray*}
&& \norm{\P_{\tilde{\Lambda}^t}(\x - \vect{w})}_2 \le
    \eta_2\norm{\vect{e}}_2
 + \rho_2\norm{\x - \hat\x^{t-1}}_2.
\end{eqnarray*}
The constants $\eta_2$ and $\rho_2$ are as defined in Theorem~\ref{thm:ACoSaMP_iter_bound}.

In the proof of the lemma we use the following Proposition.

{\em Proposition~~E.1:}
For any two given vectors $\vect{x}_1$, $\vect{x}_2$ and any constant $c>0$ it holds that
\begin{eqnarray}
\label{eq:norm2_ineq}
\norm{\vect{x}_1+\vect{x}_2}_2^2 \le (1+c)\norm{\vect{x}_1}_2^2 + \left(1+\frac{1}{c}\right)\norm{\vect{x}_2}^2.
\end{eqnarray}

The proof of the proposition is immediate using the inequality of arithmetic and geometric means.
We turn to the proof of the lemma.

{\em Proof:}
Looking at the step of finding new cosupport elements one can observe that $\Q_{\Lambda_{\Delta}}$
is a near optimal projection for $\matr{M}^*\vect{y}^{t-1}_{\text{resid}} = \matr{M}^*(\vect{y} - \matr{M}\hat\x^{t-1})$
with a constant $C_{2\cosp-p}$. The fact that $\abs{{\hat\Lambda^{t-1} \cap \Lambda}} \ge 2\cosp -p$ combined
with \eqref{eq:C_optimal_ineq} gives
\begin{eqnarray*}
\norm{(\matr{I} - \Q_{{\Lambda}_{\Delta}})\matr{M}^*(\vect{y} - \matr{M}\hat\x^{t-1})}_2^2 \le C_{2\cosp -\pdim}\norm{(\matr{I} - \Q_{\hat\Lambda^{t-1} \cap \Lambda})\matr{M}^*(\vect{y} - \matr{M}\hat\x^{t-1})}_2^2.
\end{eqnarray*}
Using simple projection properties and the fact that $\tilde{\Lambda}^t \subseteq \Lambda_{\Delta}$ with $\vect{z} = \matr{M}^*(\vect{y} - \matr{M}\hat\x^{t-1})$
we have
\begin{eqnarray}
\label{eq:ACoSaMP_QMy_Mx_ineq}
&& \hspace{-0.3in} \norm{\Q_{\tilde{\Lambda}^t}\vect{z}}_2^2 \ge \norm{\Q_{{\Lambda}_\Delta}\vect{z}}_2^2
=\norm{\vect{z}}_2^2 - \norm{(\matr{I}-\Q_{{\Lambda}_\Delta})\vect{z}}_2^2
\ge \norm{\vect{z}}_2^2 - C_{2\cosp-p}\norm{(\matr{I}-\Q_{\hat\Lambda^{t-1} \cap \Lambda})\vect{z}}_2^2
 \\ \nonumber && =\norm{\vect{z}}_2^2 - C_{2\cosp-p}\left(\norm{\vect{z}}_2^2 - \norm{\Q_{\hat\Lambda^{t-1} \cap \Lambda}\vect{z}}_2^2\right)
= C_{2\cosp-p}\norm{\Q_{\hat\Lambda^{t-1} \cap \Lambda}\vect{z}}_2^2- (C_{2\cosp-p}-1)\norm{\vect{z}}_2^2.
\end{eqnarray}

We turn to bound the LHS of \eqref{eq:ACoSaMP_QMy_Mx_ineq} from above. Noticing that
$\y = \M\x+\e$ and using \eqref{eq:norm2_ineq} with a constant $\gamma_1>0$ gives
\begin{eqnarray}
\label{eq:ACoSaMP_QMy_Mx_ineq_lhs_step1}
\norm{\Q_{\tilde{\Lambda}^t}\matr{M}^*(\vect{y} - \matr{M}\hat\x^{t-1})}_2^2 \le
\left(1+\frac{1}{\gamma_1}\right)\norm{\Q_{\tilde{\Lambda}^t}\matr{M}^*\vect{e}}_2^2 + (1+\gamma_1)\norm{\Q_{\tilde{\Lambda}^t}\matr{M}^*\matr{M}(\vect{x} - \hat\x^{t-1})}_2^2.
\end{eqnarray}
Using \eqref{eq:norm2_ineq} again, now with a constant $\alpha>0$, we have
\begin{eqnarray}
\label{eq:ACoSaMP_QMy_Mx_ineq_lhs_step2}
&&\hspace{-0.32in}\norm{\Q_{\tilde{\Lambda}^t}\matr{M}^*\matr{M}(\vect{x} - \hat\x^{t-1})}_2^2
 \le  (1+\alpha)\norm{\Q_{\tilde{\Lambda}^t}(\vect{x} - \hat\x^{t-1})}_2^2
+\left(1+\frac{1}{\alpha}\right)\norm{\Q_{\tilde{\Lambda}^t}(\matr{I}-\matr{M}^*\matr{M})(\vect{x} - \hat\x^{t-1})}_2^2
\\ \nonumber && \le (1+\alpha)\norm{\vect{x} - \hat\x^{t-1}}_2^2 - (1+\alpha)\norm{\P_{\tilde{\Lambda}^t}(\vect{x} - \hat\x^{t-1})}_2^2
+\left(1+\frac{1}{\alpha}\right)\norm{\Q_{\tilde{\Lambda}^t}(\matr{I}-\matr{M}^*\matr{M})(\vect{x} - \hat\x^{t-1})}_2^2.
\end{eqnarray}
Putting \eqref{eq:ACoSaMP_QMy_Mx_ineq_lhs_step2} into \eqref{eq:ACoSaMP_QMy_Mx_ineq_lhs_step1}
and using \eqref{eq:omega_RIP_norm_diff} and Corollary~\ref{cor:MQ_RIP_norm} gives
\begin{eqnarray}
\label{eq:ACoSaMP_QMy_Mx_ineq_lhs}
&& \hspace{-0.5in} \norm{\Q_{\tilde{\Lambda}^t}\matr{M}^*(\vect{y} - \matr{M}\hat\x^{t-1})}_2^2
 \le  \frac{(1+\gamma_1)(1+\delta_{3\cosp-2p})}{\gamma_1}\norm{\vect{e}}_2^2
 - (1+\alpha)(1+\gamma_1)\norm{\P_{\tilde{\Lambda}^t}(\vect{x} - \hat\x^{t-1})}_2^2
 \\ \nonumber && ~~~~~~~~~~~~~~~~~~~~~~~ +\left(1+\alpha +\delta_{4\cosp-3p} + \frac{\delta_{4\cosp-3p}}{\alpha}\right)(1+\gamma_1)\norm{\vect{x} - \hat\x^{t-1}}_2^2.
\end{eqnarray}

We continue with bounding the RHS of \eqref{eq:ACoSaMP_QMy_Mx_ineq} from below.
For the first element of the RHS we use an altered version of \eqref{eq:norm2_ineq} with a constant $\gamma_2>0$ and have
\begin{eqnarray}
\label{eq:ACoSaMP_QMy_Mx_ineq_rhs1_step1}
 \norm{\Q_{{\hat\Lambda^{t-1} \cap \Lambda}}\matr{M}^*(\vect{y} - \matr{M}\hat\x^{t-1})}_2^2
\ge \frac{1}{1+\gamma_2}\norm{\Q_{\hat\Lambda^{t-1} \cap \Lambda}\matr{M}^*\matr{M}(\vect{x} - \hat\x^{t-1})}_2^2
-\frac{1}{\gamma_2}\norm{\Q_{\hat\Lambda^{t-1} \cap \Lambda}\matr{M}^*\vect{e}}_2^2.
\end{eqnarray}
Using the altered form again, for the first element in the RHS of \eqref{eq:ACoSaMP_QMy_Mx_ineq_rhs1_step1}, with a constant $\beta >0$ gives
\begin{eqnarray}
\label{eq:ACoSaMP_QMy_Mx_ineq_rhs1_step2}
\norm{\Q_{\hat\Lambda^{t-1} \cap \Lambda}\matr{M}^*\matr{M}(\vect{x} - \hat\x^{t-1})}_2^2
 \ge \frac{1}{1+\beta}\norm{\vect{x} - \hat\x^{t-1}}_2^2
 -\frac{1}{\beta}\norm{\Q_{\hat\Lambda^{t-1} \cap \Lambda}(\matr{M}^*\matr{M} - \matr{I})(\vect{x} - \hat\x^{t-1})}_2^2.
\end{eqnarray}
Putting \eqref{eq:ACoSaMP_QMy_Mx_ineq_rhs1_step2} in \eqref{eq:ACoSaMP_QMy_Mx_ineq_rhs1_step1} and using
the RIP properties and \eqref{eq:omega_RIP_norm_diff} provide
\begin{eqnarray}
\label{eq:ACoSaMP_QMy_Mx_ineq_rhs1}
 \norm{\Q_{{\hat\Lambda^{t-1} \cap \Lambda}}\matr{M}^*(\vect{y} - \matr{M}\hat\x^{t-1})}_2^2
 \ge \left(\frac{1}{1+\beta}-\frac{\delta_{2\cosp-p}}{\beta}\right)\frac{1}{1+\gamma_2}\norm{\vect{x} - \hat\x^{t-1}}_2^2
-\frac{(1+\delta_{2\cosp-p})}{\gamma_2}\norm{\vect{e}}_2^2.
\end{eqnarray}
Using \eqref{eq:norm2_ineq}, with a constant $\gamma_3 >0$,  \eqref{eq:omega_RIP}, and some basic algebraic steps we have for the second element in the RHS of \eqref{eq:ACoSaMP_QMy_Mx_ineq}
\begin{eqnarray}
\label{eq:ACoSaMP_QMy_Mx_ineq_rhs2}
 \hspace{-0.3in} \norm{\matr{M}^*(\vect{y} - \matr{M}\hat\x^{t-1})}_2^2
 &\le& (1+\gamma_3)\norm{\matr{M}^*\matr{M}(\vect{x} - \hat\x^{t-1})}_2^2 + \left(1+\frac{1}{\gamma_3}\right)\norm{\matr{M}^*\vect{e}}_2^2
\\ \nonumber &\le& (1+\gamma_3)(1+\delta_{2\cosp-p})\sigma_{\matr{M}}^2\norm{(\vect{x} - \hat\x^{t-1})}_2^2 + \left(1+\frac{1}{\gamma_3}\right)\sigma_{\matr{M}}^2\norm{\vect{e}}_2^2.
\end{eqnarray}

By combining \eqref{eq:ACoSaMP_QMy_Mx_ineq_lhs}, \eqref{eq:ACoSaMP_QMy_Mx_ineq_rhs1} and \eqref{eq:ACoSaMP_QMy_Mx_ineq_rhs2} with
\eqref{eq:ACoSaMP_QMy_Mx_ineq} we have
\begin{eqnarray}
&& \hspace{-0.15in}(1+\alpha)(1+\gamma_1)\norm{\P_{\tilde{\Lambda}^t}(\vect{x} - \hat\x^{t-1})}_2^2
 \le  \frac{(1+\gamma_1)(1+\delta_{3\cosp-2p})}{\gamma_1}\norm{\vect{e}}_2^2
 +C_{2\cosp-p}\frac{(1+\delta_{2\cosp-p})}{\gamma_2}\norm{\vect{e}}_2^2
\\ \nonumber &&  + (C_{2\cosp-p}-1)\left(1+\frac{1}{\gamma_3}\right)\sigma_{\matr{M}}^2\norm{\vect{e}}_2^2
 + \left(1+\alpha +\delta_{4\cosp-3p} + \frac{\delta_{4\cosp-3p}}{\alpha}\right)(1+\gamma_1)\norm{\vect{x} - \hat\x^{t-1}}_2^2
\\ \nonumber &&
+ (C_{2\cosp-p}-1)(1+\gamma_3)(1+\delta_{2\cosp-p})\sigma_{\matr{M}}^2\norm{(\vect{x} - \hat\x^{t-1})}_2^2
- C_{2\cosp-p}\left(\frac{1}{1+\beta}-\frac{\delta_{2\cosp-p}}{\beta}\right)\frac{1}{1+\gamma_2}\norm{\vect{x} - \hat\x^{t-1}}_2^2.
\end{eqnarray}
Dividing both sides by $(1+\alpha)(1+\gamma_1)$ and gathering coefficients give
\begin{eqnarray}
&& \hspace{-0.3in}\norm{\P_{\tilde{\Lambda}^t}(\vect{x} - \hat\x^{t-1})}_2^2 \le
   \bigg(\frac{1+\delta_{3\cosp-2p}}{\gamma_1(1+\alpha)}
 +\frac{(1+\delta_{2\cosp-p})C_{2\cosp-p}}{\gamma_2(1+\alpha)(1+\gamma_1)}
 + \frac{(C_{2\cosp-p}-1)(1+\gamma_3)\sigma_{\matr{M}}^2}{(1+\alpha)(1+\gamma_1)\gamma_3}\bigg)\norm{\vect{e}}_2^2
 \\ \nonumber && \hspace{-0.3in} ~~~~~~~~~~~~~~~~~~~~~~~~~~
 + \bigg(1 +\frac{\delta_{4\cosp-3p}}{\alpha} + \frac{(C_{2\cosp-p}-1)(1+\gamma_3)(1+\delta_{2\cosp-p})\sigma_{\matr{M}}^2}{(1+\alpha)(1+\gamma_1)}
 \\ \nonumber && ~~~~~~~~~~~~~~~~~~~~~~~~~~ - \frac{C_{2\cosp-p}}{(1+\alpha)(1+\gamma_1)(1+\gamma_2)}\left(\frac{1}{1+\beta}-\frac{\delta_{2\cosp-p}}{\beta}\right)\bigg)\norm{\vect{x} - \hat\x^{t-1}}_2^2.
\end{eqnarray}
The smaller the coefficient of $\norm{\vect{x} - \hat\x^{t-1}}_2^2$, the better convergence guarantee we obtain. Thus, we choose
$\beta = \frac{\sqrt{\delta_{2\cosp-p}}}{1- \sqrt{\delta_{2\cosp-p}}}$ and
$ \small {\alpha = \frac{\sqrt{\delta_{4\cosp-3p}}}{\sqrt{\frac{C_{2\cosp-p}}{(1+\gamma_1)(1+\gamma_2)}\left(1-\sqrt{\delta_{2\cosp-p}}\right)^2 -\frac{(C_{2\cosp-p}-1)(1+\gamma_3)(1+\delta_{2\cosp-p})\sigma_{\matr{M}}^2}{1+\gamma_1}} - \sqrt{\delta_{4\cosp-3p}}}}$ so that the coefficient is minimized.
The values of $\gamma_1, \gamma_2, \gamma_3$ provide a tradeoff between the convergence rate and the size of the noise coefficient.
For smaller values we get better convergence rate but higher amplification of the noise.
We make no optimization on their values and choose them to be $\gamma_1 = \gamma_2 = \gamma_3 = \gamma$
for an appropriate $\gamma>0$. Thus, the above yields
\begin{eqnarray}
\label{eq:ACoSaMP_Pxp_bound_square}
&& \hspace{-0.3in}\norm{\P_{\tilde{\Lambda}^t}(\vect{x} - \hat\x^{t-1})}_2^2 \le
  \Bigg(\frac{1+\delta_{3\cosp-2p}}{\gamma(1+\alpha)}
 +\frac{(1+\delta_{2\cosp-p})C_{2\cosp-p}}{\gamma(1+\alpha)(1+\gamma)}
+ \frac{(C_{2\cosp-p}-1)(1+\gamma)\sigma_{\matr{M}}^2}{(1+\alpha)(1+\gamma)\gamma}\bigg)\norm{\vect{e}}_2^2
  \\ \nonumber &&
 + \left(1   -\left(\sqrt{\delta_{4\cosp-3p}}  -
 \sqrt{\frac{C_{2\cosp-p}}{(1+\gamma)^2}\left(1-\sqrt{\delta_{2\cosp-p}}\right)^2 -(C_{2\cosp-p}-1)(1+\delta_{2\cosp-p})\sigma_{\matr{M}}^2}
 \right)^2\right)\norm{\vect{x} - \hat\x^{t-1}}_2^2.
\end{eqnarray}
Since $\P_{\tilde{\Lambda}^t}\vect{w} = \P_{\tilde{\Lambda}^t}\hat\x^{t-1} =0$ the above inequality holds also for
$\norm{\P_{\tilde{\Lambda}^t}(\vect{x} - \hat\x^{t-1})}_2^2$.
Inequality \eqref{eq:ACoSaMP_Pxp_bound} follows since the right-hand side of \eqref{eq:ACoSaMP_Pxp_bound_square} is smaller than the square of the right-hand side of \eqref{eq:ACoSaMP_Pxp_bound}.

Before ending the proof, we notice that $\rho_2$, the coefficient of $\norm{\vect{x} - \hat\x^{t-1}}_2^2$ is defined only when
\begin{eqnarray}
\label{eq:rho_2_delta_ineq}
(C_{2\cosp-p}-1)(1+\deltaB)\sigma_{\matr{M}}^2 \le \frac{C_{2\cosp-p}}{(1+\gamma)^2}\left(1-\sqrt{\deltaB}\right)^2.
\end{eqnarray}
First we notice that since $1+\delta_{2\cosp-p} \ge \left(1-\sqrt{\delta_{2\cosp-p}}\right)^2$ a necessary condition for \eqref{eq:rho_2_delta_ineq} to hold is
$(C_{2\cosp-p}-1)\sigma_{\matr{M}}^2 < \frac{C_{2\cosp-p}}{(1+\gamma)^2}$ which is equivalent to \eqref{eq:C_2lp_cond}.
By moving the terms in the RHS to the LHS we get a quadratic function of $\sqrt{\deltaB}$.
The condition in \eqref{eq:C_2lp_cond} guarantees that its constant term is smaller than zero
and thus there exists a positive $\deltaB$ for which the function is smaller than zero.
Therefore, for any $\deltaB < \tilde{\delta}_{\text{ \tiny ACoSaMP}}(C_{2\cosp - \pdim}, \sigma_{\matr{M}}^2, \gamma)$ \eqref{eq:rho_2_delta_ineq} holds, where
$\tilde{\delta}_{\text{ \tiny ACoSaMP}}(C_{2\cosp - \pdim}, \sigma_{\matr{M}}^2, \gamma) >0$
is the square of the positive solution of the quadratic function.

\hfill $\Box$ \bigskip

\section*{Acknowledgment}
The authors would like to thank Jalal Fadili for fruitful discussion,
and the unknown reviewers for the important remarks that helped to improved the shape of the paper.
Without both of them, the examples of the optimal projections would not have appeared in the paper.
This research was supported by New York Metropolitan Research Fund.
R. Giryes is grateful to the Azrieli Foundation for the award of
an Azrieli Fellowship.
This work was supported in part by the EU FP7, SMALL project under FET-Open grant number 225913, and EPSRC grants EP/J015180/1 and EP/F039697/1. 
R. Gribonval acknowledges the support of the European Research Council, PLEASE project, under grant ERC-StG- 2011-277906.
MED acknowledges support of his position from the Scottish Funding Council and their support of the Joint Research Institute with the Heriot-Watt University as a component part of the Edinburgh Research Partnership.

\bibliographystyle{elsarticle-num}
\bibliography{AnalysisGreedyLike}

\begin{thebibliography}{10}
\expandafter\ifx\csname url\endcsname\relax
  \def\url#1{\texttt{#1}}\fi
\expandafter\ifx\csname urlprefix\endcsname\relax\def\urlprefix{URL }\fi
\expandafter\ifx\csname href\endcsname\relax
  \def\href#1#2{#2} \def\path#1{#1}\fi

\bibitem{Lu08Theory}
Y.~Lu, M.~Do, A theory for sampling signals from a union of subspaces, IEEE
  Trans. Signal Process. 56~(6) (2008) 2334 --2345.

\bibitem{Donoho02Optimally}
D.~Donoho, M.~Elad, Optimally sparse representation in general (nonorthogonal)
  dictionaries via $\ell^1$ minimization, Proc. Nat. Aca. Sci. 100~(5) (2003)
  2197--2202.

\bibitem{Gribonval03spars}
R.~Gribonval, M.~Nielsen, Sparse representations in unions of bases, IEEE
  Trans. Inf. Theory. 49~(12) (2003) 3320--3325.

\bibitem{Davis97Adaptive}
G.~Davis, S.~Mallat, M.~Avellaneda, Adaptive greedy approximations,
  Constructive Approximation 13 (1997) 57--98.

\bibitem{Candes06Near}
E.~J. Cand\`{e}s, T.~Tao, Near-optimal signal recovery from random projections:
  Universal encoding strategies?, IEEE Trans. Inf. Theory. 52~(12) (2006) 5406
  --5425.

\bibitem{foucart10Sparse}
S.~Foucart, Sparse recovery algorithms: sufficient conditions in terms of
  restricted isometry constants, in: Approximation Theory XIII, Springer
  Proceedings in Mathematics, 2010, pp. 65--77.

\bibitem{Mo11New}
Q.~Mo, S.~Li, New bounds on the restricted isometry constant, Applied and
  Computational Harmonic Analysis 31~(3) (2011) 460 -- 468.

\bibitem{Rauhut08Compressed}
H.~Rauhut, K.~Schnass, P.~Vandergheynst, Compressed sensing and redundant
  dictionaries, IEEE Trans. Inf. Theory. 54~(5) (2008) 2210 --2219.

\bibitem{Pati93OMP}
Y.~Pati, R.~Rezaiifar, P.~Krishnaprasad, Orthonormal matching pursuit~:
  recursive function approximation with applications to wavelet decomposition,
  in: Proceedings of the $27^{th}$ Annual Asilomar Conf. on Signals, Systems
  and Computers, 1993.

\bibitem{MallatZhang93}
S.~Mallat, Z.~Zhang, Matching pursuits with time-frequency dictionaries, IEEE
  Trans. Signal Process. 41 (1993) 3397--3415.

\bibitem{Needell09CoSaMP}
D.~Needell, J.~Tropp, {CoSaMP}: Iterative signal recovery from incomplete and
  inaccurate samples, Applied and Computational Harmonic Analysis 26~(3) (2009)
  301 -- 321.

\bibitem{Dai09Subspace}
W.~Dai, O.~Milenkovic, Subspace pursuit for compressive sensing signal
  reconstruction, IEEE Trans. Inf. Theory. 55~(5) (2009) 2230 --2249.

\bibitem{Blumensath09Iterative}
T.~Blumensath, M.~Davies, Iterative hard thresholding for compressed sensing,
  Applied and Computational Harmonic Analysis 27~(3) (2009) 265 -- 274.

\bibitem{Foucart11Hard}
S.~Foucart, Hard thresholding pursuit: an algorithm for compressive sensing,
  SIAM J. Numer. Anal. 49~(6) (2011) 2543--2563.

\bibitem{Giryes12RIP}
R.~Giryes, M.~Elad, {RIP}-based near-oracle performance guarantees for {SP},
  {CoSaMP}, and {IHT}, IEEE Trans. Signal Process. 60~(3) (2012) 1465--1468.

\bibitem{Garg09Gradient}
R.~Garg, R.~Khandekar, Gradient descent with sparsification: an iterative
  algorithm for sparse recovery with restricted isometry property, in:
  Proceedings of the 26th Annual International Conference on Machine Learning,
  ICML '09, ACM, New York, NY, USA, 2009, pp. 337--344.

\bibitem{Zhang11Sparse}
T.~Zhang, Sparse recovery with orthogonal matching pursuit under {RIP}, IEEE
  Trans. Inf. Theory. 57~(9) (2011) 6215 --6221.

\bibitem{Nam12Cosparse}
S.~Nam, M.~Davies, M.~Elad, R.~Gribonval, The cosparse analysis model and
  algorithms, Applied and Computational Harmonic Analysis.

\bibitem{Nam11CosparseConf}
S.~Nam, M.~Davies, M.~Elad, R.~Gribonval, Cosparse analysis modeling -
  uniqueness and algorithms, in: IEEE International Conference on Acoustics,
  Speech and Signal Processing (ICASSP), 2011.

\bibitem{elad07Analysis}
M.~Elad, P.~Milanfar, R.~Rubinstein, Analysis versus synthesis in signal
  priors, Inverse Problems 23~(3) (2007) 947--968.

\bibitem{Candes11Compressed}
E.~J. Cand\`{e}s, Y.~C. Eldar, D.~Needell, P.~Randall, Compressed sensing with
  coherent and redundant dictionaries, Applied and Computational Harmonic
  Analysis 31~(1) (2011) 59 -- 73.

\bibitem{Vaiter12Robust}
S.~Vaiter, G.~Peyr\`{e}, C.~Dossal, J.~Fadili, Robust sparse analysis
  regularization, submitted to IEEE Trans. on Information Theory.

\bibitem{Nam11GAPN}
S.~Nam, M.~Davies, M.~Elad, R.~Gribonval, Cosparse analysis modeling, in: 9th
  International Conference on SamplingTheory and Applications (sampta-2011),
  Singapore, 2011.

\bibitem{IRLS10Deubechies}
I.~Daubechies, R.~DeVore, M.~Fornasier, C.~S. G\"{u}nt\"{u}rk, Iteratively
  reweighted least squares minimization for sparse recovery, Communications on
  Pure and Applied Mathematics 63~(1) (2010) 1--38.

\bibitem{Rubinstein12Cosparse}
R.~Rubinstein, T.~Peleg, M.~Elad, {Analysis K-SVD: A dictionary learning
  algorithm for the analysis sparse model}, submitted to IEEE Trans. on Signal
  Processing.

\bibitem{Peleg12Performance}
T.~Peleg, M.~Elad, Performance guarantees of the thresholding algorithm for the
  {Co-Sparse} analysis model, submitted to IEEE Trans. on Information Theory.

\bibitem{Giryes11Iterative}
R.~Giryes, S.~Nam, R.~Gribonval, M.~E. Davies, Iterative cosparse projection
  algorithms for the recovery of cosparse vectors, in: {The 19th European
  Signal Processing Conference (EUSIPCO-2011)}, Barcelona, Spain, 2011.

\bibitem{Giryes12CoSaMP}
R.~Giryes, M.~Elad, {CoSaMP} and {SP} for the cosparse analysis model, in: {The
  20th European Signal Processing Conference (EUSIPCO-2012)}, Bucharest,
  Romania, 2012.

\bibitem{Blumensath09Sampling}
T.~Blumensath, M.~Davies, Sampling theorems for signals from the union of
  finite-dimensional linear subspaces, IEEE Trans. Inf. Theory. 55~(4) (2009)
  1872 --1882.

\bibitem{Mendelson08Uniform}
S.~Mendelson, A.~Pajor, N.~Tomczak-Jaegermann, Uniform uncertainty principle
  for bernoulli and subgaussian ensembles, Constructive Approximation 28 (2008)
  277--289.

\bibitem{Krahmer11New}
F.~Krahmer, R.~Ward, New and improved {J}ohnson-{L}indenstrauss embeddings via
  the restricted isometry property, SIAM J. Math. Analysis 43~(3) (2011)
  1269--1281.

\bibitem{Baraniuk08Simple}
R.~Baraniuk, M.~Davenport, R.~DeVore, M.~Wakin, A simple proof of the
  restricted isometry property for random matrices, Constructive Approximation
  28~(3) (2008) 253--263.

\bibitem{Gribonval12Projection}
R.~{G}ribonval, M.~E. Pfetsch, A.~M. Tillmann, Projection onto the k-cosparse
  set is {NP}-hard, Unpuhlished draft, 2012.

\bibitem{Han04Optimal}
T.~Han, S.~Kay, T.~Huang, Optimal segmentation of signals and its application
  to image denoising and boundary feature extraction, in: International
  Conference on Image Processing, 2004. ICIP '04., Vol.~4, 2004, pp. 2693 --
  2696 Vol. 4.

\bibitem{Tibshirani05Sparsity}
R.~Tibshirani, M.~Saunders, S.~Rosset, J.~Zhu, K.~Knight, Sparsity and
  smoothness via the fused {L}asso, Journal of the Royal Statistical Society:
  Series B (Statistical Methodology) 67~(1) (2005) 91--108.

\bibitem{Kyrillidis11Recipes}
A.~Kyrillidis, V.~Cevher, Recipes on hard thresholding methods, in:
  Computational Advances in Multi-Sensor Adaptive Processing (CAMSAP), 2011 4th
  IEEE International Workshop on, 2011, pp. 353 --356.

\bibitem{Blumensath12Accelerated}
T.~Blumensath, Accelerated iterative hard thresholding, Signal Processing
  92~(3) (2012) 752 -- 756.

\bibitem{Rudelson10Non}
M.~Rudelson, R.~Vershynin, {Non-asymptotic theory of random matrices: extreme
  singular values}, in: International Congress of Mathematicans, 2010.

\bibitem{Candes08TheRIP}
E.~J. Cand\`{e}s, The restricted isometry property and its implications for
  compressed sensing, Comptes-rendus de l'Acad\'{e}mie des Sciences, Paris,
  Series I 346~(9–10) (2008) 589 -- 592.

\bibitem{Needell12Stable}
D.~Needell, R.~Ward, Stable image reconstruction using total variation
  minimization, to appear in SIAM J. Imaging Sciences.

\bibitem{Donoho09countingfaces}
D.~L. Donoho, J.~Tanner, Counting faces of randomly-projected polytopes when
  the projection radically lowers dimension, J. of the AMS (2009) 1--53.

\end{thebibliography}

\end{document}